\documentclass{amsart}
\title{Barycentric Subdivision and Isomorphisms of Groupoids}
\author{Jasha Sommer-Simpson}

\date{\today}
\usepackage{adjustbox}    
\usepackage{amssymb}      
\usepackage{stmaryrd}     
\usepackage{array}        
\usepackage{calc}         
\usepackage{caption}      
\usepackage{color}        
\usepackage{csquotes}     
\usepackage{enumitem}     
\usepackage{float}        
\usepackage{framed}       
\usepackage{graphicx}     
\usepackage{hyperref}     
\usepackage{ifthen}       
\usepackage{longtable}    
\usepackage{mathrsfs}     
\usepackage{mathtools}    
\usepackage{mdframed}     
\usepackage{multicol}     
\usepackage{picinpar}     
\usepackage{placeins}     
\usepackage{tikz}         
\usepackage{setspace}     
\usepackage{soul}         
\usepackage{wrapfig}      
\usepackage{xparse}       
\usepackage[all]{xy}			
\entrymodifiers={+!!<0pt,\fontdimen22\textfont2>}

\tikzset{node distance=1cm, auto}
\usetikzlibrary{positioning}             
\usetikzlibrary{calc}

\newsavebox{\overlongequation}
\newenvironment{dontbotheriftheequationisoverlong}
 {\begin{displaymath}\begin{lrbox}{\overlongequation}$\displaystyle}
 {$\end{lrbox}\makebox[0pt]{\usebox{\overlongequation}}\end{displaymath}}

        \makeatletter
        \renewenvironment{abstract}{%
          \ifx\maketitle\relax
            \ClassWarning{\@classname}{Abstract should precede
              \protect\maketitle\space in AMS document classes; reported}%
          \fi
          \global\setbox\abstractbox=\vtop \bgroup
            \normalfont\Small
            \list{}{\labelwidth\z@
              \leftmargin2pc \rightmargin\leftmargin
              \listparindent\normalparindent \itemindent\z@
              \parsep\z@ \@plus\p@
              
            }%
            \item[\hskip\labelsep\scshape\abstractname.]%
        }{%
          \endlist\egroup
          \ifx\@setabstract\relax \@setabstracta \fi
        }
        \makeatother

\theoremstyle{plain}
\newtheorem{thm}{Theorem}[section]
\newtheorem{cor}[thm]{Corollary}
\newtheorem{prop}[thm]{Proposition}
\newtheorem{lem}[thm]{Lemma}
\theoremstyle{definition}
\newtheorem{defn}[thm]{Definition}
\newtheorem*{defn*}{Definition}
\newtheorem{defns}[thm]{Definitions}
\newtheorem*{defns*}{Definitions}
\newtheorem{con}[thm]{Construction}
\newtheorem*{exmp*}{Example}
\newtheorem{exmp}[thm]{Example}
\newtheorem{notn}[thm]{Notation}
\newtheorem*{notn*}{Notation}
\theoremstyle{remark}
\newtheorem{rem}[thm]{Remark}
\makeatletter
\let\c@equation\c@thm
\makeatother
\numberwithin{equation}{section}

\DeclarePairedDelimiterX \abs    [1]   {\lvert}{\rvert}{\,#1\,}
\newcommand              \Aut          {\textnormal{Aut}}
\newcommand              \B      [1][] {\categ[#1]{\mathscr B}}
\DeclarePairedDelimiterX \bangle [1]   {\langle}{\rangle}{\,#1\,}
\renewcommand            \bar    [1]   {\overline{#1}}
\newcommand              \barcons[1]   {   \scalebox{1.2}{$\prec$}
							                          #1 \scalebox{1.2}{$\succ$}}
\newcommand              \C      [1][] {\categ[#1]{\mathscr C}}
\newcommand              \Cat          {\textnormal{\textbf{Cat}}}
\newcommand              \categ  [2][] {\ifx\\#1\\
                                          #2
                                        \else
                                          \ifthenelse
                                            { #1 = 0 }
                                            { \Ob #2 }
                                            { #2_#1  }
                                        \fi                 }
\newcommand              \cod          {\textnormal{cod}}

\newcommand              \Del          {\mathbf\Delta}
\newcommand              \dom          {\textnormal{dom}}

\newcommand              \ev           {\textnormal{ev}}

\newcommand              \faces  [1]   {\ulcorner#1\urcorner}
\newcommand              \fills        {\textnormal{fill}}
\newcommand              \id           {\textnormal{id}}
\newcommand              \img          {\textnormal{im}} 

\newcommand              \G            {\mathscr G}
\renewcommand            \H            {\mathscr H}
\renewcommand            \hom    [3][0]{\ifx0#1
                                          \textnormal{hom}(#2,#3)
                                        \else
                                          #1(#2,#3)
                                        \fi                       }
\newcommand              \Mor          {\mathcal Mor\,}
\newcommand              \mt     [2][] {\textnormal{mt}_{#1}(#2)}
\newcommand              \mono         {\hookrightarrow}
\newcommand              \N            {\mathbb N}

\newcommand              \ndgEpi       {\omega}
\newcommand              \ndgRoot      {\lambda}
\makeatother
\newcommand              \Ob           {\mathcal Ob\,}

\newcommand              \op     [1]   {{#1}^{op}}

\newcommand              \Sd           {\textnormal{Sd}}
\newcommand              \Sds          {\textnormal{Sd}^s}

\newcommand              \simp   [2]   {\textnormal{simp}_{#1}(\Sd#2)}
\newcommand              \sSet         {\textbf{sSet}}

\newcommand              \Z            {\mathbb Z}

\newlength\trilength
\newlength\shiftHeight
\NewDocumentCommand{\bigtri}{ O{} O{} O{} m m m }{
  \tikz[baseline=(baase),scale=.6, auto, every node/.style={transform shape}]{
   \path (0,-.2cm) node (baase) {};
   \useasboundingbox (-.7cm,.2cm) rectangle (1cm,-.2cm);
   \path	(0,0) node (middle) {}
   	(-30:\trilength) node (a) {$#1$}
   	( 90:\trilength) node (b) {$#2$}
   	(210:\trilength) node (c) {$#3$};
   \draw (b) to node[transform canvas={yshift=\shiftHeight},swap]
                                   {$#4$} (c);
   \draw (a) to node[transform canvas={yshift=\shiftHeight},swap]
                                   {$#5$} (b);
   \draw (a) to node[transform canvas={yshift=.03cm}]  {$#6$} (c);
   }}
\NewDocumentCommand{\filltri}{ O{} O{} O{} m m }{
  \tikz[baseline=(baase),scale=.6, auto, every node/.style={transform shape}]{
   \path (0,-.8mm) node (baase) {};
   \useasboundingbox (-.7cm,.2cm) rectangle (1cm,-.2cm);
   \path	(0,0) node (middle) {}
   	(-30:\trilength) node (a) {$#1$}
   	( 90:\trilength) node (b) {$#2$}
   	(210:\trilength) node (c) {$#3$};
   \draw (b) to node[transform canvas={yshift=\shiftHeight},swap]
                                   {$#4$} (c);
   \draw (a) to node[transform canvas={yshift=\shiftHeight},swap]
                                   {$#5$} (b);
   \draw[dashed] (a) to                           (c); 
  }}
\NewDocumentCommand{\commtri}{ O{} O{} O{} m m m O{->} O{->} O{->}}{  \tikz{
   \path	(0,0) node (middle) {}
   	(-30:\trilength) node (a) {$\scriptstyle{#1}$}
   	( 90:\trilength) node (b) {$\scriptstyle{#2}$}
   	(210:\trilength) node (c) {$\scriptstyle{#3}$};
   \draw[#7] (b) to node[transform canvas={yshift=\shiftHeight},swap]
                                   {$\scriptstyle{#4}$} (c);
   \draw[#8] (a) to node[transform canvas={yshift=\shiftHeight},swap]
                                   {$\scriptstyle{#5}$} (b);
   \draw[#9] (a) to node[transform canvas={yshift=.03cm}]  {$\scriptstyle{#6}$} (c);
   }}

\newenvironment{sixgraph}
          {\begin{tikzpicture}
                    \path node                  (a) {$*$}
                          node [right=5mm of a] (b) {$*$}
                          node [below=5mm of a] (c) {$*$}
                          node [below=5mm of b] (d) {$*$}
                          node [below=5mm of c] (e) {$*$}
                          node [below=5mm of d] (f) {$*$}
                    ;
          }{\end{tikzpicture}}
\newcommand\sixadbc  {\draw[shorten >=-2mm, shorten <=-2mm](a) -- (d);  
                      \draw[shorten >=-2mm, shorten <=-2mm](b) -- (c);}
\newcommand\sixafbe  {\draw[shorten >=-1mm, shorten <=-1mm](a) -- (f);  
                      \draw[shorten >=-1mm, shorten <=-1mm](b) -- (e);}
\newcommand\sixcfde  {\draw[shorten >=-2mm, shorten <=-2mm](c) -- (f);  
                      \draw[shorten >=-2mm, shorten <=-2mm](d) -- (e);}
\newcommand\sixacbd  {\draw[shorten >=-1mm, shorten <=-1mm](a) -- (c);  
                      \draw[shorten >=-1mm, shorten <=-1mm](b) -- (d);}
\newcommand\sixcedf  {\draw[shorten >=-1mm, shorten <=-1mm](c) -- (e);  
                      \draw[shorten >=-1mm, shorten <=-1mm](d) -- (f);}
\newcommand\sixaebf  {\draw[shorten >=-1mm, shorten <=-1mm](a) to [bend right](e);  
                      \draw[shorten >=-1mm, shorten <=-1mm](b) to [bend left] (f);}
\newcommand\sixabcdef{\draw[shorten >=-1mm, shorten <=-1mm](a) -- (b);  
                      \draw[shorten >=-1mm, shorten <=-1mm](c) -- (d);
                      \draw[shorten >=-1mm, shorten <=-1mm](e) -- (f);}

\newenvironment{eightgraph}
          {\begin{tikzpicture}
            \newdimen\clen
            \clen=.8cm
            \newdimen\dlen
            \dlen=.5cm
                    \node (a) at ( 0, \clen) {$*$};
                    \node (b) at ( 0, \dlen) {$*$};
                    \node (c) at (-\clen, 0) {$*$};
                    \node (d) at (-\dlen, 0) {$*$};
                    \node (e) at ( \clen, 0) {$*$};
                    \node (f) at ( \dlen, 0) {$*$};
                    \node (g) at ( 0,-\clen) {$*$};
                    \node (h) at ( 0,-\dlen) {$*$};
          }{\end{tikzpicture}}
\newcommand\eightfggf    {\draw                                      (a) -- (b);
                          \draw                                      (c) -- (d);
                          \draw                                      (e) -- (f);
                          \draw                                      (g) -- (h);}
\newcommand\eightfgfg    {\draw[shorten >=-1.25mm, shorten <=-1.25mm](a) -- (d);
                          \draw[shorten >=-1.25mm, shorten <=-1.25mm](b) -- (c);
                          \draw[shorten >=-1.25mm, shorten <=-1.25mm](e) -- (h);  
                          \draw[shorten >=-1.25mm, shorten <=-1.25mm](f) -- (g);}
\newcommand\eightgfgf    {\draw[shorten >=-1.25mm, shorten <=-1.25mm](a) -- (f);
                          \draw[shorten >=-1.25mm, shorten <=-1.25mm](b) -- (e);
                          \draw[shorten >=-1.25mm, shorten <=-1.25mm](c) -- (h);  
                          \draw[shorten >=-1.25mm, shorten <=-1.25mm](d) -- (g);}
\newcommand\eightfgfginv {\draw[shorten >=-1.25mm, shorten <=-1.25mm](a) to [bend right] (h);
                          \draw[shorten >=-1.25mm, shorten <=-1.25mm](b) to [bend left] (g);}
\newcommand\eightfgfginvtwo {\draw[shorten >=-1.25mm, shorten <=-1.25mm](a) to [bend left] (h); 
                          \draw[shorten >=-1.25mm, shorten <=-1.25mm](b) to [bend right] (g);}
\newcommand\eightfginvfg {\draw[shorten >=-1.25mm, shorten <=-1.25mm](c) to [bend right] (f); 
                          \draw[shorten >=-1.25mm, shorten <=-1.25mm](d) to [bend left] (e);}
\newcommand\eightfginvfgtwo {\draw[shorten >=-1.25mm, shorten <=-1.25mm](c) to [bend left] (f); 
                          \draw[shorten >=-1.25mm, shorten <=-1.25mm](d) to [bend right] (e);}
\newcommand\eightffid    {\draw[shorten >=-1.25mm, shorten <=-1.25mm](a) -- (c);
                          \draw[shorten >=-1.25mm, shorten <=-1.25mm](b) -- (d);
                          \draw[shorten >=-1.25mm, shorten <=-1.25mm](e) -- (g);
                          \draw[shorten >=-1.25mm, shorten <=-1.25mm](f) -- (h);}
\newcommand\eightggid    {\draw[shorten >=-1.25mm, shorten <=-1.25mm](a) -- (e);
                          \draw[shorten >=-1.25mm, shorten <=-1.25mm](b) -- (f);
                          \draw[shorten >=-1.25mm, shorten <=-1.25mm](c) -- (g);
                          \draw[shorten >=-1.25mm, shorten <=-1.25mm](d) -- (h);}
\newcommand\eightffgginv {\draw[shorten >=-1.25mm, shorten <=-1.25mm](a) to [bend left] (g);
                          \draw[shorten >=-1.00mm, shorten <=-1.00mm](b) -- (h);}
\newcommand\eightffgg    {\draw[shorten >=-1.25mm, shorten <=-1.25mm](c) to [bend left] (e); 
                          \draw[shorten >=-1.00mm, shorten <=-1.00mm](d) -- (f);}

\newlength{\arrow}
\newcommand*{\myrightarrow}{\xrightarrow{\mathmakebox[\arrow]{}}}
\newcommand*{\myleftarrow}{\xleftarrow{\mathmakebox[\arrow]{}}}
\newcommand*{\myxrightarrow}[1]{\xrightarrow{\mathmakebox[\arrow]{#1}}}
\newcommand*{\myxleftarrow}[1]{\xleftarrow{\mathmakebox[\arrow]{#1}}}

\makeatletter
\newcommand*{\shifttext}[2]{
  \settowidth{\@tempdima}{#2}%
  \makebox[\@tempdima]{\hspace*{#1}#2}%
}
\makeatother

\newlength\dlen      
\newlength\mylen

\begin{document}
\begin{abstract}
  Given groupoids $\G$ and $\H$ as well as an isomorphism
  $\Psi:\Sd\,\G\cong\Sd\,\H$ between subdivisions,
  we construct an isomorphism $P:\G\cong\H$.
  If $\Psi$ equals $\Sd F$ for some functor $F$,
  then the constructed isomorphism $P$ is equal to $F$.
  It follows that the restriction of $\Sd$ to the category of
  groupoids is conservative.
  These results do not hold for arbitrary categories.
\end{abstract} 
\maketitle
\tableofcontents


\section{Introduction}
\numberwithin{thm}{section}
\numberwithin{equation}{section}
The categorical subdivision functor ${\Sd:\Cat\to\Cat}$ is defined
as the composite $\Pi\circ\Sds\circ N$ of the nerve functor ${N:\Cat\to\sSet}$,
the simplicial barycentric subdivision functor ${\Sds:\sSet\to\sSet}$, and the
fundamental category functor ${\Pi:\sSet\to\Cat}$ (left adjoint to $N$).
The categorical subdivision functor $\Sd$ has deficiencies:
for example, it is neither a left nor a right adjoint.
Nevertheless, $\Sd$ bears similarities to its simplicial analog $\Sds$,
and has many interesting properties.
For example, it is known that the second subdivision $\Sd^2\C$
of any small category $\C$ is a poset \cite[ch. 13]{May}.

In general, there exist small categories $\B$ and $\C$
such that $\Sd\,\B$ is isomorphic to $\Sd\,\C$ but $\B$ is not isomorphic to $\C$
or to $\op\C$.
In this paper we show that such examples do
not occur in the category of small groupoids: if $\Sd\,\G$ is isomorphic
to $\Sd\,\H$ for groupoids $\G$ and $\H$, then there exists an
isomorphism between \mbox{$\G$ and $\H$.}

In broad strokes, the argument is as follows.
There is a canonical identification of objects in $\Sd\,\G$
with non-degenerate simplices in the nerve $N\G$.
Any isomorphism between $\Sd\,\G$ and $\Sd\,\H$ induces
a bijection between the 0-simplices of $N\G$ and the
0-simplices of $N\H$, and similarly for 1-simplices.
The objects and morphisms of a category correspond (respectively)
to the 0-simplices and 1-simplices in the nerve of that category, so for any
isomorphism ${\Psi:\Sd\,\G\to\Sd\,\H}$ there are induced bijections
${\Ob(\psi):\Ob(\G)\to\Ob(\H)}$ and ${\Mor(\psi):\Mor(\G)\to\Mor(\H)}$.
We will show that the restriction of these bijections to any connected
component of $\G$ determines a possibly-contravariant functor into $\H$.
Proceeding one component at a time and using the fact that any
groupoid is isomorphic to its opposite, we can use these maps
$\psi$ to construct a covariant isomorphism between $\G$ and $\H$.

\section{Notation and overview of proof}\label{sec:notn/overview}
In this section we introduce notation used throughout the paper.
We then give a sketch of the proof.

\subsection*{Notation}
We write $[n]$ for the totally-ordered poset category
having objects the non-negative integers $0,\ldots,n$.
If $i$ and $j$ are integers satisfying $0\leq i<j\leq n$,
we write $[i<j]$ to indicate the morphism from $i$ to $j$
in the category $[n]$.

Given an object $c$ in a small category $\C$,
we shall write $\bangle c$ for the functor $[0]\to\C$
that represents $c$.
Given a sequence $f_1,\cdots,f_n$ of morphisms in $\C$ satisfying
$\dom\,f_n=\cod\,f_{n-1}$ for $0<i\leq n$,
we will write $\barcons{f_n|\cdots|f_1}$ for
the functor $[n]\to\C$ that represents the given sequence
  \raisebox{-1mm}{$\xleftarrow{f_n}\cdots\xleftarrow{f_1}$}.

\subsection*{Overview}
\begin{itemize}[leftmargin=7mm]
\item  In Section \ref{sec:constructions}, we introduce the simplex category
  $\Del$ and the nerve functor $N$.
  We define the notion of degenerate simplices, and introduce $\Sd\,\C$
  as a category whose objects are the non-degenerate simplices of $N\C$.
  
\item  In Section \ref{sec:Lemmas}, we demonstrate that for any small
  category $\C$, there is an isomorphism ${\Sd\,\C\cong\Sd(\op\C)}$
  which sends $\bangle c$ to $\bangle c$
  and $\barcons{f_n|\cdots|f_1}$ to $\barcons{f_1|\cdots|f_n}$.
  For any small groupoid $\G$, there is an automorphism $\alpha_\G$ of $\Sd\G$
  that sends $\bangle c$ to $\bangle c$
  and $\barcons{f_n|\cdots|f_1}$
  to $\barcons{f_1^{-1}|\cdots|f_n^{-1}}$.
  We also show that $\C$ is connected if
  and only if $\Sd\,\C$ is connected, and that there is an
  isomorphism ${\Sd(\amalg_i\C_i)\cong\amalg_i(\Sd\,\C_i)}$
  for any set $\{\C_i\}$ of small categories.
  These facts will be useful for reducing the problem to the case
  of connected groupoids, and for obtaining covariant isomorphisms
  from contravariant ones.

\item  In Section \ref{sec:encoding}, we identify categorical properties
  of $\Sd\,\C$ that encode some structural aspects of $\C$.
  This will be useful later when we consider isomorphisms $\Sd\,\B\cong\Sd\,\C$
  and show that the identified categorical properties are preserved.
  \begin{itemize}
  \item Given an object $y$ in $\Sd\,\C$, we write $\mt[n]y$ for the set
      \[  \coprod_{x:[n]\to\C}\Sd\,\C(x,y) \]
    of \textbf morphisms \textbf targeting $y$ and having
    domain equal to some (non-degenerate) $n$-simplex.
    We write $\mt y$ for the set $\coprod_{n\in\N}\mt[n]y$
    of all morphisms targeting $y$ in $\Sd\,\C$.
  \item  We show that for any non-degenerate $n$-simplex $y:[n]\to\C$,
    regarded as an object in $\Sd\,\C$, the morphisms of $\Sd\,\C$
    targeting $y$ are in bijection with the monomorphisms of $\Del$
    targeting $[n]$.
    Thus, an object of $\Sd\,\C$ is an $n$-simplex if and only if
    it is the target of $2^{n+1}-1$ morphisms in $\Sd\,\C$.

  \item  Given an object $y$ in $\Sd\,\C$, we define the \textit{faces} of $y$
    to be the objects $x$ such that there exists a morphism $x\to y$.
    We show that the proper faces of a non-degenerate 1-simplex $\barcons f$
    are precisely the 0-simplices $\bangle{\dom\,f}$ and $\bangle{\cod\,f}$.

  \item  Given a 2-simplex $\barcons{f|g}$ in $\Sd\,\C$,
    we show that $f$ is left-inverse to $g$ if and only if there
    are exactly two morphisms in the set $\mt[1]{\barcons{f|g}}$,
    the domains of which are the 1-simplices $\barcons f$ and $\barcons g$.
    If $f$ is not inverse to $g$, then $\mt[1]{\barcons{f|g}}$
    contains three morphisms, the domains of which are $\barcons f$,
    $\barcons g$, and $\barcons{f\circ g}$.
  \end{itemize}

\item  Supposing that $\G$ is a small groupoid, we show
  in section \ref{sec:encoding2} how properties of $\Sd\,\G$
  encode the structure of $\G$ up to opposites.
  Supposing that $f$ and $g$ are non-identity morphisms in $\G$,
  the local structure of $\Sd\G$ near $\barcons f$ and
  $\barcons g$ determines whether a given 1-simplex $\barcons h$
  satisfies one of the equations $h=f\circ g$ or $h=g\circ f$.
  The two propositions stated at the end of section \ref{sec:encoding2} (and
  proved in the Appendix) are this paper's most difficult results.

\item  Section \ref{sec:Correspondence} concerns maps
  between subdivisions of categories.
  We suppose that $\G$ and $\H$ are groupoids, and that
  $\Psi:\Sd\,\G\to\Sd\,\H$ is an isomorphism.
  The results from Sections \ref{sec:encoding}
  and \ref{sec:encoding2} are used to show
  that the structure of $\G$ encoded by $\Sd\,\G$
  must match the structure of $\H$ encoded by $\Sd\,\H$.

  \begin{itemize}
  \item  For any objects $x$ and $y$ in $\Sd\,\G$, the map $\Psi$ induces a
  bijection between $\hom[\Sd\,\G]xy$ and $\hom[\Sd\,\H]{\Psi\,x}{\Psi\,y}$,
    hence faces of $y$ are sent to faces of $\Psi\,y$.
    Because the cardinality of $\mt y$ is equal that of $\mt{\Psi\,y}$,
    $\Psi$ sends $n$-simplices to $n$-simplices.
    It follows that the elements of $\mt[n]y$ are sent bijectively
    to elements of $\mt[n]{\Psi\,y}$.
  
  \item  We define the map $\psi$ so that\\[-1mm]
      \begin{minipage}[b]{\linewidth}
      \[\Psi\bangle c=\bangle{\psi\,c} 
        \hspace{15pt}\text{and}
        \hspace{15pt}\psi\,\id_a=\id_{\psi a}
        \hspace{15pt}\text{and}
        \hspace{15pt}\Psi\barcons f=\barcons{\psi\,f}\]
        \end{minipage}\\[1mm]
    are satisfied for each object $c$ and each
    non-identity morphism $f$ in $\G$.
    We show that this map $\psi:\G\to\H$ is an isomorphism between
    the (undirected) graphs that underlie $\G$ and $\H$.
  
  \item  The structure of $\Sd\,\G$ near given 1-simplices $\barcons f$,
    $\barcons g$, and $\barcons h$ is the same as that of $\Sd\,\H$
    near $\barcons{\psi\,f}$, $\barcons{\psi\,g}$ and $\barcons{\psi\,h}$.
    Therefore, for any $f$, $g$, and $h$ in $\G$, we have $h=f\circ g$
    or $h=g\circ f$ if and only if $\psi\,h=(\psi\,f)\circ(\psi\,g)$
    or $\psi\,h=(\psi\,g)\circ(\psi\,f)$.

  \item With help from a group-theoretic result due to Bourbaki, we show that
  if $\G$ is a connected single-object groupoid (that is, a group)
    then the map $\psi:\G\to\H$ is a possibly-contravariant isomorphism.
    We give an analogous result concerning connected
    groupoids that have multiple objects.

  \item  Working with arbitrary (non-connected) groupoids $\G$ and $\H$,
    we suppose that $\G$ equals the coproduct $\G_1\amalg\G_2$,
    and that $\psi$ is contravariant on $\G_1$ and covariant on $\G_2$.
    There exist some groupoids $\H_1$ and $\H_2$ such that $\H=\H_1\amalg\H_2$,
    and such that the composite isomorphism\\[-1mm]
    \begin{minipage}[b]{\linewidth}
      \[\Sd\,\G_1\amalg\Sd\,\G_2\xrightarrow\cong\Sd(\G_1\amalg\G_2)\xrightarrow\Psi
        \Sd(\H_1\amalg\H_2)\xrightarrow\cong\Sd\,\H_1\amalg\Sd\,\H_2\]
    \end{minipage}\\[1mm]
    restricts to isomorphisms $\Psi_1:\Sd\,\G_1\to\Sd\,\H_1$
    and $\Psi_2:\Sd\,\G_2\to\Sd\,\H_2$.
    We have maps $\psi_1$ and $\psi_2$, each of which is a restriction
    of $\psi$; by using the fact that $\G_1$ is isomorphic to its opposite
    category,
    we can flip the variance of $\psi_1$ to obtain a covariant isomorphism
    $\psi_1':\G_1\to\H_1$, and thus a covariant isomorphism
    $\psi_1'\amalg\psi_2$ between $\G$ and $\H$.
  \end{itemize}
 \item 
  Appendix \ref{appx:proof_of_compoo} is the combinatorial heart of this paper.
  It is devoted to the proof
  of Section \ref{sec:encoding2}'s two most involved results,
  which concern the way that $\Sd\,\G$ encodes the relationship between
  endomorphisms in a given groupoid $\G$.
\end{itemize}


\section{Construction of $\Sd\,\C$}\label{sec:constructions}
\numberwithin{thm}{subsection}
\numberwithin{equation}{subsection}
This section defines the nerve functor, explains what (non)degenerate
simplices in the nerve are, and gives a construction of the
categorical subdivision functor $\Sd:\Cat\to\Cat$.

\subsection{The Nerve of a Category}
For each non-negative integer $n$, let $[n]$ denote the poset category
whose objects are the integers $0,\ldots,n$ and whose morphisms
are given by the usual ordering on $\Z$.
The simplex category $\Del$ is the full subcategory of $\Cat$
whose objects are the posets $[n]$.
Note that $\Del$ is a concrete category:
for any $m$ and $n$, the morphisms $[m]\to[n]$ in $\Del$
can be identified with the order-preserving functions
$\{0,\ldots,m\}\to\{0,\ldots,n\}$.
Note also that epimorphisms in $\Del$ are just order-preserving surjections,
and monomorphisms in $\Del$ are just order-preserving injections.

The \textit{nerve} $N\C$ of a category $\C$ is the
restriction from $\op\Cat$ to $\op\Del$ of the hom-functor
$\hom[\Cat]-\C$.
The nerve construction can be made into a functor $N:\Cat\to\sSet$
by sending a functor $F:\B\to\C$ to the natural transformation
$N\B\Rightarrow N\C$ given by $\hom[\Cat]-F$.
The elements of $\hom[\Cat]{[n]}\C$ are
called the \textit{$n$-simplices} of $N\C$.
Given a morphism $\mu:[m]\to[n]$ in $\Del$, write $\mu^*$ for the function
  \[  \hom[\Cat]{[n]}\C\to\hom[\Cat]{[m]}\C  \]
defined by $\mu^*x=x\circ\mu$ for each $n$-simplex $x$ of $N\C$.

A simplex $y:[m]\to\C$ of $N\C$ is said to be \textit{degenerate}
if there exists some simplex $x$ of $N\C$ and some non-identity
epimorphism $\sigma$ in $\Del$ satisfying $\sigma^*x=y$.
If this simplex $x$ is not degenerate, then it is
called the \textit{non-degenerate root} of $y$.

\begin{prop}
  Given an $m$-simplex $y:[m]\to\C$ of $N\C$, the non-degenerate
  root of $y$ is unique.
  \begin{proof}
  This is a corollary of the standard general unique decomposition of
  a morphism in $\Del$ as a composite of an epimorphism and a monomorphism.
    Suppose that $x_1:[n_1]\to[m]$ and $x_2:[n_2]\to[m]$ are distinct
    simplices satisfying $\sigma_1^*x_1=y=\sigma_2^*x_2$ for some
    epimorphisms $\sigma_1$ and $\sigma_2$.
    We will prove that one of $x_1$ and $x_2$ must be degenerate.

    At least one of $\sigma_1$ and $\sigma_2$ must be a non-identity morphism,
    for otherwise $x_1=y=x_2$.
    Assume without loss of generality that $\sigma_2\neq\id$,
    and let $\nu_1:[m]\to[n_1]$ be right-inverse to $\sigma_1$.
    Then there is equality
      \[x_1=(\sigma_1\nu_1)^*x_1=\nu_1^*\sigma_1^*x=\nu_1^*y
        =\nu_1^*\sigma_2^*x_2=(\sigma_2\nu_1)^*x_2.\]
    In general, any morphism in $\Del$ has a unique decomposition
    as the composite of an epimorphism and a monomorphism.
    Therefore, we have $\sigma_2\nu_1=\nu'\sigma'$ for some monomorphism $\nu'$
    and some non-identity epimorphism $\sigma'$,
    and the resulting equality
      \[x_1=(\sigma_2\nu_1)^*x_2=(\nu'\sigma')^*x_2=\sigma'^*(\nu'^*x_2)\]
    proves that $x_1$ is degenerate.
  \end{proof}
\end{prop}

The above proposition holds for simplices in arbitrary simplicial
sets, not just for simplices in the nerve of a category.
We will use the notation $\lambda(y)$ to denote
the non-degenerate root of $y$.

No 0-simplex of $N\C$ is degenerate, and each 0-simplex $[0]\to\C$ may be
regarded as picking out a single object of $\C$.
For a positive integer $m$, each $m$-simplex $x:[m]\to\C$ is
identified with the sequence
     \begin{equation}\label{diag:chains_in_C}
     x(m)\myxleftarrow{f_m}x({m-1})\myleftarrow
      \cdots\myleftarrow x(1)\myxleftarrow{f_1}x(0)
     \end{equation}
of morphisms obtained by setting each $f_i$ equal to the image under $x$
of the morphism $[i-1<i]$ in $[m]$.
Such an $m$-simplex $x$ is non-degenerate if and only if
all of the arrows $f_i$ are non-identity morphisms of $\C$.

\subsection{Construction of $\Sd\,\C$}
Let $\C$ be a small category.
The objects of $\Sd\,\C$ are the non-degenerate simplices of $N\C$.
Given objects $x:[m]\to\C$ and $y:[n]\to\C$ of $\Sd\,\C$, the morphisms
$x\to y$ in $\Sd\,\C$ are the pairs $(\sigma,\nu)$ of morphisms in
$\Del$ satisfying $\cod\,\sigma=[m]$ and $\cod\,\nu=[n]$ such that
  \begin{enumerate}
    \item $\nu$ is a monomorphism,
    \item $\sigma$ is an epimorphism, and
    \item there is equality $\sigma^*x=\nu^*y$.
  \end{enumerate}

\begin{wrapfigure}[5]{r}{.175\textwidth}
    \vspace{-15pt}
    \centering
    \tikz{
      \path node                  (k) {$[k]$}
            node [right=7mm of k] (n) {$[n]$}
            node [below=7mm of k] (m) {$[m]$}
            node [right=7mm of m] (C) {$\C$};
              \begingroup \everymath{\scriptstyle}
      \draw[>->] (k) -- node                 {$\nu$}    (n);
      \draw[->>] (k) -- node [pos=.45, swap] {$\sigma$} (m);
      \draw[->]  (n) -- node [pos=.45]       {$y$}      (C);
      \draw[->]  (m) -- node [swap]          {$x$}      (C);
              \endgroup
    }
\end{wrapfigure}
\noindent
The identity morphism $x\to x$ is given by the pair $(\id_{[m]},\id_{[m]})$.
Condition (3) above requires that $\dom\,\sigma=[k]=\dom\,\nu$ for
some $k$, as in the commutative diagram to the right.
Note that for each morphism $(\sigma,\nu):x\to y$ in $\Sd\,\C$,
there is a multivalued function $\{0,\ldots,m\}\to\{0,\ldots,n\}$
given by $\nu\circ\sigma^{-1}$.
Here $\sigma^{-1}$ denotes the subset of $\{0,\ldots,m\}\times\{0,\ldots,k\}$
that is inverse to $\sigma$ (as a binary relation), and ${\nu\circ\sigma^{-1}}$
is the composite of binary relations $\nu$ and $\sigma^{-1}$.
Composition of morphisms in $\Sd\,\C$ is induced by composition
of such multivalued functions.
See Figure \ref{fig:a_morphism_in_SdC} on page \pageref{fig:a_morphism_in_SdC} 
for visualization of
a test morphism in $\Sd\,\C$.
Figure \ref{fig:two_morphisms_in_SdC} (on page
\pageref{fig:two_morphisms_in_SdC}) depicts two
composable morphisms in $\Sd\,\C$, and Figure
\ref{fig:composite_morphism_in_SdC} depicts the composite.

	\begin{exmp}\label{exmp:non-isomorphic-categories-w/-isomorphic-subds}
Let $\B$ be the category $(\cdot\to\cdot\leftarrow\cdot)$ with three objects
and two non-identity arrows, and let $\C$ be the category
$(\cdot\leftarrow\cdot\to\cdot)$ opposite to $\B$.
We see that $\B$ and $\C$ are not isomorphic, yet the barycentric subdivisions
$\Sd\,\B$ and $\Sd\,\C$ are both isomorphic to the category
$(\cdot\to\cdot\leftarrow\cdot\to\cdot\leftarrow\cdot)$ with five objects and
four non-identity arrows.
We will later prove that $\Sd\,\C\cong\Sd(\op\C)$ for any category $\C$.
\end{exmp}

The following example is due to
Jonathan Rubin:\footnote{Private communication.}
\begin{exmp}\label{exmp:SdB=SdC,B=/=C,B=/=Cop}
Consider $\N$ and $\Z$ as poset categories.
There is no isomorphism between $\N$ and $\Z$,
yet the subdivisions $\Sd\,\N$ and $\Sd\,\Z$ are isomorphic.

Generally, if $T$ is a totally-ordered set (regarded as a poset category)
then the $n$-simplices of $\Sd\,T$ are all the linearly-ordered
subsets of $T$ having $n+1$ elements.
If $x$ and $y$ are objects in $\Sd\,T$ and
if $\bar x$ and $\bar y$ are the corresponding subsets of $T$,
then there exists a morphism $x\to y$ in $\Sd\,T$ if and only if
$\bar x$ is a subset of $\bar y$.
Thus $\Sd\,T$ is the poset
  \[  \Bigl(\{F\subseteq T~|~0<\abs F<\infty\},\subseteq\Bigr)\]
of finite non-empty subsets of $T$, ordered by inclusion.

If $T$ and $T'$ are totally-ordered sets
then any set bijection between $T$ and $T'$ induces an isomorphism
$\Sd\,T\cong\Sd\,T'$ of categories, regardless of whether
the given bijection between $T$ and $T'$ is order-preserving.
\end{exmp}

\subsection{Concerning $\Sd F$}

Let $\B$ and $\C$ be small categories,
and let $F$ be a functor $\B\to\C$.
The functor $\Sd\,F$ sends each 0-simplex $y:[0]\to\B$
to the 0-simplex $F\circ y:[0]\to\C$.
Using the notation from Section \ref{sec:notn/overview},
this means that $\bangle b$ is sent to $\bangle{F(b)}$ for each object
$b$ of $\B$.
For larger $m$, $\Sd F$ sends each $m$-simplex (\ref{diag:chains_in_C}) to
the subsequence of
    \[ F\bigl(x(m)\bigr)\myxleftarrow{F(f_m)}F\bigl(x({m-1})\bigr)\myleftarrow
      \cdots\myleftarrow F\bigl(x(1)\bigr)\myxleftarrow{F(f_1)}F\bigl(x(0)\bigr) \]
consisting of non-identity arrows.
Formally, each $m$-simplex $x:[m]\to\B$ is sent by $\Sd F$
to the non-degenerate root of $F\circ x$.
Multivalued functions can be used to
calculate where $\Sd F$ sends the morphisms of $\Sd\,\B$.
For the purposes of this paper, a precise formulation of $\Sd F$'s value
on morphisms of $\Sd\,\B$ is not necessary.

\begin{exmp}
Below is a depiction of a sample morphism $x\to y$ in $\Sd\,\C$.
\nopagebreak
\begin{figure}[H]
\centering
  	\begin{tikzpicture}[scale=0.9]
	\draw
			(-1.5, 0) node (high) {$[4]$}
			(-1.5,-3) node (low)  {$[1]$}
      (7.5 ,-3) node (x)    {$x$}
      (7.5 , 0) node (y)    {$y$}
      (   0, 0) node (a4)   {$d$}
			( 1.5, 0) node (a3)   {$c$}
			(   3, 0) node (a2)   {$a$}
			( 4.5, 0) node (a1)   {$b$}
			(   6, 0) node (a0)   {$a$}
			(2.25,-3) node (c4)   {$d$}
			(3.75,-3) node (c0)   {$a$}
	    ($(c4)!0.5!(a4)$) node (b4) {$d$}
			($(c0)!0.5!(a2)$) node (b2) {$a$}
			($(c0)!0.5!(a0)$) node (b0) {$a$}
			($(high)!0.5!(low)$) node (mid) {$[2]$};
          \begingroup \everymath{\scriptstyle}
	\draw[->] (a0) to node [swap] {$f$}      (a1);
	\draw[->] (a1) to node [swap,pos=.4] {$f^{-1}$} (a2);
	\draw[->] (a2) to node [swap] {$g$}      (a3);
	\draw[->] (a3) to node [swap] {$h$}      (a4);
	\draw[->] (c0) to node {$h\circ g$} (c4);
	\draw[->] (b0) to node [swap] {$\id_a$} (b2);
	\draw[->] (b2) to node [swap] {$h\circ g$} (b4);
	\draw[dashed]	(c0) to (b0) to (a0)
					(c0) to (b2) to (a2)
					(c4) to (b4) to (a4);
	\draw[>->] (mid) to node {$\nu$} (high);
	\draw[->>] (mid) to node [swap] {$\sigma$} (low);
  \draw[->] (x) to node [swap] {$(\sigma,\nu)$} (y);
          \endgroup
	\end{tikzpicture}
   \caption{\label{fig:a_morphism_in_SdC}}
\end{figure}
   \newcommand*{\farrow}[1]{\xleftarrow{\mathmakebox[\widthof{\scriptsize$f^{-1}$}]{#1}}}
   \noindent
    Here ${x:[1]\to\C}$ picks out the morphism $h\circ g$,
    and $y:[4]\to\C$ picks out the sequence
    \vspace{-3pt}
    \begin{center}
      $\cdot\farrow h \cdot\farrow g \cdot\farrow{f^{-1}} \cdot\farrow f \cdot$
    \end{center}
    \vspace{-3pt}
    of arrows in $\C$.
    The dashed lines illustrate how $h$ and $g$ compose to $h\circ g$, and how
    $f^{-1}$ and $f$ compose to $\id_a$.
    The injection $\nu:[2]\to[4]$ has image $\{0,2,4\}$, and
    $\sigma:[2]\to[1]$ is the surjection satisfying $\sigma(0)=0=\sigma(1)$.
    We have assumed that $h$ is not left-inverse to $g$,
    and that none of $f$, $f^{-1}$, $g$, $h$ are identity morphisms.
    The dashed lines determine a multivalued function from
    $\{0,1\}$ to $\{0,1,2,3,4\}$.
\end{exmp}

\begin{exmp}
  Below, two composable morphisms $x\to y$ and $y\to z$
  in $\Sd\,\C$ are depicted end-to-end.
\nopagebreak
\begin{figure}[H]
\centering
	\begin{tikzpicture}[scale=0.9]
	\draw	(0   ,2)    node (z6) {}
			  (1.5 ,2)    node (z5) {}
			  (3   ,2)    node (z4) {}
			  (4.5 ,2)    node (z3) {}
			  (6   ,2)    node (z2) {}
			  (7.5 ,2)    node (z1) {}
			  (9   ,2)    node (z0) {}
			  (10.5,2)    node (zx) {}
			  (2.25,0)    node (a3) {}
			  (3.75,0)    node (a2) {}
			  (5.25,0)    node (a1) {}
			  (6.75,0)    node (a0) {}
			  (8.25,0)    node (ax) {}
			  (3.75,-1.5) node (c4) {}
			  (5.25,-1.5) node (c0) {}
        (11.5,2)    node (z) {$z$}
        (11.5,0)    node (y) {$y$}
        (11.5,-1.5) node (x) {$x$}
        ;
          \begingroup \everymath{\scriptstyle}
	\draw[->] (ax) to node [swap,pos=.4] {$f$}        (a0);
	\draw[->] (a0) to node [swap,pos=.4] {$g$}        (a1);
	\draw[->] (a1) to node [swap,pos=.6] {$k$}        (a2);
	\draw[->] (a2) to node [swap]        {$k^{-1}$}   (a3);
	\draw[->] (c0) to node [swap,pos=.4] {$g$}        (c4);
	\draw[->] (zx) to node [swap]        {$f$}        (z0);
	\draw[->] (z0) to node [swap]        {$g$}        (z1);
	\draw[->] (z1) to node [swap,pos=.4] {$h^{-1}$}   (z2);
	\draw[->] (z2) to node [swap]        {$h$}        (z3);
	\draw[->] (z3) to node [swap]        {$l\circ k$} (z4);
	\draw[->] (z4) to node [swap,pos=.4] {$l^{-1}$}   (z5);
	\draw[->] (z5) to node [swap,pos=.4] {$k^{-1}$}   (z6);
	\draw[dashed]
          (c0) to (a0)
					(c4) to (a1)
					(c4) to (a3)
					(ax) to (zx)
					(a0) to (z0)
					(a1) to (z1)
					(a1) to (z3)
					(a2) to (z5)
					(a3) to (z6);
  \draw[->] (x) to (y);
  \draw[->] (y) to (z);
          \endgroup
	\end{tikzpicture}
	\caption{\label{fig:two_morphisms_in_SdC}}
  \end{figure}
  \noindent
  The multivalued functions associated to these morphisms
  $x\to y$ and $y\to z$ can be used
  to calculate the composite arrow in $\Sd\,\C$.
  This composite is depicted below.
\nopagebreak
\begin{figure}[H]
\centering
	\begin{tikzpicture}[scale=0.9]
	\draw	(0   ,2)    node (z6) {}
			  (1.5 ,2)    node (z5) {}
			  (3   ,2)    node (z4) {}
			  (4.5 ,2)    node (z3) {}
			  (6   ,2)    node (z2) {}
			  (7.5 ,2)    node (z1) {}
			  (9   ,2)    node (z0) {}
			  (10.5,2)    node (zx) {}
			  (3.75,-1.5) node (c4) {}
			  (5.25,-1.5) node (c0) {}
        (11.5,2)    node (z) {$z$}
        (11.5,-1.5) node (x) {$x$}
      ;
          \begingroup \everymath{\scriptstyle}
	\draw[->] (c0) to node [swap,pos=.4] {$g$}        (c4);
	\draw[->] (zx) to node [swap]        {$f$}        (z0);
	\draw[->] (z0) to node [swap]        {$g$}        (z1);
	\draw[->] (z1) to node [swap,pos=.4] {$h^{-1}$}   (z2);
	\draw[->] (z2) to node [swap]        {$h$}        (z3);
	\draw[->] (z3) to node [swap]        {$l\circ k$} (z4);
	\draw[->] (z4) to node [swap,pos=.4] {$l^{-1}$}   (z5);
	\draw[->] (z5) to node [swap,pos=.4] {$k^{-1}$}   (z6);
	\draw[dashed]	(c0) to (z0)
					(c4) to (z1)
					(c4) to (z3)
					(c4) to (z6);
  \draw[->] (x) to (z);
          \endgroup
	\end{tikzpicture}
	\caption{\label{fig:composite_morphism_in_SdC}}
\end{figure}
\end{exmp}

\begin{exmp}
  The figure below depicts a sample calculation.
  The functor $\Sd F:\Sd\,\B\to\Sd\,\C$ is applied
  to a morphism $x\to y$ of $\Sd\,\B$.
\nopagebreak
\begin{figure}[H]
\newlength\SdFV
\SdFV=.5in
  \begin{tikzpicture}[scale=.5, every node/.style={scale=0.6}]
  \path               node (y0) {}       			
        +(225:\SdFV)  node (y1) {}
        +(225:2\SdFV) node (y2) {}
        +(225:3\SdFV) node (y3) {}
        +(225:4\SdFV) node (y4) {}
        +(225:5\SdFV) node (y5) {}
        +(45:\SdFV) node (yA) {}
        ($(225:5\SdFV)+(180:\SdFV)$) node (yS) {}
        ($(y0)+(270:3.5in)+(225:1.5\SdFV)$) node (x0) {}	
        +(225: \SdFV)     node (x1) {}
        +(225:2\SdFV)     node (x2) {}
        +(225:3\SdFV)     node (x3) {}
        ($(yS)+(270:3.5in)$) node (xS) {}
        (0:2in)     node (Foy0) {}             			
        +(225:\SdFV) node (Foy1) {} 
        +(225:2\SdFV)  node (Foy2) {}
        +(225:3\SdFV) node (Foy3) {}
        +(225:4\SdFV) node (Foy4) {}
        +(225:5\SdFV) node (Foy5) {}
        +(225:-\SdFV) node (FoyA) {}
        +(225:6\SdFV) node (FoyB) {}
        ($(x0)+(0:2in)$)   node (Fox0)  {}			
        +(225:\SdFV)   node (Fox1)  {}
        +(225:2\SdFV)   node (Fox2)  {}
        +(225:3\SdFV)   node (Fox3)  {}
        +(225:4\SdFV) node (FoxB) {}
        ($(FoyA)+(330:2in)+(225:.5\SdFV)$) node (LFoyA) {}	
        ++(225:\SdFV)   node (LFoy0) {} 
        ++(225:\SdFV) node (LFoy1) {}
        +(225:\SdFV)  node (LFoy2) {}
        +(225:2\SdFV) node (LFoy3) {}
        +(225:3\SdFV) node (LFoy4) {}
        +(225:4\SdFV) node (LFoyB) {}
        ($(yS)+(0:6in)$) node (LFoyS) {}
		($(LFoyS)+(270:3.5in)$) node (LFoxS) {}
        ($(FoxB)+(330:2in)+(45:.5\SdFV)$) node (LFoxB) {}	
        ++(45:3\SdFV) node (LFox0) {}
        +(225:\SdFV)  node (LFox1) {}
        +(225:2\SdFV) node (LFox2) {}
        ;
  \draw[->] (y0) -- node [swap] {$k$} (y1);			
  \draw[->] (y1) -- node [swap] {$j$} (y2);
  \draw[->] (y2) -- node [swap] {$i$} (y3);
  \draw[->] (y3) -- node [swap] {$h$} (y4);
  \draw[->] (y4) -- node [swap] {$g$} (y5);
  \draw[->] (Foy0) -- node {$Fk$} (Foy1);			
  \draw[->] (Foy1) -- node {$Fj$} (Foy2);
  \draw[->] (Foy2) -- node {$\id$} (Foy3);
  \draw[->] (Foy3) -- node {$Fh$} (Foy4);
  \draw[->] (Foy4) -- node {$Fg$} (Foy5);
  \draw (y0) -- (Foy0)                   			
        (y1) -- (Foy1)
        (y2) -- (Foy2)
        (y3) -- (Foy3)
        (y4) -- (Foy4)
        (y5) -- (Foy5);
  \draw[->] (LFoy0) -- node [swap] {$Fk$} (LFoy1);		
  \draw[->] (LFoy1) -- node [swap] {$Fj$} (LFoy2);
  \draw[->] (LFoy2) -- node [swap] {$Fh$} (LFoy3);
  \draw[->] (LFoy3) -- node [swap] {$Fg$} (LFoy4);
  \draw[dashed] (Foy0) -- (LFoy0)  				
                (Foy1) -- (LFoy1)
                (Foy2) -- (LFoy2)
                (Foy3) -- (LFoy2)
                (Foy4) -- (LFoy3)
                (Foy5) -- (LFoy4);
  \draw[->] (LFox0) -- node [transform canvas={yshift=.1cm}] {$Fj$} (LFox1); 			
  \draw[->] (LFox1) -- node [transform canvas={yshift=.1cm}] {$F(g\circ h)$} (LFox2);
  \draw[dashed] (LFoy1) -- (LFox0) 				
                (LFoy2) -- (LFox1)
                (LFoy4) -- (LFox2);
  \draw[->] (Fox0) -- node [swap,transform canvas={xshift=.-.1cm, yshift=-.17cm}] {$Fj$} (Fox1); 
  \draw[->] (Fox1) -- node [swap,transform canvas={xshift=.-.1cm, yshift=-.17cm}] {$\id$} (Fox2);
  \draw[->] (Fox2) -- node [swap,transform canvas={xshift=.-.1cm, yshift=-.17cm}] {$F(g\circ h)$} (Fox3);
  \draw[dashed] (Fox0) -- (LFox0)  
                (Fox1) -- (LFox1)
                (Fox2) -- (LFox1)
                (Fox3) -- (LFox2);
  \draw[->] (x0) -- node [transform canvas={xshift=.1cm, yshift=.17cm}] {$j$} (x1); 
  \draw[->] (x1) -- node [transform canvas={xshift=.1cm, yshift=.17cm}] {$i$} (x2);
  \draw[->] (x2) -- node [transform canvas={xshift=.1cm, yshift=.17cm}] {$g\circ h$} (x3);
  \draw (x0) -- (Fox0)  
                (x1) -- (Fox1)
                (x2) -- (Fox2)
                (x3) -- (Fox3);
  \draw[dashed] (y1) -- (x0)  
                (y2) -- (x1)
                (y3) -- (x2)
                (y5) -- (x3);
  \begingroup
  \path ($(y5)   +(205:1.25in   )$) node (y)    {$\scalebox{1.5}{$y$}$}    	
        ($(x3)   +(205:1.25in   )$) node (x)    {$\scalebox{1.5}{$x$}$}
        ($(Foy0) +(20 :1.25in)$)    node (Foy)  {$\scalebox{1.5}{$F\circ y$}$}
        ($(LFoy0)+(20 :1.25in)$)    node (LFoy) {$\scalebox{1.5}{$\ndgRoot({F\circ y})$}$}
        ($(LFox0)+(20 :1.25in)$)    node (LFox) {$\scalebox{1.5}{$\ndgRoot({F\circ x})$}$}
        ($(Fox3) +(205:1.25in)$)    node (Fox)  {$\scalebox{1.5}{$F\circ x$}$}
        ;
  \endgroup
  \draw[->,thick] (y)    to [out=0, in=225]  (y5)   ; 
  \draw[->,thick] (Foy)  to [out=180, in=45] (Foy0) ; 
  \draw[->,thick] (LFoy) to [out=180, in=45] (LFoy0); 
  \draw[->,thick] (x)    to [out=0, in=225]  (x3)   ; 
  \draw[->,thick] (Fox)  to [out=0, in=225]  (Fox3) ; 
  \draw[->,thick] (LFox) to [out=180, in=45] ($(LFox0)+(.3cm,.3cm)$); 
  \begingroup
    \draw[->,thick,transform canvas={xshift=.17cm,yshift=.17cm}] ($(FoxB)!.2!(LFoxB)$)
          -- node [swap] {$\scalebox{1.5}{$\ndgEpi_1$}$} ($(FoxB)!.8!(LFoxB)$);
    \draw[->,thick,transform canvas={xshift=.17cm,yshift=.17cm}] ($(FoyB)!.2!(LFoyB)$)
          -- node [swap] {$\scalebox{1.5}{$\ndgEpi_2$}$} ($(FoyB)!.8!(LFoyB)$);
  \endgroup
  \end{tikzpicture}
  \caption{\label{fig:SdF(V)}}
\end{figure}
\noindent
  The given morphism $x\to y$
  is sent to a morphism $\ndgRoot({F\circ x})\to\ndgRoot({F\circ y})$
  in $\Sd\,\C$, where $\ndgRoot({F\circ x})$ is the non-degenerate root of $x$
  (and similarly for $\ndgRoot({F\circ y})$).
  The epimorphisms $\ndgEpi_1$ and $\ndgEpi_2$ are the unique
  maps in $\Del$ satisfying $\ndgEpi_1^*\ndgRoot({F\circ x})=F\circ x$
  and $\ndgEpi_2^*\ndgRoot({F\circ y})=F\circ y$.
  We suppose here that $Fi$ is an identity morphism in $\C$,
  and that $Fg$ is not left inverse to $Fh$.
\end{exmp}


\section{$\Sd$ preserves coproducts and identifies opposite categories}
\label{sec:Lemmas}
\numberwithin{thm}{section}
\numberwithin{equation}{section}

This section states two lemmas that will be used later in the paper.
In particular, we show in Lemma \ref{result:SdC=SdCop} that
the functor $\Sd$ does not distinguish between a category $\C$
and its opposite category $\op\C$, and in Lemma
\ref{result:connected/coproducts_preserved} that $\Sd$ preserves coproducts.
The latter lemma allows for reduction of this paper's main theorem
to the case of connected groupoids.
The former highlights a fundamental issue: for any small category
$\C$ we have an isomorphism $\Sd\,\C\to\Sd(\op\C)$.
Therefore, a naive attempt to construct a map $\B\to\C$ from a given
isomorphism $\Sd\,\B\to\Sd\,\C$ could result in contravariance.
As mentioned in the introduction, this pitfall can be avoided
when working with groupoids, as any groupoid is isomorphic to its opposite
via inversion.

\begin{lem}\label{result:SdC=SdCop}
  Let $\C$ be a small category. There is a canonical isomorphism
  \begin{equation}\label{eqn:SdC=SdCop}
    \Sd\,\C\xrightarrow\cong\Sd(\op\C)
  \end{equation}
  between the subdivision of $\C$ and the subdivision of the opposite category.
  \begin{proof}
    The claimed isomorphism sends each $m$-simplex
    (\ref{diag:chains_in_C}) in $\Sd\,\C$ to the $m$-simplex
      \begin{equation*}
      x(m)\myxrightarrow{f_m}x({m-1})\myrightarrow
       \cdots\myrightarrow x(1)\myxrightarrow{f_1}x(0)
      \end{equation*}
    in $\Sd(\op\C)$.
    If $\mu$ is a arrow $[m]\to[n]$ in $\Del$, write $\mu'$ for the
    arrow $[m]\to[n]$ defined by $\mu'(m-i)=n-\mu(i)$.
    The claimed map (\ref{eqn:SdC=SdCop}) is defined on morphisms of $\Sd\,\C$
    by $(\sigma,\nu)\mapsto(\sigma',\nu')$.
  \end{proof}
\end{lem}

Thus, for any statement about the relationship between $\Sd\,\C$ and $\C$,
there is a dual statement about $\Sd\,\C$ and $\op\C$.

\begin{exmp}
For any small groupoid $\G$ we have the inversion isomorphism
${\op\G\to\G}$ defined on morphisms by $f\mapsto f^{-1}$.
The subdivision of this isomorphism is a map $\Sd(\op\G)\to\Sd\,\G$,
which can be composed with the map $\Sd\,\G\to\Sd(\op\G)$ from Lemma
\ref{result:SdC=SdCop} to obtain an isomorphism $\alpha_\G:\Sd\,\G\to\Sd\,\G$.
This functor $\alpha_\G$ sends each $m$-simplex (\ref{diag:chains_in_C})
of $\Sd\,\G$ to the $m$-simplex
    \[
      x(m)\myxrightarrow{f^{-1}_m}x({m-1})\myrightarrow
       \cdots\myrightarrow x(1)\myxrightarrow{f^{-1}_1}x(0).
    \]

Note that $\Sd$ is faithful, as any functor $F:\B\to\C$
is determined by where $\Sd F$ sends the 0 and 1-simplices of $\Sd\,\B$.
Unless the groupoid $\G$ is a discrete category,
the map $\alpha_\G$ is not equal to the subdivision of any
automorphism $\G\to\G$.
Thus, for any automorphism $\xi$ of a non-discrete groupoid $\G$
we find an automorphism $\alpha_\G\circ\Sd\,\xi$ of the category $\Sd\,\G$.
Therefore, we have the inequality
  \[  \abs{\Aut(\Sd\,\G)}\geq2\cdot\abs{\Aut(\G)}.   \]
in the case where $\G$ is not discrete.
\end{exmp}

\begin{lem}\label{result:connected/coproducts_preserved}
The functor $\Sd$ preserves coproducts.
A category $\C$ is connected if and only
if its subdivision $\Sd\,\C$ is connected.
  \begin{proof}
  Subdivision $\Sd$ is given by the composite $\Pi\circ\Sds\circ N$.
  Here $\Pi:\sSet\to\Cat$ is left adjoint to $N$,
  and $\Sds:\sSet\to\sSet$ denotes barycentric subdivision of simplicial sets.
  The functors $\Pi$ and $\Sds$ are left adjoints,
  and the nerve functor $N$ preserves coproducts, so $\Sd$ does too.
  Therefore, if $\Sd\,\C$ is connected then $\C$ must be connected, for otherwise
  we would have $\C=\amalg_i\C_i$ and hence $\Sd\,\C=\amalg_i\Sd\,\C_i$.

  Suppose now that $\C$ is a connected category.
  To prove that $\Sd\,\C$ is connected we will introduce some new notation.
  For an object $c$ in $\C$, let $\bangle c$ denote the 0-simplex
  $[0]\to\C$ sending $0$ to $c$.
  For a morphism $f$ in $\C$, let $\barcons f$ denote
  the 1-simplex $[1]\to\C$ that sends the morphism $[0<1]$ to $f$.

  For any non-identity morphism $f$ in $\C$ there are arrows
  \[\bangle{\dom\,f}\xrightarrow{(\id_{[0]},\delta)}\barcons f
      ~\hspace{7pt}\text{ and }~\hspace{7pt}
    \bangle{\cod\,f}\xrightarrow{(\id_{[0]},\delta')}\barcons f\]
    in $\Sd\,\C$, where $\delta$ and $\delta'$ are the
    monomorphisms $[0]\to[1]$ having respective images $\{0\}$ and $\{1\}$.
    Thus, for any sequence
      \[\tikz{\path (0,0) node (a) {$a$}
                          node [right=7mm of a] (b) {$b$}
                          node [right=7mm of b] (c) {$c$}
                          node [right=7mm of c] (d) {$\cdots$};
                          \begingroup\everymath{\scriptstyle}
              \draw (a) -- node {$f$} (b);
              \draw (b) -- node {$g$} (c);
              \draw (c) -- node {$h$} (d);
                          \endgroup
      }\]
  of objects in $\C$, there is a sequence
    \[
       \bangle a \rightarrow \barcons f \leftarrow
       \bangle b \rightarrow \barcons g \leftarrow
       \bangle c \rightarrow \barcons h \leftarrow
       \cdots
    \]
  of objects in $\Sd\C$.
  Therefore, all 0-simplices are in the same connected component.
  To complete the proof, let $\delta:[0]\to[m]$ demote the morphism
  that sends $0$ to $0$, and note that for
  any $m$-simplex $x$ in $\Sd\,\C$ there is a morphism $(\id_{[0]},\delta)$
  whose domain is the 0-simplex $\delta^*x$ and whose codomain is $x$.
  Thus all object of $\Sd\,\C$ belong to a single connected component.
  \end{proof}
\end{lem}

Note that the single-object category $[0]$
is isomorphic to its subdivision $\Sd\,[0]$.
It follows from the previous lemma that any discrete
category $\C$ is isomorphic to its own subdivision.


\section{$\Sd\,\C$ encodes objects and arrows of $\C$}\label{sec:encoding}
\numberwithin{thm}{subsection}
\numberwithin{equation}{subsection}

This section concerns the relationship between a category and its subdivision.
We identify categorical properties of $\Sd\,\C$
that correspond to certain structures in the category $\C$.
The identified properties of $\Sd\,\C$ are preserved by
isomorphism:
supposing that $\Sd\,\B\cong\Sd\,\C$ for some category $\B$,
the structures of $\C$ will appear the same was in $\B$.
For example, $\B$ is a groupoid if and only if $\C$ is a groupoid.

\begin{notn}\label{notn:bangle,barcons,S*barcons,ndg_root}
As in Section \ref{sec:constructions}, for any objects $x$ and $y$
in $\Sd\,\C$, we regard the morphisms $x\to y$ in $\Sd\,\C$ as pairs
$(\sigma,\nu)$ of morphisms in $\Del$ satisfying
  \begin{enumerate}
    \item $\nu$ is a monomorphism,
    \item $\sigma$ is an epimorphism, and
    \item there is equality $\sigma^*x=\nu^*y$.
  \end{enumerate}

As in Section \ref{sec:Lemmas}, for any object $a$ of $\C$
we let $\bangle a:[0]\to\C$ denote the 0-simplex that represents $a$.
Given a morphism $f$ in $\C$, write $\barcons f:[1]\to\C$
for the 1-simplex that represents $f$.
  We extend this notation as follows:
  given a sequence $f_1,\ldots,f_m$ of morphisms in $\C$
  satisfying $\dom\,f_{n+1}=\cod\,f_n$ for $1\leq n<m$,
  write $\barcons{f_m|\cdots|f_1}$ for the $m$-simplex
  $x:[m]\to\C$ given by the diagram
    \[ x(m)\myxleftarrow{f_m}x({m-1})\myleftarrow
      \cdots\myleftarrow x(1)\myxleftarrow{f_1}x(0) \]
  in $\C$, where $f_i$ is equal to $x([i-1<i])$ for each $i$ above.
  This notation is inspired by the ``bar construction'' on groups.

       \begin{window}[0,r,
  \xymatrix{ [2] \ar[d]_x & [1] \ar[l]_{\{0,1\}} \ar[dl]^{\barcons f} \\ \C }
             ,{}] 
  For a natural number $m$, identify each nonempty subset
  $S$ of $\{0,\ldots,m\}$ with the monomorphism in $\Del$ having
  domain $\bigl[\abs S -1\bigr]$, codomain $[m]$, and image $S$.
  For example, if $x$ equals $\barcons{g|f}$ for some morphisms
  $f:a\to b$ and $g:b\to c$, then we can identify $\{0,1\}$ and $\{0,2\}$
  with injections $[1]\to[2]$ to obtain equalities $\{0,1\}^*x=\barcons f$
  and $\{0,2\}^*x=\barcons{g\circ f}$.
    \end{window}

  Recall that the \textit{non-degenerate root} of
  an $m$-simplex $y:[m]\to\C$ is the unique non-degenerate simplex
  $\ndgRoot({y}):[n]\to\C$ that satisfies $\ndgEpi({y})^*\ndgRoot({y})=y$
  for some surjection $\ndgEpi({y}):[m]\twoheadrightarrow[n]$.
  Setting $y=\barcons{f_m|\cdots|f_1}$, we obtain $\ndgRoot({y})$
  by omitting identity arrows from the sequence $f_1,\ldots,f_m$.
  The map $\ndgEpi({y})$ is uniquely determined by the equations
    \begin{equation}
     \ndgEpi({y})(0)=0
    \end{equation}
  and
    \begin{equation}\label{eqn:ndgEpi}
     \ndgEpi({y})(i+1)=\begin{cases}
        \ndgEpi({y})(i)   &\text{if $f_{i+1}$ is an identity morphism}\\
        \ndgEpi({y})(i)+1 &\text{if $f_{i+1}$ is a non-identity morphism}
       \end{cases}
    \end{equation}
  where (\ref{eqn:ndgEpi}) above holds for all $i<m$.
  An $m$-simplex $y:[m]\to\C$ may be regarded as an object
  of $\Sd\,\C$ precisely when it is non-degenerate.
\end{notn}

\begin{defn}\label{defn:mt(y)}
  Given an object $y$ of $\Sd\,\C$, let $\mt y$ denote the set
    \[ \{f\in\Mor(\Sd\,\C)~|~\cod\,f=y\} \]
  of \textbf morphisms \textbf targeting $y$ in $\Sd\,\C$.
\end{defn}

\begin{prop}\label{result:inj=mt}
  Given an $m$-simplex $y$, regarded as an object of $\Sd\,\C$,
  there is a bijection between the set
  of nonempty subsets of $\{0,\ldots,m\}$
  and the set $\mt y$ of morphisms targeting $y$.
  This bijection sends a subset $S$ of $\{0,\ldots,m\}$ to the morphism
    $  (\ndgEpi({S^*y}),S):\ndgRoot({S^*y})\to y
    $
  in $\Sd\,\C$.
  \begin{proof}
    Identifying each subset $S$ with the monomorphism
    $\bigl[\abs S -1\bigr]\to[m]$ having image $S$,
    it suffices to show that the map
      \[ \bigl\{\text{monomorphisms of $\Del$ targeting $[m]$}\bigr\}\to\mt y
      \]
    given by $S\mapsto(\ndgEpi({S^*y}),S)$ is bijective.
    This pair $(\ndgEpi({S^*y}),S)$ is to be regarded as a morphism from
    $\ndgRoot({S^*y})$ to $y$,
    witnessed by the equality $\ndgEpi({S^*y})^*\ndgRoot({S^*y})=S^*y$.
    Note that $\ndgEpi({S^*y})$ is an identity map
    if and only if $S^*y$ is non-degenerate.

    The given function is injective: if $S_1\neq S_2$
    then $(\ndgEpi({S_1^*y}),S_1)\neq(\ndgEpi({S_2^*y}),S_2)$.
    To check surjectivity, suppose that $(\sigma,\nu)$
    is a morphism $z\to y$ in $\Sd\,\C$.
    Writing $S=\img(\nu)$ so that $S^*y=\nu^*y$, we
    obtain the equality $S^*y=\sigma^*z$ which demonstrates
    that $z$ is the non-degenerate root of $S^*y$.
    Therefore we have $z=\ndgRoot({S^*y})$ and $\sigma=\ndgEpi({S^*y})$.
  \end{proof}
\end{prop}

\begin{lem}\label{result:2^m+1-1}
  An object $y$ of $\Sd\,\C$ is an $m$-simplex
  if and only if the set $\mt y$ has cardinality $2^{m+1}-1$.
  \begin{proof}
    Every object of $\Sd\,\C$ is a $m$-simplex for some $m$.
    There are $2^{m+1}-1$ nonempty subsets of $\{0,\ldots,m\}$,
    and just as many morphisms targeting each $m$-simplex.
  \end{proof}
\end{lem}

Thus the categorical structure of $\Sd\,\C$ encodes the
specific dimension of simplices.
Note that if $x$ is an $m$-simplex and $y$ is an $n$-simplex
such that there exists some non-identity morphism $x\to y$ in $\Sd\,\C$,
then we must have inequality $m<n$.
This is because if $(\sigma,\nu)$ is an epi-mono pair satisfying
$\sigma^*x=\nu^*y$, then we have $\dom\,\sigma=[k]=\dom\,\nu$ for
some $k$ satisfying $m\leq k\leq n$.
We obtain a strict inequality $m<n$ if either $\sigma$ or $\nu$
is a non-identity morphism.

\begin{exmp}[Automorphisms of $\Sd{\,[n]}$]
The category $[n]$ is a finite poset.
Therefore, as discussed in Example \ref{exmp:SdB=SdC,B=/=C,B=/=Cop},
there is a canonical isomorphism between the subdivision $\Sd\,[n]$ and
the poset of non-empty subsets of $\{0,\ldots,n\}$,
with order given by subset inclusion.
Under this isomorphism, each $0$-simplex $\bangle k$ of $\Sd\,[n]$ is sent to
the singleton subset $\{k\}$ of $\{0,\ldots,n\}$.
It follows that each automorphism of $\Sd[n]$ is determined by a permutation
of the set $\bangle 0,\ldots,\bangle n$ of $0$-simplices.
Therefore, we have an isomorphism
	$$\Aut\bigl(\Sd[n]\bigr)\cong S_{n+1}$$
between the group of automorphisms of
$\Sd[n]$ and the symmetric group on $n+1$ elements.
For comparison, the group $\Aut([n])$ is trivial for each natural number $n$.
\end{exmp}

\begin{defns}\label{defn:faces}
  Let $x$ and $y$ be objects in $\Sd\,\C$ such that the
  hom-set $\Sd\,\C(x,y)$ is non-empty.
  Then we say that $x$ is a \textit{face} of $y$.
  If, in addition, $x$ and $y$ are distinct,
  then $x$ is a \textit{proper face} of $y$.
  The \textit{category of faces} of $y$, written $\faces y$,
  is the full subcategory of $\Sd\,\C$ whose objects are the faces of $y$.
\end{defns}

The faces of an $m$-simplex $y$ are
precisely the simplices $\ndgRoot(S^*y)$.
This is because, by Proposition \ref{result:inj=mt},
the maps into $y$ are precisely of the form $(\ndgEpi({S^*y}),S)$ for $S$ a
non-empty subset of $\{0,\ldots,m\}$.
The proper faces of $y$ are those simplices $\ndgRoot(S^*y)$
where $S$ is a proper non-empty subset of $\{0,\ldots,m\}$.

  Here is an example.
  Let $f$, $g$, and $h$ be distinct non-identity morphisms in $\C$,
  and assume that $h=g\circ f$.
  Write $y=\barcons{g|f}$, and set $a=\dom\,f$,
  $b=\cod\,f=\dom\,g$, and $c=\cod\,g$.
  If the objects $a,b,c$ are all distinct then the category
  of faces $\faces y$ is as in Figure \ref{fig:cat_carried_by_y}. 
  If $a$ equals $c$ and is distinct from $b$, then
  $\faces y$ is as in Figure \ref{fig:cat_carried_by_y2}.
  If $a=b=c$ then $\faces y$ is as in Figure \ref{fig:cat_carried_by_y3}.
\nopagebreak
\begin{figure}[H]
\centering
\begin{minipage}[b]{.32\textwidth}
    \centering
    \begin{tikzpicture}[baseline=(current  bounding  box.center),
      scale=.7, every node/.style={transform shape}]
      \path  (0,0)      node (mid) {$         y  $}
        +(-150:2*1.3cm) node (c)   {$\bangle  c  $}
        +( -90:1*1.3cm) node (gf)  {$\barcons h  $}
        +( -30:2*1.3cm) node (a)   {$\bangle  a  $}
        +(  30:1*1.3cm) node (f)   {$\barcons g  $}
        +(  90:2*1.3cm) node (b)   {$\bangle  b  $}
        +( 150:1*1.3cm) node (g)   {$\barcons f  $}
        ;
        \draw[->] (a)  -- (gf) ;
        \draw[->] (a)  -- (mid);
        \draw[->] (a)  -- (f)  ;
        \draw[->] (b)  -- (f)  ;
        \draw[->] (b)  -- (mid);
        \draw[->] (b)  -- (g)  ;
        \draw[->] (c)  -- (g)  ;
        \draw[->] (c)  -- (mid);
        \draw[->] (c)  -- (gf) ;
        \draw[->] (f)  -- (mid);
        \draw[->] (g)  -- (mid);
        \draw[->] (gf) -- (mid);
      \end{tikzpicture}
      \caption{}\label{fig:cat_carried_by_y}
\end{minipage}
\begin{minipage}[b]{.32\textwidth}
    \centering
    \begin{tikzpicture}[baseline=(current  bounding  box.center),
      scale=.7, every node/.style={transform shape}]
      \path  (0,0)       node (mid) {$y     $}
        +(  -90:1*1.3cm) node (gf)  {$\barcons h$}
        +(  -90:2*1.3cm) node (a)   {$\bangle  a        $}
      +(   0:1.5*1.3cm) node (f)   {$\barcons g        $}
      +( 90:   1.3cm) node (b)   {$\bangle  b        $}
      +( 180:1.5*1.3cm) node (g)   {$\barcons f        $}
      ;
      \draw[->] (a. 15+90)  --                        (gf.-20-90);
      \draw[->] (a.-15+90)  --                        (gf. 20-90);
      \draw[->] (a) .. controls (-.60*1.3cm,-1*1.3cm) ..  (mid)      ;
      \draw[->] (a) .. controls ( .60*1.3cm,-1*1.3cm) ..  (mid)      ;
      \draw[->] (a.90-50)         --                        (f)        ;
      \draw[->] (a.90+50)         --                        (g)        ;
      \draw[->] (b)         --                        (f)        ;
      \draw[->] (b)         --                        (mid)      ;
      \draw[->] (b)         --                        (g)        ;
      \draw[->] (f)         --                        (mid)      ;
      \draw[->] (g)         --                        (mid)      ;
      \draw[->] (gf)        --                        (mid)      ;
    \end{tikzpicture}
      \caption{}\label{fig:cat_carried_by_y2}
\end{minipage}
\begin{minipage}[b]{.32\textwidth}
    \begin{tikzpicture}[baseline=(current  bounding  box.center),
      scale=.7, every node/.style={transform shape},
            cross line/.style={preaction={draw=white,shorten >=0.05cm, -,
               line width=4pt}}]
      \path  (0,0)           node (mid) {$ y     $}
     	+(   0:1.5*1.3cm) node (g)   {$\barcons g        $}
     	+( 180:1.5*1.3cm) node       {$\phantom{\barcons g        }$}
    	+(-135:1.5*1.3cm) node (f)   {$\barcons f        $}
    	+( 135:1.5*1.3cm) node (gf)  {$\barcons h$}
    	+( -90:2  *1.3cm) node (a)   {$\bangle  a        $}
    	;
      \draw[->] (a.113+3)             -- (f.180+113-13)  ;
      \draw[->] (a.113-3)             -- (f.180+113+13) ;
      \draw[->] (a.63+3)             -- (g.180+63-13)  ;
      \draw[->] (a.63-3)             -- (g.180+63+13) ;
      \draw[->] (a.100+2)             -- (gf.180+100-13)  ;
      \draw[->] (a.100-2)             -- (gf.180+100+13)  ;
      \draw[->] (a)             -- (mid);
      \draw[->] (a. 5+90)             -- (mid.-25-90);
      \draw[->] (a.-5+90)             -- (mid. 25-90);
      \draw[->] (a)             -- (mid);
      \draw[->,cross line] (f)  -- (mid.-130);
      \draw[->] (g)             -- (mid);
      \draw[->] (gf)            -- (mid);
    \end{tikzpicture}
      \caption{}\label{fig:cat_carried_by_y3}
\end{minipage}
\end{figure}
Note that identity arrows have been omitted from the above diagrams.

\begin{prop}\label{result:proper_faces=dom/cod}
  Let $\barcons f:[1]\to\C$ be a 1-simplex in $\Sd\,\C$.
  Then the proper faces of $\barcons f$ are $\bangle{\cod\,f}$
  and $\bangle{\dom\,f}$.
  \begin{proof}
    By Proposition \ref{result:inj=mt}, the proper faces of $\barcons f$
    are the non-degenerate roots $\ndgRoot({\{0\}^*\barcons f})$
    and $\ndgRoot({\{1\}^*\barcons f})$.
    Every 0-simplex is non-degenerate, so we have
\begin{gather*}
    \lambda\left(\{0\}^*\barcons f\right)=\{0\}^*\barcons f=\bangle{\dom\,f}\\
    \text{ and }\\
    \lambda\left(\{1\}^*\barcons f\right)=\{1\}^*\barcons f=\bangle{\cod\,f}.
  \qedhere
\end{gather*}
  \end{proof}
\end{prop}

\begin{cor}
A non-identity morphism $f$ is an endomorphism if and only if
the 1-simplex $\barcons f$ has just one proper face.
\end{cor}


\begin{wrapfigure}{r}{0.35\textwidth}
\centering
  \tikz{\path node (f) {$\barcons f$}
              node [left =5mm of f] (0) {$\bangle 0$}
              node [right=5mm of f] (1) {$\bangle 1$};
        \draw[->](0)--(f);\draw[->](1)--(f);}
\end{wrapfigure}

Note that $\Sd\,\C$ does not distinguish
which face of $\barcons f$ corresponds to domain and which to codomain.
For example, write $f$ for the arrow $0\to 1$ in $[1]$, and consider
the category $\Sd\,[1]$ pictured to the right.
Note that the non-identity automorphism of $\Sd\,[1]$ switches
$\bangle 0$ with $\bangle 1$.
Thus, $\Sd$ introduces symmetry, illustrated generally
by the canonical isomorphism $\Sd\,\C\cong\Sd(\op\C)$.
At best, we can hope for $\Sd\,\C$ to encode $\C$ ``up to opposites.''

\subsection{Encoding triangles}

We have seen that the proper faces of 1-simplices correspond
to their domain and codomain; we now go up one dimension
to see how 2-simplices in $\Sd\,\C$ codify relationships
among morphisms in $\C$.

\begin{notn}\label{notn:mt_n}
  Given an object $y$ of $\Sd\,\C$, write $\mt[n]y$ for the set
    \[ \{f\in\Mor(\Sd\,\C)~|
      ~\text{$\cod\,f=y$ and $\dom\,f$ is an $n$-simplex}\}
    \]
  of morphisms targeting $y$ that have source equal to some
  non-degenerate $n$-simplex.
\end{notn}

Some of the coming proofs will require explicit calculation of
the sets $\mt[0]y$ and $\mt[1]y$ for 2-simplices $y$.
Below we state some general facts that will make
such calculation easier.

Elements of $\mt[n]y$ are morphisms
$(\ndgEpi({S^*y}),S):\ndgRoot({S^*y})\to y$
such that $\ndgRoot({S^*y})$ is an $n$-simplex.
Note that the sets $\mt[n]y$ partition $\mt y$.
Also, note that $\ndgRoot({S^*y})$ has dimension less than or equal to $S^*y$,
which in turn has dimension equal to $|S|-1$.
Therefore, to find $\mt[n]y$ it is enough to consider only those $S$ with
$\abs S-1\geq n$.
Indeed, if $\abs S-1<n$ then
  \[\dim(\ndgRoot({S^*y}))\leq\dim(S^*y)\leq\abs S-1<n,\]
hence the dimension of $\ndgRoot({S^*y})$ is less than $n$, and
the morphism $(\ndgEpi({S^*y}),S)$ cannot be in $\mt[n]y$.
Finally, note that if $y$ is an $m$-simplex and if $S$ is equal to the
full set $\{0,\ldots,m\}$,
then $S^*y=y$ is non-degenerate, in which case we have equality
  \[\dim(\ndgRoot({S^*y}))=\dim(y)=m.\]
Therefore, to find $\mt[n]y$ in the case where $n<m$, it is good
enough to consider the proper subsets $S$ of $\{0,\ldots,m\}$.

The next order of business will be to show how $\Sd$ encodes inversion.
This will allow us to determine by looking at $\Sd\,\C$ whether
the category $\C$ is a groupoid.
Moreover, if $\G$ is a groupoid and $\barcons f$ is a 1-simplex in $\Sd\,\G$,
being able to find $\barcons{f^{-1}}$ will
help determine whether given 1-simplices
$\barcons g$ and $\barcons h$ satisfy $f\circ g=h$.

\begin{prop}\label{result:fg=id<=>|mt1<f|g>|=2}
  Let $f,g$ be non-identity morphisms in $\C$ satisfying\/ $\dom\,f=\cod\,g$.
  Set $y$ equal to the 2-simplex $\barcons{f|g}$.
  If $f$ is left-inverse to $g$ then $\mt[0]y$ has four elements
  and $\mt[1]y$ has two elements.
  On the other hand, if the composite $f\circ g$ is not an identity arrow
  then $\mt[0]y$ and $\mt[1]y$ have three elements each.
  \begin{proof}
    For any object $x$ of $\Sd\,\C$, the hom-set $\Sd\,\C(x,x)$
    contains only one morphism.
    We have the identity $\{0,1,2\}^*y=y$, therefore the map from Proposition
    \ref{result:inj=mt} restricts to a bijection
    between the set of non-empty proper subsets of $\{0,1,2\}$
    and the union $\mt[0]y\sqcup\mt[1]y$.
    By counting subsets of $\{0,1,2\}$
    we find $\abs{\mt[0]y}+\abs{\mt[1]y}=6$.

    For each singleton subset $\{i\}$ we have a morphism
    with domain $\{i\}^*y$, hence $\mt[0]y$ has at least three elements.
    The subsets $\{0,1\}$ and $\{1,2\}$ correspond to morphisms
    $\barcons g\to y$ and $\barcons f\to y$,
    thus $\mt[1]y$ has size at least two.

    It remains to consider the subset $\{0,2\}$.
    We have $\{0,2\}^*y=\barcons{f\circ g}$,
    and thus the bijection from Proposition \ref{result:inj=mt}
    sends $\{0,2\}$ to an arrow
      \[\ndgRoot({\barcons{f\circ g}})\to y\]
    in $\Sd\,\C$.
    If the composite $f\circ g$ is an identity arrow
    then $\ndgRoot({\barcons{f\circ g}})$ is a 0-simplex.
    On the other hand, if $f$ is not left-inverse to $g$
    then we have a third element of $\mt[1]y$.
  \end{proof}
\end{prop}

Given the identification of objects in $\Sd\,\C$ with simplices of $N\C$,
we may think of each 2-simplex in $\Sd\,\C$ as witnessing a
triangle-shaped commutative diagram in $\C$.

\begin{defn}\label{defn:classification_of_2splx}
\onehalfspacing
  Let $f$, $g$, and $h$ be morphisms in $\C$,
  and let $y$ be a 2-simplex of $N\C$.
  Say $y$ is of the form $\bigtri fgh$
  if $y$ represents the composite of some
  pair among $f,g,h$ yielding the third morphism.
  Explicitly, $y$ is of the form $\bigtri fgh$ if one
  of the following is satisfied:
  \begin{multicols}{2}
    \begin{itemize}
    \item[1)] $\mathmakebox{f\circ g=h~\text{ and }~y\,=\barcons{f|g}}$,
    \item[2)] $\mathmakebox{g\circ f=h~\text{ and }~y\,=\barcons{g|f}}$,
    \item[3)] $\mathmakebox{f\circ h=g~\text{ and }~y\,=\barcons{f|h}}$,
    \item[4)] $\mathmakebox{h\circ f=g~\text{ and }~y\,=\barcons{h|f}}$,
    \item[5)] $\mathmakebox{h\circ g=f~\text{ and }~y\,=\barcons{h|g}}$,
    \item[6)] $\mathmakebox{g\circ h=f~\text{ and }~y\,=\barcons{g|h}}$.
    \end{itemize}
  \end{multicols}
\end{defn}

{\onehalfspacing
The expression $\bigtri fgh$ does not encode the order of $f$, $g$, and $h$;
the forms $\bigtri fgh$ and $\bigtri ghf$ are logically equivalent.
We will sometimes label vertices of these triangles.
For example, if $y:[2]\to\C$ represents the commutative triangle
\begin{equation}\label{commutative_triangle_demo}
\begin{minipage}[m]{.5\textwidth}\centering
  \tikz{
   \path	(0,0) node (middle) {}
   	(-30:\trilength) node (a) {$\scriptstyle{a}$}
   	( 90:\trilength) node (b) {$\scriptstyle{b}$}
   	(210:\trilength) node (c) {$\scriptstyle{c}$};
   \draw[->] (b) to node[transform canvas={yshift=\shiftHeight},swap]
                                   {$\scriptstyle f$} (c);
   \draw[->] (a) to node[transform canvas={yshift=\shiftHeight},swap]
                                   {$\scriptstyle g$} (b);
   \draw[->] (a) to node[transform canvas={yshift=.03cm}]  {$\scriptstyle h$} (c);
   }.
\end{minipage}
\end{equation}
in $\C$, then we will say that $y$ is of the form $\bigtri[a][b][c]fgh$.
Note that if there is a non-degenerate 2-simplex of the form $\bigtri fgh$
then the two morphisms $y[0<1]$ and $y[1<2]$ must be non-identity,
where $[i<j]$ denotes the morphism $i\to j$ in the category $[2]$.
Therefore, if $y$ is 2-simplex in $\Sd\,\C$ of the form $\bigtri fg\id$,
then $y$ satisfies either
    \begin{itemize}
      \item[1)] $y=\barcons{f|g}$ and $f$ is left inverse to $g$, or
      \item[2)] $y=\barcons{g|f}$ and $g$ is left inverse to $f$.
    \end{itemize}
\par}

The following Lemma shows how
this triangle notation summarizes information
concerning 2-simplices in $\Sd\,\C$.

\begin{lem}\label{result:faces_of_Sd2C_and_triangular_forms}
\onehalfspacing
  Let $y:[2]\to\C$ be a non-degenerate 2-simplex.
  Then $y$ is of the form $\bigtri gfh$ for some non-identity
  arrows $f$, $g$, and $h$ if and only if there are three morphisms in the set
  $\mt[1]y$, and these morphisms have respective domains
  $\barcons f$, $\barcons g$, and $\barcons h$.
  On the other hand, $y$ is of the form $\bigtri fg\id$ if and only if
  there are two morphisms in the set $\mt[1]y$, and these morphisms
  have respective domains $\barcons f$ and $\barcons g$.
  \begin{proof}
    Let $y$ be a non-degenerate 2-simplex $[2]\to\C$.
    By definition, $y$ is of the form $\bigtri fgh$ for some
    morphisms $f$, $g$, and $h$ in $\C$.
    Assume without loss of generality that $y=\barcons{f|g}$.
    Then we have $h=g\circ f$, again by definition.
    Let $a,b,c$ satisfy
    \[
    a=\dom\,g
    \hspace{15pt}\text{and}\hspace{15pt}
    \cod\,g=b=\dom\,f
    \hspace{15pt}\text{and}\hspace{15pt}
    \cod\,f=c
    \]
    as in the commutative triangle (\ref{commutative_triangle_demo}) above.
    By Proposition \ref{result:inj=mt}, the non-identity elements of $\mt y$
    are in bijection with the proper non-empty subsets of $\{0,1,2\}$.

    Suppose first that all three of $f$, $g$, and
    $h$ are non-identity morphisms.
    Then the non-identity elements of $\mt y$ are
    the morphisms displayed below:
      \begin{align*}
        \bangle a&\xrightarrow{(\id_{[0]},\{0\})} y  & \barcons g&\xrightarrow{(\id_{[1]},\{0,1\})} y \\
        \bangle b&\xrightarrow{(\id_{[0]},\{1\})} y  & \barcons f&\xrightarrow{(\id_{[1]},\{1,2\})} y \\
        \bangle c&\xrightarrow{(\id_{[0]},\{2\})} y  & \barcons h&\xrightarrow{(\id_{[1]},\{0,2\})} y .
      \end{align*}
    On the other hand, if $f$ is left-inverse to $g$ then
    then we must have $a=c$ and $f\circ g=\id_a$.
    In this case, writing $\sigma$ for the unique epimorphism $[1]\to[0]$, the
    non-identity elements of $\mt y$ are the morphisms displayed below:
      \begin{align*}
        \bangle a&\xrightarrow{(\id_{[0]},\{0\})} y  & \barcons g&\xrightarrow{(\id,\{0,1\})} y \\
        \bangle b&\xrightarrow{(\id_{[0]},\{1\})} y  & \barcons f&\xrightarrow{(\id,\{1,2\})} y \\
        \bangle a&\xrightarrow{(\id_{[0]},\{2\})} y  & \bangle a&\xrightarrow{(\sigma,\{0,1\})} y .
      \end{align*}
    Thus we have $|\mt[0]y|=4$ and $|\mt[1]y|=2$ if $f$ is left-inverse to $g$,
    whereas $|\mt[0]y|=3=|\mt[1]y|$ if $f\circ g$ is non-identity.
  \end{proof}
\end{lem}

\begin{exmp}\label{exmp:faces_of_<g|f>-1}
  By Lemma \ref{result:2^m+1-1} we have $|\mt[0]y|+|\mt[1]y|=6$
  for any non-degenerate 2-simplex $y$.
  If $y$ is of the form
    \[ \bigtri[a][b][c]fgh
    \]
  and satisfies $\abs{\mt[0]y}=3$, then the morphisms in
  $\mt[0]y$ have respective domains $\bangle a$, $\bangle b$, and $\bangle c$,
  as in Figures \ref{fig:cat_carried_by_y} through \ref{fig:cat_carried_by_y3}.

  Suppose instead that $y$ is of the form
    \[  \bigtri[a][b][a]fg{\id_a},
    \]
  satisfying $\abs{\mt[0]y}=4$.
  Letting $\sigma$ denote the unique epimorphism $[1]\to[0]$,
  the elements of $\mt[0]y$ are given by the pairs
    \begin{align*}
      (\sigma,\{0,2\})  & : \bangle a\to y,            \\
      (\id_{[0]},\{0\}) & : \bangle a\to y,            \\
      (\id_{[0]},\{1\}) & : \bangle b\to y,\text{ and} \\
      (\id_{[0]},\{2\}) & : \bangle a\to y.            \\
    \end{align*}
  If $a$ and $b$ are distinct, then $\faces y$
  appears as in Figure \ref{fig:faces-degen-1}.
  If $a=b$ but $f\neq g$, then $\faces y$
  is as in Figure \ref{fig:faces-degen-2}.
  If $f=g$ then $\faces y$
  is as in Figure \ref{fig:faces-degen-3}.
\nopagebreak
\begin{figure}[H]
\centering
\begin{minipage}[b]{.32\textwidth}
    \centering
    \begin{tikzpicture}[baseline=(current  bounding  box.center),
      scale=.7, every node/.style={transform shape}]
      \path  (0,0)           node (mid) {$y$}
      	+( 180:1.5\dlen) node (f)   {$\barcons f        $}
    	+(   0:1.5\dlen) node (-f)  {$\barcons{g}  $}
      	+( -90:1.5\dlen) node (a)   {$\bangle  b        $}
    	+(  90:1.5\dlen) node (b)   {$\bangle  a        $}
    	;
      \draw[->] (b)  --              (f)  ;
      \draw[->] (a)  --              (f)  ;
      \draw[->] (b)  --              (-f) ;
      \draw[->] (a)  --              (-f) ;
      \draw[->] (b)  --              (mid);
      \draw[->] (b.-90+25)  --              (mid.90-25);
      \draw[->] (b.-90-25)  --              (mid.90+25);
      \draw[->] (a)  --              (mid);
      \draw[->] (f)  --              (mid);
      \draw[->] (g)  --              (mid);
    \end{tikzpicture}
      \caption{}\label{fig:faces-degen-1}
 \end{minipage}
\begin{minipage}[b]{.32\textwidth}
\centering
    \begin{tikzpicture}[baseline=(current  bounding  box.center),
      scale=.7, every node/.style={transform shape}]
      \path  (0,0)          node (mid) {$y$}
            +(  0:1.5\dlen) node (f-)  {$\barcons{g}$}
    	+(180:1.5\dlen) node (f)   {$\barcons f   $}
    	+(-90:1.5\dlen) node (a)   {$\bangle  a   $}
    	;
      \draw[->] (f)  -- (mid);
      \draw[->] (f-) -- (mid);
      \draw[->] (a.90+34) -- (mid.-90-34);
      \draw[->] (a.90+12) -- (mid.-90-12);
      \draw[->] (a.90-12) -- (mid.-90+12);
      \draw[->] (a.90-34) -- (mid.-90+34);
      \draw[->] (a.45-6)    -- (f-.180+45+6);
      \draw[->] (a.45+6)    -- (f-.180+45-6);
      \draw[->] (a.90+45-6)    -- (f.-45+6);
      \draw[->] (a.90+45+6)    -- (f.-45-6);
    \end{tikzpicture}
\caption{}
\label{fig:faces-degen-2}
\end{minipage}
\begin{minipage}[b]{.32\textwidth}
\centering
    \begin{tikzpicture}[baseline=(current  bounding  box.center),
      scale=.7, every node/.style={transform shape}]
      \path  (0,0)          node (mid) {$y$}
    	+(180:1.5\dlen) node (f)   {$\barcons{f}$}
    	+(  0:1.5\dlen) node       {$\phantom{\barcons{f}}$}
    	+(-90:1.5\dlen) node (a)   {$\bangle  a   $}
    	;
      \draw[transform canvas={yshift=0.42ex},->] (f) -- (mid);
      \draw[transform canvas={yshift=-0.42ex},->] (f) -- (mid);
      \draw[->] (a.90+34) -- (mid.-90-34);
      \draw[->] (a.90+12) -- (mid.-90-12);
      \draw[->] (a.90-12) -- (mid.-90+12);
      \draw[->] (a.90-34) -- (mid.-90+34);
      \draw[->] (a.90+45-6)    -- (f.-45+6);
      \draw[->] (a.90+45+6)    -- (f.-45-6);
    \end{tikzpicture}
\caption{}
\label{fig:faces-degen-3}
\end{minipage}
\end{figure}
\end{exmp}

{\onehalfspacing
  \begin{lem}[Inverse criterion]\label{result:inv-criterion}
    A non-identity morphism $f$ in $\C$ is self-inverse if and only if
    there exists a 2-simplex of the form $\bigtri ff\id$ in $\Sd\,\C$.
    For distinct non-identity morphisms $f$ and $g$ in $\C$,
    $f$ is the two-sided inverse to $g$ if and only if
    there exist two distinct 2-simplices of the form $\bigtri fg\id$.
    \begin{proof}
      If there is a 2-simplex of the form $\bigtri ff\id$
      then it must be equal to $\barcons{f|f}$.
      In this case $f$ is left inverse to $f$
      (by Proposition \ref{result:fg=id<=>|mt1<f|g>|=2}).
      Conversely, if $f$ is self-inverse then the 2-simplex
      $\barcons{f|f}$ is of the form $\bigtri ff\id$.
  
      Suppose now that $f$ and $g$ are distinct and that $\Sd\,\C$ contains
      two 2-simplices of the form $\bigtri fg\id$.
      These 2-simplices  must equal $\barcons{f|g}$ and $\barcons{g|f}$,
      so $f$ is left inverse to $g$ and vice versa.
      Conversely, if $f$ is the two-sided inverse to $g$
      then $\barcons{f|g}$ and $\barcons{g|f}$ are of the form $\bigtri fg\id$.
    \end{proof}
  \end{lem}
}

Thus, the subdivision $\Sd\,\C$ encodes whether
the morphisms in $\C$ are invertible.

\begin{exmp}[Automorphisms of $\Sd\,D_3$]\label{exmp:AutD3}
This example illustrates how the categorical
structure of a subdivision $\Sd\,\C$ might fail to
distinguish between the composites $f\circ g$ and $g\circ f$
of endomorphisms $f$ and $g$ in $\C$.

Let $r$ and $s$ be generators for the six-element dihedral group $D_3$.
  \[ D_3 = \bangle{ r,s~|~r^3,s^2,rsrs   }
  \]
Let $\phi$ denote the automorphism $D_3\to D_3$ that sends
$r$ to $r^2$ and $s$ to $s$.
Let $\alpha_{D_3}$ denote the map from Section \ref{sec:Lemmas} defined
by $\barcons{f_m|\cdots|f_1}\mapsto\barcons{f_1^{-1}|\cdots|f_m^{-1}}$.
Below is a comparison of the maps $\alpha_{D_3}$ and $\Sd\,\phi$ with the
composite $\alpha_{D_3}\circ(\Sd\,\phi)$.
\begin{align*}
\alpha_{D_3}:\Sd D_3&\to\Sd D_3   & \Sd\,\phi:\Sd D_3&\to\Sd D_3         & \alpha_{D_3}\circ\Sd\,\phi:D_3&\to D_3 \\
\barcons r  &\mapsto\barcons{r^{-1}} & \barcons r&\mapsto\barcons{r^{-1}}    & \barcons r&\mapsto\barcons r         \\
\barcons s  &\mapsto\barcons s    & \barcons s&\mapsto\barcons s       & \barcons s&\mapsto\barcons s         \\
\barcons{rs}&\mapsto\barcons{rs}  & \barcons{rs}&\mapsto\barcons{sr} & \barcons{rs}&\mapsto\barcons{sr}
\end{align*}
Note in particular that the composite $\alpha_{D_3}\circ\Sd\,\phi$
sends $r$ to $r$ and $s$ to $s$, but does not send $rs$ to $rs$.

Given 1-simplices $\barcons f$ and $\barcons g$ corresponding to composable
endomorphisms $f$ and $g$ in a groupoid $\G$, we will show in the following
section that the categorical structure of $\Sd\,\G$ allows to pick out
the set $\{\barcons{f\circ g},\barcons{g\circ f}\}$ of 1-simplices
corresponding to the composites of $f\circ g$ and $g\circ f$.
Of course, these composites may or may not be distinct.
\end{exmp}

\section{$\Sd\,\G$ encodes composition in $\G$ up to opposites}\label{sec:encoding2}

{\onehalfspacing
This section builds on the previous one.
As Lemma \ref{result:possible_3rd-sides} will demonstrate,
the assumption that the morphisms $f$ and $g$ are invertible
will make it easier to count 2-simplices of the form $\bigtri fgh$.
Therefore, we now restrict our
attention from small categories to small groupoids.
We will show how the categorical structure of $\Sd\,\G$ determines
which triples $(\barcons f,\barcons g,\barcons h)$
satisfy one of the equations $f\circ g=h$ and $g\circ f=h$.

  \begin{defn}\label{defn:filler,3rd_side}
    Let $f,g$ be non-identity arrows in a groupoid $\G$,
    and let $y$ be a 2-simplex in $\Sd\,\G$.
    Say that $y$ is a \textit{filler} for the triangle $\filltri fg$ if
    $y$ is of the form $\bigtri fg\id$ or of the form $\bigtri fgh$
    for some morphism $h$ in $\G$.
    Such an $h$ is called a \textit{third side} of $\filltri fg$.
  \end{defn}
}

Because a fillers of $\filltri fg$ are objects in $\Sd\G$,
any such filler $y$ must be non-degenerate.
The third sides of a given triangle are classified as follows.

\begin{lem}[Possible third sides]\label{result:possible_3rd-sides}
\onehalfspacing
  Let $f$ and $g$ be non-identity morphisms in a groupoid $\G$,
  and suppose that $h$ is a third side of the triangle $\filltri fg$.
  Then $h$ must equal one of the following six composites:
    \begin{equation}\label{list:possible_3rd-sides}
        f\circ g,      \ \ \ 
        g\circ f,      \ \ \ 
        f^{-1}\circ g,  \ \ \ 
        g\circ f^{-1},  \ \ \ 
        f\circ g^{-1},  \ \ \ 
        g^{-1}\circ f.
   \end{equation}
  \begin{proof}
    Assuming that $f$ and $g$ are invertible, the composites
    (\ref{list:possible_3rd-sides}) are the values of $h$ corresponding to each
    of the cases in Definition \ref{defn:classification_of_2splx}.
    Explicitly, any filler of $\filltri fg$ must be one of the six diagrams
    $[2]\to\G$ displayed below:
    \begin{dontbotheriftheequationisoverlong}
        \newcolumntype{C}[1]{>{\centering\let\newline\\\arraybackslash\hspace{0pt}}m{#1}}
        \begin{array}{ C{2cm} | C{2cm} | C{2cm} | C{2cm} | C{2cm} | C{2cm} }
              \commtri fg{fg     }[->][->][->]         &
              \commtri fg{gf     }[<-][<-][<-]         &
              \commtri fg{f^{-1}g}[<-][->][->]         &
              \commtri fg{gf^{-1}}[->][<-][<-]         &
              \commtri fg{fg^{-1}}[->][<-][->]         &
              \commtri fg{g^{-1}f}[<-][->][<-]          \\&&&&&\\
               $ \barcons{f|g}       $&
               $ \barcons{g|f}       $&
               $ \barcons{f|f^{-1}g} $&
               $ \barcons{gf^{-1}|f} $&
               $ \barcons{fg^{-1}|g} $&
               $ \barcons{g|g^{-1}f} $
        \end{array}
    \end{dontbotheriftheequationisoverlong}
  \end{proof}
\end{lem}

Note that some or all of the formal composites
(\ref{list:possible_3rd-sides}) may be undefined, depending
on how the domain and codomain of $f$ and $g$ match up.

\subsection{Composites in Groupoids}

Suppose that $f$ and $g$ are endomorphisms of an object in a groupoid $\G$,
and that the composites $f\circ g$ and $g\circ f$ are distinct.
The categorical structure of $\Sd\,\G$ need not distinguish between
the 1-simplices $\barcons{f\circ g}$ and $\barcons{g\circ f}$,
as in Example \ref{exmp:AutD3} on page \pageref{exmp:AutD3}.
It is possible, however, to pick the 1-simplices $\barcons{f\circ g}$
and $\barcons{g\circ f}$ out from among the other 1-simplices in $\Sd\,\G$.
Generally, given non-identity morphisms $f$ and $g$ in $\G$ satisfying
$\dom\,f=\cod\,g$ or $\dom\,g=\cod\,f$, the local structure of $\Sd\,\G$ near
$\barcons f$ and $\barcons g$ determines whether a given 1-simplex $\barcons h$
satisfies $h=f\circ g$ or $h=g\circ f$.
This will be key in proving functorality of
the map $\psi$ mentioned in the introduction.

To achieve this result, we introduce some terminology
concerning relationships between arrows in $\G$.
Each definition below can be formulated in terms of the proper faces
$\{\bangle{\dom\,f},\bangle{\cod\,f}\}$ of 1-simplices $\barcons f$.

\begin{defns}\label{defns:ete,estes,eteo,eoteo,ur}\hfill
\begin{enumerate}
  \item \textit{End-to-end morphisms}:
       $\tikz[baseline=(a.base),every loop/.style={looseness=5}]{
          \path node [inner sep=1pt] (a) {$\cdot$}
                node [inner sep=1pt, right=.75cm of a] (b) {$\cdot$}
                node [inner sep=1pt, right=.75cm of b] (c) {$\cdot$};
          \draw[->] (a) -- node {} (b);
          \draw[->] (b) -- node {} (c);
             }$
     or $\tikz[baseline=(a.base),every loop/.style={looseness=5}]{
            \path node [inner sep=1pt] (a) {$\cdot$}
                  node [inner sep=1pt, right=.75cm of a] (b) {$\cdot$}
                  node [inner sep=1pt, right=.75cm of b] (c) {$\cdot$};
            \draw[->] (a) -- node {} (b);
            \draw[<-] (b) -- node {} (c);
               }$
     or $\tikz[baseline=(a.base),every loop/.style={looseness=5}]{
            \path node [inner sep=1pt] (a) {$\cdot$}
                  node [inner sep=1pt, right=.75cm of a] (b) {$\cdot$}
                  node [inner sep=1pt, right=.75cm of b] (c) {$\cdot$};
            \draw[<-] (a) -- node {} (b);
            \draw[->] (b) -- node {} (c); }$\\
      Morphisms $f$ and $g$ are end-to-end if neither $f$ nor $g$ is
      an endomorphism and the intersection
      $\{\dom\,f,\cod\,f\}\cap\{\dom\,g,\cod\,g\}$ has one element.
      Note that this implies the three dots are necessarily distinct.\\
  \item \textit{Ends-to-ends morphisms}:
       $\tikz[baseline=(a.base),every loop/.style={looseness=5}]{
        \path node [inner sep=1pt] (a) {$\cdot$}
              node [inner sep=1pt, right=.75cm of a] (b) {$\cdot$};
        \draw[->,transform canvas={yshift=0.4ex}] (a) -- node {} (b);
        \draw[->,transform canvas={yshift=-0.4ex}] (b) -- node {} (a); }$
    or
       $\tikz[baseline=(a.base),every loop/.style={looseness=5}]{
         \path node [inner sep=1pt] (a) {$\cdot$}
               node [inner sep=1pt, right=.75cm of a] (b) {$\cdot$};
         \draw[->,transform canvas={yshift=0.4ex}] (a) -- node {} (b);
         \draw[<-,transform canvas={yshift=-0.4ex}] (b) -- node {} (a);
          }$\\
      Morphisms $f$ and $g$ are ends-to-ends if neither $f$ nor $g$ is an
      endomorphism and there is equality
      $\{\dom\,f,\cod\,f\}=\{\dom\,g,\cod\,g\}$.
      This means that the intersection
      $\{\dom\,f,\cod\,f\}\cap\{\dom\,g,\cod\,g\}$ has two elements.\\
  \item \textit{End-to-endo morphisms}:
       $\tikz[baseline=(a.base),every loop/.style={looseness=5}]{
          \path node [inner sep=1pt] (a) {$\cdot$}
                node [inner sep=1pt, right=.75cm of a] (b) {$\cdot$};
          \draw[->] (a) -- node {} (b);
          \draw[->] (b) to [out=35,in=-35,loop] (b);
             }$
     or $\tikz[baseline=(a.base),every loop/.style={looseness=5}]{
        \path node [inner sep=1pt] (a) {$\cdot$}
              node [inner sep=1pt, right=.75cm of a] (b) {$\cdot$};
        \draw[<-] (a) -- node {} (b);
        \draw[->] (b) to [out=35,in=-35,loop] (b);
           }$\\
      Morphisms $f$ and $g$ are end-to-endo if
      $f$ is not an endomorphism, $g$ is an endomorphism, and the intersection
          $\{\dom\,f,\cod\,f\}\cap\{\dom\,g,\cod\,g\}$ has one element.\\
  \item \textit{Endo-to-endo morphisms}:
       $\tikz[baseline=(a.base),every loop/.style={looseness=5}]{
            \path node [inner sep=1pt] (a) {$\cdot$};
            \draw[->] (a) to [out=35,in=-35,loop] (a);
            \draw[->] (a) to [out=215,in=145,loop] (a);
               }$\\
      Morphisms $f$ and $g$ are endo-to-endo if they
      are both endomorphisms of a common object.\\
  \item \textit{Unrelated morphisms}\\
      Morphisms $f$ and $g$ are unrelated if the sets $\{\dom\,f,\cod\,f\}$
      and $\{\dom\,g,\cod\,g\}$ have no elements in common.
\end{enumerate}
\end{defns}

For example, non-endomorphisms $f$ and $g$ in $\G$ are end-to-end
if and only if $\barcons f$ and $\barcons g$ have one face in common.
A non-identity morphism $f$ is an endomorphism if and only if the
1-simplex $\barcons f$ has exactly one proper face.

\begin{defn}\label{defn:sequential,coinitial,coterminal}
Let $f$ and $g$ be end-to-end morphisms.
We say that $f$ and $g$ are \textit{sequential} if $\dom\,f=\cod\,g$ or
$\dom\,g=\cod\,f$.
We say that $f$ and $g$ are \textit{coinitial} if $\dom\,f=\dom\,g$,
and that $f$ and $g$ are \textit{coterminal} if $\cod\,f=\cod\,g$.
\end{defn}

\begin{defn}\label{defn:parallel,opposed}
Let $f$ and $g$ be ends-to-ends morphisms.
We say that $f$ and $g$ are \textit{parallel} if $\dom\,f=\dom\,g$ and
$\cod\,f=\cod\,g$.
We say that $f$ and $g$ are \textit{opposed} if $\dom\,f=\cod\,g$ and
$\dom\,g=\cod\,f$.
\end{defn}

Observe that end-to-end morphisms are sequential if and only if
they can be composed in some order.
Similarly, ends-to-ends morphisms are opposed if and only if
they can be composed in either order.
Note that unrelated morphisms are never composable, and that
end-to-endo and endo-to-endo morphisms are always composable.
Below are criteria for the composability of
end-to-end and ends-to-ends morphisms.

\begin{prop}\label{result:composability-criterion-for-end-to-end}
{\onehalfspacing
  Let $f$ and $g$ be end-to-end morphisms in $\G$.
  There exists a unique filler for the triangle $\filltri fg$
  if and only if $f$ and $g$ are sequential.
  There is more than one filler for $\filltri fg$ if and only if
  $f$ and $g$ are coinitial or coterminal.\par}
  \singlespacing
    \begin{proof}
    It will suffice to show that if $f$ and $g$ are coinitial or coterminal
    then there are exactly two fillers for $\filltri fg$, and that if $f$ and
    $g$ are sequential, then there is exactly one filler for $\filltri fg$.

    Given distinct elements $i$ and $j$ of the set $\{0,1,2\}$, we will
    write $[i<j]$ for the morphism from $i\to j$ in the category $[2]$.
    If $y$ is a 2-simplex in $\Sd\G$, i.e. a functor $[2]\to\G$
    that sends $[0<1]$ and $[1<2]$ to non-identity arrows, then:
    \begin{enumerate}
      \item the domain of $y[0<1]$ equals the domain of $y[0<2]$,
      \item the codomain of $y[1<2]$ equals the codomain of $y[0<2]$, and
      \item $y[1<2]$ and $y[0<1]$ are composable.
    \end{enumerate}

    Given a 2-simplex $y$ in $\Sd\G$ that fills $\filltri fg$,
    we must have $f=y[i<j]$ and $g=y[k<l]$ for some distinct
    morphisms $[i<j]$ and $[k<l]$ in $[2]$.
    Functors preserve domain and codomain;
    by looking at the source and target of $f$ and $g$,
    we can rule out combinations of $i,j,k,l$.

    Suppose first that $f$ and $g$ are coinitial.
    Then any filler $y$ of $\filltri fg$
    must satisfy either
    \begin{itemize}
      \item $f=y[0<1]$ and $g=y[0<2]$, or
      \item $f=y[0<2]$ and $g=y[0<1]$.
    \end{itemize}
    To see this, note that $y$ sends coinitial pairs in $[2]$ to coinitial
    pairs in $\G$, and similarly for sequential and coterminal pairs.
    Write $f=y[i<j]$ and $g=y[k<l]$, and note that
    $f$ and $g$ are neither sequential nor coterminal.
    By contraposition, $[i<j]$ and
    $[k<l]$ are neither sequential nor coterminal.
    Thus the preimages of $f$ and $g$ must be $[0<1]$ and $[0<2]$.

    If $f=y[0<1]$ and $g=y[0<2]$ then we must have $(y[1<2])\circ f=g$,
    hence $y$ is the 2-simplex $\barcons{gf^{-1}|f}$.
    If $f=y[0<2]$ and $g=y[0<1]$ then we have $(y[1<2])\circ g=f$,
    hence $y$ is the 2-simplex $\barcons{fg^{-1}|g}$.
    Thus, there are exactly two fillers for $\filltri fg$.

    The argument is similar supposing that $f$ and $g$ are coterminal.
    Of all end-to-end morphism pairs in the category [2], only
    $[1<2]$ and $[0<2]$ are neither coinitial nor sequential.
    Therefore, if $y$ fills $\filltri fg$
    then we must have $f=y[1<2]$ and $g=y[0<2]$, or vice versa.
    These two possibilities correspond to the cases
    $y=\barcons{f|f^{-1}g}$ and $y=\barcons{g|g^{-1}f}$.
    The coterminal case is dual to the coinitial case
    in a sense made precise by the isomorphism $\Sd\G\to\Sd(\op\G)$
    from Lemma \ref{result:SdC=SdCop},
    which sends each $n$-simplex $\barcons{f_1|\cdots|f_n}$ in $\Sd\G$
    to the $n$-simplex $\barcons{f_n|\cdots|f_1}$ in $\Sd(\op\G)$.

    Finally, suppose that $f$ and $g$ are sequential.
    If $\dom\,f=\cod\,g$ and if $y$ fills $\filltri fg$ then we have
    $f=y[1<2]$ and $g=y[0<1]$, hence $y$ equals $\barcons{f|g}$.
    Similarly, if $\dom\,g=\cod\,f$ and if $y$ fills $\filltri fg$ then we have
    $g=y[1<2]$ and $f=y[0<1]$, hence $y$ equals $\barcons{g|f}$.
    In either case, there is only one possible filler $y$.
    \end{proof}
\end{prop}

Note that the above result can fail if $f$ and $g$ are not invertible.

\begin{cor}\label{result:encoding_sequential_fg=h}
\onehalfspacing
Let $f$ and $g$ be end-to-end morphisms in $\G$.
Then $f$ and $g$ are sequential if and only if there is a unique
third side of the triangle $\filltri fg$.
This third side is necessarily equal to the composite of $f$ and $g$.
  \begin{proof}
   Suppose first that $f$ and $g$ are sequential.
   The previous proposition shows that there is a unique 2-simplex filler,
   and thus a unique third side, for the given triangle $\filltri fg$.

   For the reverse implication, suppose that $f$ and $g$ are not sequential;
   we will show that there are two distinct third sides for the given triangle.
   Note that non-sequential end-to-end morphisms must be either
   coinitial or coterminal.

   Supposing first that $f$ and $g$ are coinitial, we have two
   2-simplex fillers for $\filltri fg$, namely $\barcons{gf^{-1}|f}$ and
   $\barcons{fg^{-1}|g}$.
   The corresponding third sides are $gf^{-1}$ and $fg^{-1}$, which
   must be distinct because $\cod\,f\neq\cod\,g$.

   Similarly, if $f$ and $g$ are coterminal then we have third sides
   $g^{-1}f$ and $f^{-1}g$ of $\filltri fg$. These third sides must
   be distinct because $\dom\,g\neq\dom\,f$.
  \end{proof}
\end{cor}


The following result is analogous to the previous
proposition, concerning ends-to-ends morphisms.

\begin{prop}\label{result:composability-criterion-for-ends-to-ends}
 \onehalfspacing
  Let $f$ and $g$ be distinct ends-to-ends morphisms in $\G$.
  There are four fillers for $\filltri fg$ if and only if
  $f$ and $g$ are parallel.
  There are two fillers for the triangle $\filltri fg$
  if and only if $f$ and $g$ are opposed.
        \begin{proof}
    It will suffice to show that if $f$ and $g$ are parallel then there
    are exactly four fillers for $\filltri fg$, and that if $f$ and $g$
    are opposed, then there are exactly two filler for $\filltri fg$.
    As in the previous proof, we will write $[i<j]$ for the
    morphism $i\to j$ in $[2]$.

    Suppose first that $f$ and $g$ are parallel,
    and that $y$ fills $\filltri fg$.
    Then we must have $f=y[i<j]$ and $g=y[k<l]$ for
    some distinct morphisms $[i<j]$ and $[k<l]$ in $[2]$.
    Moreover, the morphisms $[i<j]$ and $[k<l]$ cannot be composable
    because $f$ and $g$ are not composable.
    This means that $i\neq l$ and $j\neq k$.
    Therefore, if $y$ fills $\filltri fg$ then there are four possibilities:
    \begin{enumerate}
      \item $f=y[0<1]$ and $g=y[0<2]$,
      \item $g=y[0<1]$ and $f=y[0<2]$,
      \item $f=y[0<2]$ and $g=y[1<2]$, or
      \item $g=y[0<2]$ and $f=y[1<2]$.
    \end{enumerate}
    Thus, there are at most four fillers.

    The cases (1)-(4) above correspond (respectively) to the 2-simplices
    \[
    y_1=\barcons{gf^{-1}|f},              \hspace{15pt}
    y_2=\barcons{fg^{-1}|g},              \hspace{15pt}
    y_3=\barcons{f|f^{-1}g},               \hspace{15pt}
    y_4=\barcons{g|g^{-1}f}.
    \]
    To prove that there are exactly four fillers for $\filltri fg$,
    we must show that the above fillers $y_i$ are all distinct.
    By looking at the values $y(0)$, $y(1)$, $y(2)$, we see
    that the only pairs among the fillers $y_i$ that could be equal are
    $y_1,y_2$ and $y_3,y_4$.
    We have $y_1\neq y_2$ because
      \[y_1[0<1]=f\neq g=y_2[0<1],\]
    and similarly we have $y_3\neq y_4$ because $y_3[1<2]\neq y_4[1<2]$.
    Thus, there are exactly four fillers for $\filltri fg$.

    Suppose now that $f$ and $g$ are opposed.
    The 2-simplices $\barcons{f|g}$ and $\barcons{g|f}$
    are of the form $\bigtri fg{fg}$ and $\bigtri gf{gf}$, respectively.
    We must show that these are these are the
    \textit{only} two 2-simplices that fill $\filltri fg$.
    By the proof of Lemma \ref{result:possible_3rd-sides},
    which lists all potential fillers of $\filltri fg$,
    any filler which is not equal to $\barcons{f|g}$
    or $\barcons{g|f}$ must be equal to one of the four fillers
    \[
    \barcons{gf^{-1}|f},              \hspace{15pt}
    \barcons{fg^{-1}|g},              \hspace{15pt}
    \barcons{f|f^{-1}g},               \hspace{15pt}
    \barcons{g|g^{-1}f}.
    \]
    But $f$ and $g$ are opposed, so we have
    $\dom\,f\neq\dom\,g$ and $\cod\,f\neq\cod\,g$.
    Therefore none of the composites $gf^{-1}$, $fg^{-1}$,
    $f^{-1}g$, and $g^{-1}f$ are valid,
    hence none of the four fillers above are defined.
    It follows that there are exactly two fillers for $\filltri fg$.
    \end{proof}
\end{prop}                                                                     

\begin{cor}\label{result:encoding_opposed_fg=h}
\onehalfspacing
If $f$ and $g$ are opposed ends-to-ends morphisms in $\G$,
then there are exactly two third sides of $\filltri fg$,
namely $f\circ g$ and $g\circ f$.
\end{cor}

\begin{prop}\label{result:comp->o/<-o}
  Let $f$ and $g$ be end-to-endo morphisms in $\G$.
  Then $h$ is the composite
  of $f$ and $g$ if and only if
    \begin{enumerate}
      \item $h$ and $f$ are parallel ends-to-ends morphisms, and
      \item $h$ is a third side of the triangle
            $\shifttext{3mm}{\filltri {f^{-1}}g}$.
    \end{enumerate}
\end{prop}
{\onehalfspacing
        \begin{window}[0,r,\raisebox{-30pt}{\fbox{
            \tikz[baseline=(a.base),every loop/.style={looseness=5}]{
              \path node [inner sep=1pt] (a) {$\cdot$}
                    node [inner sep=1pt, right=.75cm of a] (b) {$\cdot$};
              \draw[->] (a) -- node {$f$} (b);
              \draw[->] (b) to [out=35,in=-35,loop] node {$g$} (b);
     } } },{}]
  \textit{Proof.}
  There are two cases:
  either $\cod\,f=\dom\,g=\cod\,g$,
  or $\dom\,f=\dom\,g=\cod\,g$.
  Assume first that $\cod\,f=\dom\,g=\cod\,g$, as in the diagram to the right.
  \end{window}
  
  One implication is clear:
  the composite $g\circ f$ is parallel to $f$ because $g$
  is an endomorphism, and if we set $h=g\circ f$ then
  the 2-simplex $\barcons{h|f^{-1}}$
  is witness to $h$ being a third side of $\shifttext{3mm}{\filltri {f^{-1}}g}$.

  For the reverse implication, suppose that $h$ is a third side
  of $\shifttext{3mm}{\filltri {f^{-1}}g}$ and that $h$ is parallel to $f$.
  Then $h$ and $f^{-1}$ are opposed ends-to-ends morphisms.
  By the definition of third sides, there must exist some
  non-degenerate 2-simplex of the form
  $\shifttext{3mm}{\bigtri {f^{-1}}gh}$ in $\Sd\G$.
  Thus, $g$ is a third side of the triangle $\shifttext{3mm}{\filltri {f^{-1}}h}$,
  and it follows from Corollary \ref{result:encoding_opposed_fg=h}
  above that $g$ is equal to $f^{-1}\circ h$
  or to $h\circ f^{-1}$.
  But $g$ is an endomorphism of the object $\cod\,f$,
  whereas $f^{-1}\circ h$ is an endomorphism of $\dom\,f$.
  Because $\dom\,f\neq\cod\,f$,
  we must have $g=h\circ f^{-1}$ and therefore $g\circ f=h$.

  If we suppose instead that $\dom\,f=\dom\,g=\cod\,g$,
  then the proof follows by a similar argument.
\qed
}\\

The remainder of this section establishes a result analogous to Propositions
\ref{result:composability-criterion-for-end-to-end}-\ref{result:comp->o/<-o},
pertaining to the case where $f$ and $g$ are endo-to-endo.
For the time being we will drop the composition symbol,
writing (for example) $fg$ for the composite $f\circ g$
and $f^2$ for the composite $f\circ f$.
If $f$ is self-inverse, then we will write $f^2=\id$.
The lemma below concerns composites $f^2$.

\begin{lem}[Square criterion]\label{result:composite_g^2}
  \onehalfspacing
  Let $f$ and $h$ be non-identity
  endomorphisms in $\G$.
  Then $f^2=h$ if and only if $\Sd\G$ contains a
  non-degenerate 2-simplex of the form $\bigtri ffh$.
  \begin{proof}
    If $f^2=h$ then the 2-simplex $\barcons{f|f}$ is of the form
    $\bigtri ffh$.
    Conversely, suppose that $y:[2]\to\G$ is of the form $\bigtri ffh$.
    By Definition \ref{defn:classification_of_2splx},
    which lists all possible 2-simplices of the form $\bigtri ffh$,
    we must have $f^2=h$ or $fh=f$ or $hf=f$.
    If $fh=f$ or if $hf=f$ then we have $h=\id$, which contradicts
    our assumption that $h$ is non-identity.
  \end{proof}
\end{lem}

\begin{figure}[H]
    \centering
    \begin{tikzpicture}[baseline=(current  bounding  box.center),
      scale=.7, every node/.style={transform shape}]
      \path  (0,0)          node (mid) {$y$}
            +(  0:1.5\dlen) node (f2)  {$\barcons{f^2}$}
    	+(180:1.5\dlen) node (f)   {$\barcons f   $}
    	+(-90:1.5\dlen) node (a)   {$\bangle  a   $}
    	;
      \draw[transform canvas={yshift=0.42ex},->] (f) -- (mid);
      \draw[transform canvas={yshift=-0.42ex},->] (f) -- (mid);
      \draw[->] (f2) -- (mid);
      \draw[->] (a.90+26) -- (mid.-90-40);
      \draw[->] (a)       -- (mid);
      \draw[->] (a.90-26) -- (mid.-90+40);
      \draw[transform canvas={yshift=0.30ex,xshift=-0.30ex},->] (a) -- (f2);
      \draw[transform canvas={yshift=-0.30ex,xshift=0.30ex},->] (a) -- (f2);
      \draw[transform canvas={yshift=0.30ex,xshift=0.30ex},->] (a) -- (f);
      \draw[transform canvas={yshift=-0.30ex,xshift=-0.30ex},->] (a) -- (f);
    \end{tikzpicture}
  \caption{Let $f:a\to a$ be some endomorphism that is
  not self-inverse, and set $y=\protect\barcons{f|f}$.
  Then $\faces y$ is given by the diagram above.
  }
\label{fig:f^2=h}
\end{figure}

\begin{rem}
Given a non-degenerate 2-simplex $y:[2]\to\G$,
the category of faces $\faces y$ is given by one of the Figures
\ref{fig:cat_carried_by_y} through \ref{fig:f^2=h}.
This is to say, these seven Figures classify the possible
categories of faces of 2-simplices.
This is true for groupoids, but not for arbitrary categories.
For example, in an arbitrary category we might have $f^2=f$
for some non-identity endomorphism $f$.
\end{rem}

\begin{notn}
 \onehalfspacing
  Given morphisms $f$, $g$, and $h$ in $\G$, we will write
  $\exists\bigtri fgh$ if there exists a non-degenerate
  2-simplex of the form $\bigtri fgh$ in $\Sd\,\G$.
  The notation $\exists_n\bigtri fgh$ means that there are exactly
  $n$ distinct non-degenerate 2-simplices of that form.
  Similarly, $\nexists\bigtri fgh$
  means that no such 2-simplices exist,
  and $\exists_{\geq n}\bigtri fgh$ means that there are at least $n$
  such 2-simplices.
\end{notn}

\begin{lem}\label{lem:simultinaety}
  \onehalfspacing
  Suppose $f$, $g$, $h$, and $h'$ are endomorphisms in $\G$,
  the morphisms $f$, $g$, and $h$ are distinct, and
  a 2-simplex $y$ is simultaneously of the form $\bigtri fgh$ and $\bigtri fg{h'}$.
  Then $h=h'$.
  \begin{proof}
    The morphism $h'$ is one of the $y[i<j]$ for some $0\leq i <j\leq 1$.
    Therefore $h'$ equals $f$ or $g$ or $h$.
    If $h'$ equals $f$ or $g$, then for two distinct pairs
    $i<j$, $k<l$ we have $y[i<j]=y[k<l]$
    since $y$ is of the form $\bigtri fg{h'}$.
    On the other hand, since $y$ is of the form $\bigtri fgh$
    for $f$, $g$, $h$ all distinct, this is impossible. Therefore $h=h'$.
  \end{proof}
\end{lem}

\begin{lem}
  \doublespacing
  Let $f$ and $h$ be non-identity endomorphisms in $\G$,
  and suppose that $f\neq f^{-1}$ and $f^2\neq f^{-1}$.
  If $h=f^{-1}$, then $h=f^3$ if and only if $\exists_4\bigtri f{f^2}h$.
  If $h\neq f^{-1}$, then $h=f^3$ if and only if $\exists_2\bigtri f{f^2}h$.
  \begin{proof}
  Because $f$ is non-identity and $f\neq f^{-1}$ and $f^2\neq f^{-1}$,
  the morphisms
    \[\id, \hspace{15pt} f, \hspace{15pt} f^2, \hspace{15pt}\text{and}\hspace{15pt} f^3\]
  are all distinct.
  The cases $h= f^{-1}$ and $h\neq f^{-1}$ above correspond to whether
  or not $f^4$ equals $\id$.
    There are five distinct fillers of $\filltri f{f^2}$, displayed below:
    \vspace{-10pt}
      \[  \barcons{f|f^2}       \ \ \ \  
          \barcons{f^2|f}       \ \ \ \ 
          \barcons{f|f}         \ \ \ \ 
          \barcons{f^{-1}|f^2}  \ \ \ \ 
          \barcons{f^2|f^{-1}}  .
      \]
    These fillers are of the form $\bigtri f{f^2}{f^3}$,
    $\bigtri f{f^2}{f^3}$, $\bigtri f{f^2}f$,
    $\bigtri f{f^2}{f^{-1}}$, and $\bigtri f{f^2}{f^{-1}}$ (respectively).

    If $f^3=h=f^{-1}$, then there are four distinct 2-simplices
    of the form $\bigtri f{f^2}h$, namely
      $  \barcons{f|f^2}$,
          $\barcons{f^2|f}$,    
          $\barcons{f^{-1}|f^2}$, and
          $\barcons{f^2|f^{-1}}$.
    Conversely, if $h=f^{-1}$ and there are exactly four distinct
    2-simplices of the form $\bigtri f{f^2}h$,
    then at least one of these 2-simplices must be of the form $\bigtri f{f^2}{f^3}$.
    This 2-simplex is simultaneously of the form $\bigtri f{f^2}h$ and $\bigtri f{f^2}{f^3}$,
    so it follows by the previous Lemma that $h=f^3$.

    On the other hand, if $h=f^3\neq f^{-1}$ then there are two distinct 2-simplices
    of the form $\bigtri f{f^2}h$, namely
     $  \barcons{f|f^2}$ and
          $\barcons{f^2|f}$.  
    Conversely, if $h\neq f^{-1}$ and there are exactly two distinct
    2-simplices of the form $\bigtri f{f^2}h$,
    then these 2-simplices must be equal to $\barcons{f|f^2}$
    and $\barcons{f^2|f}$ (for otherwise the previous Lemma
    would give $f=h=f^3$ or $h=f^{-1}$, contradicting our assumptions).
    The 2-simplices $\barcons{f|f^2}$ and $\barcons{f^2|f}$
    are simultaneously of the form $\bigtri f{f^2}h$ and $\bigtri f{f^2}{f^3}$,
    and it follows from the previous lemma that $h$ equals $h^3$.
  \end{proof}
\end{lem}

The corollary below follows directly from the lemma above.
\begin{cor}[Cube criterion]\label{result:composite_g^3}
  \onehalfspacing
  Let $f$ and $h$ be non-identity endomorphisms in $\G$,
  and suppose that $f\neq f^{-1}$ and $f^2\neq f^{-1}$.
  Then $h=f^3$ if and only if
  \vspace{5pt}
    \begin{center}
            either $h=f^{-1}$ and $\exists_4\bigtri f{f^2}h$, or
            $h\neq f^{-1}$ and $\exists_2\bigtri f{f^2}h$.
    \end{center}
\end{cor}

{\onehalfspacing
Given a 1-simplex $\barcons f$ in $\Sd\,\G$ such that $f$ is
an endomorphism satisfying
$f^2\neq\id$ and $f^3\neq\id$, the previous results can be used
(for example) to find the 1-simplices $\barcons{f^2}$ and $\barcons{f^3}$.

Recall from Lemma \ref{result:possible_3rd-sides}
that if a given 2-simplex is of the form $\bigtri fgh$,
then $h$ must be one of the composites
    \begin{equation}\label{possible_3rd_sides_number_2}
        f\circ g,      \ \ \ 
        g\circ f,      \ \ \ 
        f^{-1}\circ g,  \ \ \ 
        g\circ f^{-1},  \ \ \ 
        f\circ g^{-1},  \ \ \ 
        g^{-1}\circ f.
   \end{equation}
We shall now define notation to set the stage for Lemma
\ref{result:end-to-end:multiplicity=composites}, which will state that
under certain conditions on $f$ and $g$,
the number of 2-simplices of the form $\bigtri fgh$ is equal to the number
of composites above that are equal to $h$.
\par}

\begin{notn}
\onehalfspacing
  Let $f$ and $g$ be non-identity endomorphisms
  of some object in $\G$.
  Suppose that $f\neq g$
  Write $C(f,g)$ for the set of quadruples $(k,s,l,t)$ where
    \begin{enumerate}
      \item $(k,l)$ is equal to $(f,g)$ or $(g,f)$, and
      \item $(s,t)$ is equal to $(1,1)$ or $(1,-1)$ or $(-1,1)$.
    \end{enumerate}
\end{notn}

  The six elements of $C(f,g)$ are all distinct,
  whereas the six composites (\ref{possible_3rd_sides_number_2})
  might not all be distinct.
  The evaluation map
  $\ev:C(f,g)\to\Mor(\G)$, defined by $\ev(k,s,l,t)=k^sl^t$,
  gives a correspondence between the elements of $C(f,g)$
  and the composites (\ref{possible_3rd_sides_number_2}).
  Given this correspondence, we can think of $C(f,g)$
  as a set of ``formal composites''.
  We will later make use of the sets $C(f,g)$ to keep track of 

  {\onehalfspacing
  The following Lemma shows that, if we assume $f$ and $g$ are endomorphisms
  satisfying $f^2\neq g$ and $f\neq g^2$,
  there is a bijection between $C(f,g)$ and the set of 2-simplex fillers
  for $\filltri fg$, sending each formal composite $(k,s,l,t)$ to
  to a filler whose third side equals $k^sl^t$.
  We will later use the set $C(f,g)$ to simplify a
  counting argument that involves keeping track of the
  relationship between 2-simplex fillers and third sides.
  \par}

\begin{lem}\label{result:end-to-end:multiplicity=composites}
\onehalfspacing
  Let $f$ and $g$ be non-identity endomorphisms
  of some object in $\G$.
  Suppose that $f\neq g$ and $f^2\neq g$ and $f\neq g^2$.
  For any morphism $h$ in $\G$, the number of quadruples $(k,s,l,t)\in C(f,g)$
  satisfying $k^s\circ l^t=h$ is equal to the number
  of 2-simplices of the form $\bigtri fgh$.
  \begin{proof}
Note that for any non-identity $f$ and $g$ satisfying $f^2\neq g$
and $f\neq g^2$ as above, if $h$ is a third side for
$\filltri fg$ then $h$ must be distinct from $f$ and $g$.
\par    Every 2-simplex filler of $\filltri fg$ must be one of
    the six 2-simplices
                \begin{equation*}
                \barcons{f|g}       ,\hspace{5pt}
                \barcons{g|f}       ,\hspace{5pt}
                \barcons{f|f^{-1}g} ,\hspace{5pt}
                \barcons{gf^{-1}|f} ,\hspace{5pt}
                \barcons{fg^{-1}|g} ,\hspace{5pt}
                \barcons{g|g^{-1}f} 
                \end{equation*}
    displayed in the proof of Lemma
    \ref{result:possible_3rd-sides}.
    Given the present assumptions on $f$ and $g$,
    these 2-simplices are all distinct.
    For example, $\barcons{f|g}$ is distinct from $\barcons{f|f^{-1}g}$
    because $g^{-1}$ is a non-identity morphism,
    and $\barcons{f|f^{-1}g}$ is distinct from $\barcons{gf^{-1}|f}$
    because $f^2\neq g$.
    Thus we have a bijection
      \begin{align*}
        (f,1,g,1)&\mapsto\barcons{f|g}
      & (g,1,f,1)&\mapsto\barcons{g|f} \\
        (f,-1,g,1)&\mapsto\barcons{f|f^{-1}g}
      & (g,-1,f,1)&\mapsto\barcons{g|g^{-1}f} \\
        (f,1,g,-1)&\mapsto\barcons{fg^{-1}|g}
      & (g,1,f,-1)&\mapsto\barcons{gf^{-1}|f} 
      \end{align*}
    between $C(f,g)$ and the set of 2-simplex fillers for $\filltri fg$.
    We will let
      \[\zeta:C(f,g)\to\{\text{2-simplex fillers for $\filltri fg$}\}\]
    denote this bijection.
    Note that $\zeta$ sends each quadruple $(k,s,l,t)$
    to a 2-simplex of the form
    {\setlength{\belowdisplayskip}{12pt}
    \[\tikz[baseline=(baase), auto, every node/.style={transform shape}]{
      \path (0,-.2cm) node (baase) {};
      \path	(0,0) node (middle) {}
      	(-30:\trilength) node (a) {}
      	( 90:\trilength) node (b) {}
      	(210:\trilength) node (c) {};
      \draw (b) to node[transform canvas={yshift=\shiftHeight},swap] {$f$} (c);
      \draw (a) to node[transform canvas={yshift=\shiftHeight},swap] {$g$} (b);
      \draw (a) to node[transform canvas={yshift=.03cm}]  {$k^s\circ l^t$} (c);
    }.\]}
\par    The map $\zeta$ restricts to an injection
      \[
      \zeta|_h:\{\gamma\in C(f,g)~|~\ev(\gamma)=h\}\to\{\text{2-simplices of the form $\bigtri fgh$}\},
      \]
    so there are at least as many 2-simplices of the form $\bigtri fgh$
    as there are quadruples $(k,s,l,t)$ satisfying $k^s\circ l^t=h$.
    Because $\zeta$ is a bijection, we can prove that $\zeta|_h$ is a
    surjection by noting $\ev(\gamma)=h$ whenever $\zeta(\gamma)$ is of the form
    $\bigtri fgh$.
    Indeed, $\zeta(\gamma)$ is of the form $\bigtri fg{\ev(\gamma)}$,
    so if $\zeta(\gamma)$ is also of the form $\bigtri fgh$
    then equality $\ev(\gamma)=h$ follows from Lemma \ref{lem:simultinaety}.
  \end{proof}
\end{lem}

We now have a way to keep track of 2-simplices of the form $\bigtri fgh$
by using formal composites of $f$, $g$, $f^{-1}$, and $g^{-1}$.
We will later define an equivalence relation on formal composites
by
  \[
  (k,s,l,t)\sim(k',s',l',t')\iff\ev(k,s,l,t)=\ev(k',s',l',t').
  \]
This equivalence relation gives a graph structure on the set $C(f,g)$,
with edges between equivalent elements.
Such graphs will be used to make easier the computations that underlie
the proofs of Propositions \ref{result:commutativity} and \ref{result:compoo}
below.
These proofs, found in Appendix \ref{appx:proof_of_compoo},
are the combinatorial heart of this paper.

Under the assumption that $f$ and $g$ satisfy
   \begin{equation*}
   f\neq g          ,\hspace{10pt}
   f\neq g^{-1}     ,\hspace{10pt}
   f^2\neq g        ,\hspace{10pt}
           \text{and}\hspace{10pt}
   f\neq g^2        ,
   \end{equation*}
Proposition \ref{result:commutativity} below gives conditions that are necessary and sufficient
for commutativity $fg=gf$.
Assuming that
   \begin{equation}\label{list:technical conditions}
   f\neq g          ,\hspace{10pt}
   f\neq g^{-1}     ,\hspace{10pt}
   f^2\neq g        ,\hspace{10pt}
   f\neq g^2        ,\hspace{10pt}
   f^2\neq g^{-1}   ,\hspace{10pt}
   f^{-1}\neq g^2   ,
   \end{equation}
Proposition \ref{result:compoo} establishes criteria necessary and
sufficient for a given endomorphism $h$ to satisfy
one (or both) of the equations $h=fg$ and $h=gf$.

As usual, we aim to encode relationship among
$f$, $g$, and $h$ in terms of the structure of the category $\Sd\,\G$
in a neighborhood of its objects $\barcons f$, $\barcons g$, and $\barcons h$.

\begin{prop}[Commutativity criterion]\label{result:commutativity}
\onehalfspacing
    Let $f$ and $g$ be non-identity endomorphisms in $\G$
    satisfying $f\neq g$ and $f\neq g^{-1}$ and $f^2\neq g$ and $f\neq g^2$.
    Then $fg$ equals $gf$ if and only if
    for every every third side $h$ of the triangle $\filltri fg$
    there are an even number of 2-simplices of the form $\bigtri fgh$.
\end{prop}

{\onehalfspacing
To prove the above, we define a graph $G(f,g)$ whose vertices
are the formal composites $C(f,g)$;
the graph is defined to have edge between
distinct formal composites $\gamma_1$ and $\gamma_2$
whenever $\ev(\gamma_1)$ equals $\ev(\gamma_2)$.
By Lemma \ref{result:end-to-end:multiplicity=composites},
the number of 2-simplices of the form $\bigtri fgh$
is equal to the size of the connected component of $G(f,g)$
whose elements $(k,s,l,t)$ satisfy $h=k^s\circ l^t$.
By a combinatorial argument, we show that
$fg$ equals $gf$ if and only if every connected component
of $G(f,g)$ has even cardinality.
A full proof is in Appendix \ref{appx:proof_of_compoo}.
\par}

\begin{prop}\label{result:compoo}
  Let $f$ and $g$ be non-identity endomorphisms of some object in $\G$
  satisfying $f\neq g$ and $f\neq g^{-1}$ and $f^2\neq g$
  and $f\neq g^2$ and $f^2\neq g^{-1}$ and $f^{-1}\neq g^2$.
  Let $h$ be another non-identity endomorphism of the same object in $\G$.
  The cases below give criteria under which $h$ is equal to $fg$ or to $gf$.
  Cases 1-2 apply when $f^2=g^2$, and cases 3-4 apply
  when $f^2\neq g^2$. All possibilities are exhausted.

{\doublespacing
  Case 1: Suppose that $f^2=\id=g^2$.
  \begin{itemize}
    \item If $fg=gf$, then $h=fg=gf$ if and only if $\exists_6\bigtri fgh$.
    \item If $fg\neq gf$, then $h$ equals $fg$ or $gf$ if and only if $\exists_3\bigtri fgh$.
  \end{itemize}

  Case 2: Suppose that $f^2=g^2$ and $f^2\neq\id$ and $g^2\neq\id$.
  \begin{itemize}
    \item If $fg=gf$, then $h=fg=gf$ if and only if $\exists_2\bigtri fgh$.
    \item If $fg\neq gf$, then $h$ equals $fg$ or $gf$ if and only if $\exists_1\bigtri fgh$ or $\exists_3\bigtri fgh$.
  \end{itemize}

  Case 3: Suppose that $f^2\neq g^2$, and either $f^2=\id$ or $g^2=\id$.
  \begin{itemize}
    \item If $fg=gf$, then $h=fg=gf$ if and only if $\exists_4\bigtri fgh$.
    \item If $fg\neq gf$, then $h$ equals $fg$ or $gf$
      if and only if $\exists_2\bigtri fgh$ or $\exists_3\bigtri fgh$.
  \end{itemize}

  Case 4: Suppose that $f^2\neq g^2$ and $f^2\neq\id$ and $g^2\neq\id$.
  \begin{itemize}
    \item If $fg=gf$, then $h=fg=gf$ if and only if \[\exists\bigtri fgh~\text{ and }~\exists\shifttext{2mm}{\bigtri{f^{-1}}gh}~\text{ and }~\exists\bigtri f{g^{-1}}h.\]
    \item If $fg\neq gf$ and $f^2=g^{-2}$,
      then $h$ equals $fg$ or $gf$ if and only if
      $\exists_2\bigtri{f^{-1}}gh$ and
      $\exists_2\bigtri f{g^{-1}}h$.
    \item If $fg\neq gf$ and $f^2\neq g^{-2}$, then $h$
    equals $fg$ or $gf$ if and only if one of the following is satisfied:
      \begin{enumerate}
        \item $\exists_1\bigtri fgh$ and $\exists_1\shifttext{0pt}{\bigtri{f^{-1}}gh}$
          and $\exists_1\bigtri f{g^{-1}}h$, or
        \item $\exists_1\shifttext{0pt}{\bigtri{f^{-1}}{g^{-1}}h}$
          and two of the following three existential statements hold:
          \[\exists_2\bigtri fgh,
              \hspace{25pt}
            \exists_2\shifttext{0pt}{\bigtri {f^{-1}}gh},
              \hspace{25pt}
            \exists_2\bigtri f{g^{-1}}h.
          \]
      \end{enumerate}
  \end{itemize}
\par}
\end{prop}

See Appendix \ref{appx:proof_of_compoo} for proof.


\section{Construction of $\psi$ and proof of the main theorem}\label{sec:Correspondence}
\numberwithin{thm}{subsection}
\numberwithin{equation}{subsection}
This section describes how to construct an isomorphism $\G\to\H$
given an isomorphism $\Sd\,\B\to\Sd\,\C$.

\begin{notn*}
For convenience, we will write $\simp m\C$ for the set of non-degenerate
$m$-simplices $[m]\to\C$, regarded as objects in $\Sd\,\C$.
\end{notn*}

This section's results are predicated on the fact that any isomorphism
$\Psi:\Sd\B\to\Sd\C$ sends $n$-simplices to $n$-simplices
(Proposition \ref{result:simp_n|->simp_n}), so we have induced maps
  \[\Ob(\B) \xrightarrow\cong\simp 0\B\xrightarrow\Psi\simp 0\C\xrightarrow\cong\Ob(\C)\]
and
  \[\B_1\xrightarrow\cong\simp 1\B\xrightarrow\Psi\simp 1\C\xrightarrow\cong\C_1,\]
where $\B_1$ denotes the set of non-identity morphisms in $\B$,
and similarly for $\C_1$.
Together these maps determine an isomorphism $\psi$ between
the undirected graphs that underlie $\B$ and $\C$.

We will show that the map $\psi$ is ``well behaved'' if $\B$ and $\C$
are groupoids: its restriction to any connected component
is a (possibly contravariant) isomorphism of categories.
This will be proved in two steps: we first consider
groups, that is, groupoids with one object
(Subsection \ref{subsec:singleObject}), and then
prove an analogous result for connected groupoids that have multiple objects
(in Subsection \ref{subsec:ConnectedMultiObj}).

In the general (non-connected) case,
we can selectively change the variance of $\psi$ on those connected
components where it was originally contravariant,
obtaining a map that is covariant everywhere.
To do this, we use the fact that any groupoid is isomorphic to its opposite.

\subsection{Defining an isomorphism $\B\to\C$ of undirected graphs}

\begin{lem}\label{result:pFac|->pFac}
Let $y$ be an object of $\Sd\,\B$.
Any isomorphism $\Psi:\Sd\,\B\xrightarrow\cong\Sd\,\C$ sends the proper faces
of $y$ bijectively to the proper faces of $\Psi\,y$.
  \begin{proof}
    Let $x$ be distinct from $y$.
    Then $x$ is a proper face of $y$ if and only if
    the hom-set $\hom[\Sd\,\B]xy$ is non-empty.
    Given the bijection $\hom[\Sd\,\B]xy\cong\hom[\Sd\,\C]{\Psi\,x}{\Psi\,y}$
    induced by $\Psi$, we see that $x$ is a proper face of $y$ if and only if
    $\Psi\,x$ is a proper face of $\Psi\,y$.
  \end{proof}
\end{lem}

It is a corollary that any isomorphism $\Psi:\Sd\,\B\to\Sd\,\C$
restricts to an isomorphism $\faces y\cong\faces{\Psi\,y}$
for each object $y$ of $\Sd\,\B$.
This follows because a category of faces
is a full subcategory.

\begin{prop}\label{result:simp_n|->simp_n}
  Any isomorphism $\Psi:\Sd\,\B\to\Sd\,\C$ restricts
  to a bijection $\simp m\B\cong\simp m\C$ for each $m$.
  \begin{proof}
    Suppose that $y$ is an object of $\Sd\,\B$.
    We shall prove that $y$ is an $m$-simplex only if $\Psi\,y$
    is an $m$-simplex.
    For each $x$ there is a bijection
      \begin{equation}\label{eqn:SdB(x,y)=SdC(Px,Py)}
        \hom[\Sd\,\B]xy\cong\hom[\Sd\,\C]{\Psi\,x}{\Psi\,y}
      \end{equation}
    induced by $\Psi$, and therefore we have bijection
      \[        \mt y
          =     \coprod_{x}\hom[\Sd\,\B]xy
          \cong \coprod_{\Psi\,x}\hom[\Sd\,\C]{\Psi\,x}{\Psi\,y}
          =     \mt{\Psi\,x}.
      \]
    Recall from Lemma \ref{result:2^m+1-1} that $y$ is an $m$-simplex
    if and only if the set $\mt{y}$ has $2^{m+1}-1$ elements.
    The desired result follows from
    the equality $\abs{\mt{y}}=\abs{\mt{\Psi\,y}}$,
    demonstrating that $y$ has the same dimension as $\Psi\,y$.
  \end{proof}
\end{prop}

\begin{cor}
If $\Psi:\Sd\,\B\to\Sd\,\C$ is an isomorphism and $y$ is an object of $\Sd\,\B$,
then for each $n$ there is a bijection $\mt[n]y\cong\mt[n]{\Psi\,y}$
induced by $\Psi$.
\end{cor}

\begin{con}\label{con:psi_from_Psi}
  Let $\Psi:\Sd\,\B\to\Sd\,\C$ be any isomorphism.
  Write $\B[1]$ for the set of non-identity morphisms of $\B$,
  and define $\C[1]$ similarly.
  There is a unique map $\psi:\B\to\C$ whose components
  make the following diagram commute:
         \begin{dontbotheriftheequationisoverlong}
         \raisebox{5pc}{ \xymatrix@=4em{
               \simp0\B	      \ar[d]_{\Psi}
                              \ar[r]^{\shifttext{2mm}{$\scriptstyle\bangle b\,\mapsto\,b$}}
               &	\B[0]	      \ar@{-->}[d]^{\Ob(\psi)}
               			          \ar[r]^{\id}
               &	\Mor\B      \ar@{-->}[d]^{\Mor(\psi)}
               &	\B[1]       \ar@{_{(}->}[l]
               &	\simp1\B    \ar[l]_{\shifttext{-4mm}{$\scriptstyle f\,\mapsfrom\,\barcons f$}}
               			          \ar[d]^{\Psi}
                                                   \\
               	\simp0\C	    \ar[r]^{\shifttext{2mm}{$\scriptstyle\bangle c\,\mapsto\,c$}}
               &	\C[0]
               		       	    \ar[r]^{\id}
               &	\Mor\C
               &	\C[1]  	    \ar@{_{(}->}[l]
               &	\simp1\C    \ar[l]_{\shifttext{-4mm}{$\scriptstyle f\,\mapsfrom\,\barcons f$}}
               }}.
        \end{dontbotheriftheequationisoverlong}
  The maps $\simp0\B\to\Ob\B$ and $\simp0\C\to\Ob\C$ above are the
  bijections sending 0-simplices $\bangle a$ to objects $a$.
  Similarly, the maps $\simp1\B\to\B[1]$
  and $\simp1\C\to\C[1]$ are the bijections sending 1-simplices $\barcons f$ to
  non-identity morphisms $f$.
  Note that $\Mor\B$ and $\Mor\C$ are isomorphic
  to $\Ob\B\amalg\B[1]$ and to $\Ob\C\amalg\C[1]$, respectively.
  The components of $\psi$ are uniquely determined by the above commutative
  diagram; note that $\psi$ might fail to be functorial.
  This map $\psi$ has the following properties:
  \begin{itemize}
    \item For any object $b$ in $\B$, we have $\Psi\bangle b=\bangle{\psi\,b}$
          and $\psi(\id_b)=\id_{\psi_b}$.
    \item For any non-identity morphism $f$ in $\B$,
    we have $\Psi\barcons f=\barcons{\psi\,f}$.
  \end{itemize}
\end{con}

The following result shows that this map $\psi$ is an isomorphism
between the undirected graphs that underlie $\B$ and $\C$.

\begin{lem}\label{result:dom/cod_preserved}
  Let $\Psi:\Sd\,\B\to\Sd\,\C$ be an isomorphism.
  Then $\psi$ is bijective on objects and arrows, satisfying
   \[
   \{\psi(\dom\,f),\psi(\cod\,f)\}
         =\{\dom(\psi\,f),\cod(\psi\,f)\}
    \]
  for each morphism $f$ of $\B$.
  \begin{proof}
    For any object $b$ in $\B$ we have
      \begin{align*}
         &\{\psi(\dom\,\id_b),\psi(\cod\,\id_b)\}
         =\{\psi\,b,\psi\,b\}\\
        &=\{\dom(\id_{\psi\,b}),\cod(\id_{\psi\,b})\}
        =\{\dom(\psi\,\id_b),\cod(\psi\,\id_b)\},
      \end{align*}
    so our result is true for identity morphisms.

    Suppose that $f$ is a non-identity morphism.
    By Lemma \ref{result:proper_faces=dom/cod},
    the proper faces of $\barcons f$ are
    $\bangle{\dom\,f}$ and $\bangle{\cod\,f}$,
    and the proper faces of $\barcons{\psi\,f}$ are
    $\bangle{\dom(\psi\,f)}$ and $\bangle{\cod(\psi\,f)}$.
    We have $\Psi\,\barcons f=\barcons{\psi\,f}$.
    By Lemma \ref{result:pFac|->pFac},
    the isomorphism $\Psi$ sends proper faces of $\barcons f$ bijectively
    to proper faces of $\barcons{\psi\,f}$, therefore
      \[
      \{\bangle{\psi(\dom\,f)},\bangle{\psi(\cod\,f)}\}
            =\{\Psi\bangle{\dom\,f},\Psi\bangle{\cod\,f}\}
            =\{\bangle{\dom(\psi\,f)},\bangle{\cod(\psi\,f)}\}.
            \qedhere
       \]
  \end{proof}
\end{lem}

If $\Phi:\Sd\,\C\to\Sd\,\B$ is the inverse isomorphism to $\Psi$,
and if we write $\phi$ for the map $\C\to\B$
obtained from $\Phi$ as per Construction \ref{con:psi_from_Psi},
then $\psi$ is inverse to $\phi$
(when these maps are regarded as morphisms of graphs).

It is now possible to prove that the restriction of $\Sd$ to the category
of small groupoids is conservative.
Suppose that $F:\G\to\H$ is any functor between
groupoids, write $\Psi$ for the functor $\Sd F:\Sd\,\G\to\Sd\,\H$,
and suppose that $\Psi$ is an isomorphism.
The induced map $\psi:\G\to\H$ satisfies
$\bangle{\psi\,a}=\Psi\bangle a=\bangle{F\,a}$ and
$\barcons{\psi\,g}=\Psi\barcons g=\barcons{F\,g}$
for each object $a$ and non-identity morphism $g$ in $\G$,
hence we have $\psi\,a=F\,a$ and $\psi\,g=F\,g$
by injectivity of $\bangle-$ and $\barcons-$.
Thus, we have $\psi\,a=F\,a$ and $\psi\,g=F\,g$ for
all objects $a$ and morphisms $g$ of $\G$.

By construction, we have $\psi\,\id_a=\id_{\psi a}=\id_{Fa}=F\,\id_a$
for each object $a$.
Because $F$ is equal to $\psi$ on each object and morphism in $\G$,
and because $\psi$ is bijective on objects and morphisms,
$F$ must be an isomorphism.

\subsection{Encoding $\B$ with $\Sd\,\B$, decoding $\C$ from $\Sd\,\C$.}

Suppose that some property of $\B$, e.g. invertibility of morphisms,
is encoded by the subdivision $\Sd\,\B$.
Assuming as before that $\Psi:\Sd\,\B\to\Sd\,\C$ is an isomorphism,
we can show that $\Sd\,\C$ encodes the same property in $\C$.
This style of argument will be used again and again:
the point of sections \ref{sec:encoding} and \ref{sec:encoding2}
was to encode properties of a given category
so that they could be translated to another category
by use of an isomorphism between the categories' subdivisions.

The following result shows that $\psi$ assigns morphisms in $\B$ to
morphisms in $\C$ in a way that respects faces of 2-simplices.

\begin{prop}\label{result:Phi(2-simp)=2-simp(Phi)}
\onehalfspacing
 Let $\Psi:\Sd\,\B\to\Sd\,\C$ be an isomorphism.
 Then $\Psi$ sends 2-simplices of the form $\bigtri fgh$
 bijectively to 2-simplices of the form
 $\shifttext{1mm}{\bigtri {\psi\,f}{\psi\,g}{\psi\,h}}$.
 \begin{proof}
   Because $\Psi$ is invertible, it will be good enough to show that each
   2-simplex of the form $\bigtri fgh$ is sent by $\Psi$
   to a 2-simplex of the form
   $\shifttext{1mm}{\bigtri {\psi\,f}{\psi\,g}{\psi\,h}}$.
   The idea of proof is that the form of a 2-simplex $y$
   is determined by the set $\mt[1]y$ of morphisms having target $y$
   and source equal to some 1-simplex.

   By Lemma \ref{result:faces_of_Sd2C_and_triangular_forms},
   $y$ is of the form $\bigtri fgh$ (for some non-identity morphisms
   $f$, $g$, and $h$) if and only if
   there are three morphisms in the set $\mt[1]y$,
   and these morphisms have respective domains $\barcons f$,
   $\barcons g$, and $\barcons h$.
   The map $\Psi$ gives a bijection
     \[
      \hom[\Sd\,\B]{\barcons f}y\xrightarrow\cong\hom[\Sd\,\C]{\Psi\,\barcons f}{\Psi\,y}
          =\hom[\Sd\,\C]{\barcons{\psi\,f}}{\Psi\,y},
          \]
   and similar bijections exist for $g$ and $h$.
   Therefore, $\Psi$ sends morphisms in $\mt[1]y$
   bijectively to morphisms in $\mt[1]{\Psi\,y}$, which have respective domains
   $\barcons{\psi\,f}$, $\barcons{\psi\,g}$, and $\barcons{\psi\,h}$.
   It follows that $y$ is of the form
   $\shifttext{1mm}{\bigtri {\psi\,f}{\psi\,g}{\psi\,h}}$.

   Implicit in this line of reasoning is, for example,
   the claim that if $y$ is of the
   form $\bigtri ffh$ for some distinct non-identity arrows $f$ and $h$,
   then we have
    \[|\hom[\Sd\,\B]{\barcons f}y|=2=|\hom[\Sd\,\C]{\barcons{\psi\,f}}{\Psi\,y}|\]
   and
    \[|\hom[\Sd\,\B]{\barcons h}y|=1=|\hom[\Sd\,\C]{\barcons{\psi\,h}}{\Psi\,y}|,\]
   hence $\Psi\,y$ is of the form $\bigtri{\psi\,f}{\psi\,f}{\psi\,h}$.
   Generally, it is important to keep track of how many arrows in
   $\mt[1]y$ have domain equal to each given 1-simplex.

   By Lemma \ref{result:faces_of_Sd2C_and_triangular_forms}, a 2-simplex
   $y$ is of the form $\bigtri fg\id$ if and only if there
   are two morphisms in $\mt[1]y$, having respective domains $\barcons f$
   and $\barcons g$.
   By the same argument as above, the elements
   of $\mt[1]y$ are sent bijectively to those of $\mt[1]{\Psi\,y}$,
   hence $\Psi\,y$ is of the form $\bigtri {\psi\,f}{\psi\,g}\id$.
 \end{proof}
\end{prop}


\vspace{0pt}

\begin{cor}\label{result:fillers/3rd|->fillers/3rd}
  Let $f$ and $g$ be morphisms in $\B$, and let $\Psi:\Sd\,\B\to\Sd\,\C$
  be an isomorphism.
  The fillers of $\filltri fg$ are sent bijectively
  by $\Psi$ to the fillers of
  $\shifttext{2mm}{\filltri{\psi\,f}{\psi\,g}}$.
  Similarly, the third sides of $\filltri fg$ are sent bijectively
  by $\psi$ to the third sides of
  $\shifttext{2mm}{\filltri{\psi\,f}{\psi\,g}}$.
  \onehalfspacing
  \begin{proof}
    By invertibility of $\Psi$ and $\psi$, it suffices to show that each filler
    or third side of $\filltri fg$ is sent to a filler or third side of
    $\shifttext{2mm}{\filltri{\psi\,f}{\psi\,g}}$.

    If $y$ is a filler for $\filltri fg$, then $y$ is of the form
    $\bigtri fgh$ for some $h$.
    The 2-simplex $\Psi\,y$ is of then of the form
    $\bigtri{\psi\,f}{\psi\,g}{\psi\,h}$, so $\Psi\,y$ is a filler for
    $\shifttext{2mm}{\filltri{\psi\,f}{\psi\,g}}$.

    If $h$ is a third side of the triangle $\filltri fg$, then there exists
    some 2-simplex $y$ of the form $\bigtri fgh$.
    Because $\Psi\,y$ is of the form $\bigtri{\psi\,f}{\psi\,g}{\psi\,h}$,
    the morphism $\psi\,h$ is a third side of the
    triangle $\shifttext{1mm}{\filltri{\psi\,f}{\psi\,g}}$.
  \end{proof}
\end{cor}

The result below follows from the fact that a graph isomorphism
preserves relationships between edges.

\begin{prop}\label{result:psi_preserves_ete/estes/eteo/eoteo/ur}
Let $\Psi:\Sd\,\B\to\Sd\,\C$ be an isomorphism.
An arrow $f$ of $\B$ is an endomorphism if
and only if $\psi\,f$ is an endomorphism.
Arrows $f$ and $g$ in $\B$ are end-to-end if and only if
$\psi\,f$ and $\psi\,g$ are end-to-end,
and similarly for ends-to-ends, end-to-endo, endo-to-endo, and unrelated
pairs of morphisms $f,g$.
In the end-to-endo case, $\psi$ preserves which morphism is the endomorphism.
\onehalfspacing
  \begin{proof}
    Because $\psi$ is an isomorphism of graphs, we have
    a set bijection
      \[\{\dom\,f,\cod\,f\}\cong\{\dom(\psi\,f),\cod(\psi\,f)\}.\]
    Thus, $f$ is an endomorphism if and only if
    $\abs{\{\dom\,f,\cod\,f\}}=1$,
    if and only if $\abs{\{\dom(\psi\,f),\cod(\psi\,f)\}}=1$,
    if and only if $\psi\,f$ is an endomorphism.

    Because it is bijective, $\psi$ commutes with $\cap$.
    Therefore, we have a set bijection
    \[\{\dom\,f,\cod\,f\}\cap\{\dom\,g,\cod\,g\}
    \cong\{\dom(\psi\,f),\cod(\psi\,f)\}\cap\{\dom(\psi\,g),\cod(\psi\,g)\}.\]
    The cardinalities of these sets determine the
    relationship between $f$ and $g$ (in the sense made
    precise by the Definitions \ref{defns:ete,estes,eteo,eoteo,ur} of
    end-to-end, ends-to-ends, end-to-endo,
    endo-to-endo, and unrelated morphisms).

    For example, $f$ and $g$ are end-to-end if and only if
    neither $f$ nor $g$ is an endomorphism
    and the intersection $\{\dom\,f,\cod\,f\}\cap\{\dom\,g,\cod\,g\}$
    has one element.
    This holds if and only if $\psi\,f$ and $\psi\,g$
    are non-endomorphisms satisfying
    \[\abs{\{\dom(\psi\,f),\cod(\psi\,f)\}\cap\{\dom(\psi\,g),\cod(\psi\,g)\}}=1.\]
    The other cases are similar.
  \end{proof}
\end{prop}

\begin{prop}\label{result:Psi(-f)=-Psi(f)}
  Let $\Psi:\Sd\,\B\to\Sd\,\C$ be an isomorphism, and let $f$
  be a morphism of $\B$.
  Then $f$ is invertible if and only if $\psi\,f$ is invertible.
  If $f$ is invertible then $\psi(f^{-1})$ is equal to $(\psi\,f)^{-1}$.
  \begin{proof}
    \onehalfspacing
    By construction, $\psi$ is bijective on identity arrows.
    Therefore, $f$ is an identity morphism in $\B$
    if and only if $\psi\,f$ is an identity morphism in $\C$.

    Now, suppose that $f$ is a non-identity morphism.
    The following are equivalent:
    \begin{enumerate}
    \setlength\itemsep{3mm}
      \item $f^2=\id$
      \item There is a 2-simplex in $\Sd\,\B$ of the form $\bigtri ff\id$
      \item There is a 2-simplex in $\Sd\,\C$ of the form $\bigtri{\psi\,f}{\psi\,f}\id$
      \item $(\psi\,f)^2=\id$.
    \end{enumerate}
    Equivalence of points (2) and (3) above follows from Proposition
    \ref{result:Phi(2-simp)=2-simp(Phi)}.
    Equivalence of points (1) and (2) and of points (3) and (4)
    follow from Lemma \ref{result:inv-criterion}.
    Thus, we have equality $\psi(f^{-1})=\psi\,f=(\psi\,f)^{-1}$.
    
    Finally, supposing that $f$ and $g$ are distinct non-identity morphisms
    in $\B$, the following are equivalent:
    \begin{enumerate}
    \setlength\itemsep{3mm}
      \item $f=g^{-1}$
      \item There are two 2-simplices in $\Sd\,\B$ of the form $\bigtri fg\id$
      \item There are two 2-simplices in $\Sd\,\C$ of the form $\bigtri{\psi\,f}{\psi\,g}\id$
      \item $\psi\,f=(\psi\,g)^{-1}$.
    \end{enumerate}
    As before, equivalence of (2) and (3) follows from Proposition
    \ref{result:Phi(2-simp)=2-simp(Phi)},
    and Lemma \ref{result:inv-criterion} gives the other equivalences.
    We conclude that $\psi(f^{-1})=\psi\,g=(\psi\,f)^{-1}$.
  \end{proof}
\end{prop}

\begin{cor}\label{result:SdG=SdC=>C_a_groupoid}
  Let $\G$ be a groupoid and let $\C$ be a category.
  if $\Sd\,\G$ is isomorphic to $\Sd\,\C$ then $\C$ is a groupoid.
  \begin{proof}
    Take an isomorphism $\Psi:\Sd\,\G\to\Sd\,\C$, and write $\psi$ for
    the induce isomorphism of undirected graphs.
    Every morphism $g$ of $\G$ is invertible,
    so every morphism $\psi\,g$ of $\C$ is invertible,
    having inverse equal to $\psi(g^{-1})$.
  \end{proof}
\end{cor}

{\onehalfspacing
The proof of Proposition \ref{result:Psi(-f)=-Psi(f)}
uses the fact that 2-simplices of the form
$\bigtri fg\id$ encode whether $f$ and $g$ are inverse morphisms.
An analogous statement about $\psi\,f$ and $\psi\,g$ is derived
by considering 2-simplices of the form $\bigtri{\psi\,f}{\psi\,g}\id$.
Generally, the relationship between morphisms $f$ and $g$ is determined
by the 2-simplex fillers for $\filltri fg$,
as well as by the fillers for $\shifttext{1mm}{\filltri{f^{-1}}g}$ and
$\filltri f{g^{-1}}$ and $\shifttext{1mm}{\filltri{f^{-1}}{g^{-1}}}$.

We now restrict our attention to groupoids.
The construction below concerns a map $\psi'$
that satisfies $\dom(\psi\,f)=\cod(\psi'\,f)$ and $\cod(\psi\,f)=\dom(\psi'\,f)$
for each morphism $f$.
This will be useful later, as there is no a priori guarantee
that the map $\psi$ should be covariant.
\par}

\begin{con}\label{result:psi'}
\onehalfspacing
  Let $\G$ and $\H$ be groupoids and let
  $\Psi:\Sd\,\G\to\Sd\,\H$ be an isomorphism.
  Write $\Psi':\Sd\,\G\to\Sd\,\H$
  for the composite $\Psi\circ\alpha_\G$,
  where $\alpha_\G:\Sd\,\G\to\Sd\,\G$ is the map
  $\barcons{f_m|\cdots|f_1}\mapsto\barcons{f_1^{-1}|\cdots|f_m^{-1}}$
  defined in Section \ref{sec:Lemmas}.
  Write $\psi'$ for the map obtained from $\Psi'$ per Construction
  \ref{con:psi_from_Psi}.
\end{con}

\begin{prop}\label{psi'=psi-1}
Define $\psi'$ as above.
  We have $\psi'a=\psi\,a$
  for each object $a$ of $\G$,
  and $\psi'f=\psi\,f^{-1}$
  for each morphism $f$ of $\G$.
  \begin{proof}
    We have
      \[\barcons{\psi'\,a}=\Psi'\barcons a=\Psi\barcons a=\barcons{\psi\,a}\]
    and
      \[\bangle{\psi'\,f}=\Psi'\bangle f=\Psi\bangle{f^{-1}}=\bangle{\psi(f^{-1})}\]
    for each object $a$ and non-identity morphism $f$ of $\G$.
    By injectivity of the maps sending $a$ to $\bangle a$
    and $f$ to $\barcons f$,
    we have $\psi'\,a=\psi\,a$ and $\psi'\,f=\psi\,f^{-1}$.
    \end{proof}
\end{prop}

Note that the map $\psi'$ defined above is equal to the composite
$\psi\circ\gamma$, where $\gamma$ is the (covariant) functor $\op\G\to\G$
defined on morphisms by $g\mapsto g^{-1}$.

\subsection{Single-object groupoids}\label{subsec:singleObject}

It is now possible to prove this paper's main
result in the special case of groups.
We will first state two Lemmas concerning squares
and commutativity in groupoids that have one object.
Next, we will show that if $\G$ and $\H$ are groups
and if $\Psi:\Sd\,\G\to\Sd\,\H$ is an isomorphism,
then the identity
\begin{equation}\label{eqn:biglem_for_groups}
  \psi(\{fg,gf\})=\{(\psi\,f)(\psi\,g),(\psi\,g)(\psi\,f)\}
\end{equation}
is satisfied for all $f$ and $g$ in $\G$.
Finally, a group-theoretic result of Bourbaki will be used
to show that the condition (\ref{eqn:biglem_for_groups})
is sufficient to establish an isomorphism between $\G$ and $\H$.

We begin by proving that $\psi(f^2)=(\psi\,f)^2$ and $\psi(f^3)=(\psi\,f)^3$.

\begin{lem}\label{result:psi(f2)=(psif)2}
{\singlespacing
  Let $\G$ and $\H$ be groups,
  let $\Psi:\Sd\,\G\to\Sd\,\H$ be an isomorphism,
  and let $f$ and $h$ be any morphisms of $\G$.
  Then $f^2=h$ if and only if $(\psi\,f)^2=\psi\,h$.
  \par}
  \onehalfspacing
    \begin{proof}
      By Lemma \ref{result:Psi(-f)=-Psi(f)}, we have $f^2=\id$ if and only
      if $(\psi\,f)^2=\id$.
      Thus, if $f^2=\id$ then
      \[f^2=h\iff h=\id\iff\psi\,h=\id\iff(\psi\,f)^2=\psi\,h\]
      because $\psi$ sends the identity morphism of $\G$
      to the identity morphism of $\H$.

      Now, suppose that $f^2$ is a non-identity morphism.
      Then the result follows from Lemma \ref{result:composite_g^2},
      which states that $f^2$ equals $h$ if and only if
      there exists a 2-simplex in $\Sd\,\G$ of the form $\bigtri ffh$.
      The following are equivalent:
      \begin{enumerate}
      \setlength\itemsep{3mm}
        \item $f^2=h$
        \item There is a 2-simplex in $\Sd\,\G$ of the form $\bigtri ffh$
        \item There is a 2-simplex in $\Sd\,\H$ of the form $\bigtri {\psi\,f}{\psi\,f}{\psi\,h}$
        \item $(\psi\,f)^2=\psi\,h$.
        \qedhere
      \end{enumerate}
    \end{proof}
\end{lem}

It follows in particular that $\psi(f^2)=\psi(f)^2$.

\begin{lem}\label{result:psi(f3)=(psif)3}
{\singlespacing
  Let $\G$ and $\H$ be groups,
  let $\Psi:\Sd\,\G\to\Sd\,\H$ be an isomorphism,
  and let $f$ and $h$ be any morphisms of $\G$.
  Then $f^3=h$ if and only if $(\psi\,f)^3=\psi\,h$.
  \par}
  \onehalfspacing
    \begin{proof}
      The following are equivalent:
        \begin{enumerate}
          \item $f^2$ is inverse to $f$
          \item $\psi(f^2)$ is inverse to $\psi\,f$
          \item $(\psi\,f)^2$ is inverse to $\psi\,f$.
        \end{enumerate}
      Therefore,
        \[f^3=\id\iff f^2=f^{-1}\iff (\psi\,f)^2=(\psi\,f)^{-1}\iff(\psi\,f)^3=\id.\]
      Now, suppose that $f^3$ is not the identity morphism of $\G$.
      We cite Corollary \ref{result:composite_g^3}, which gives criteria
      for determining whether $f^3=h$ by counting 2-simplices
      of the form $\bigtri f{f^2}h$,
      to show that the following are equivalent:
      \begin{enumerate}
      \setlength\itemsep{3mm}
        \item $f^3=h$
        \item either $h=f^{-1}$ and $\exists_4\bigtri f{f^2}h$, or
            $h\neq f^{-1}$ and $\exists_2\bigtri f{f^2}h$
        \item either $\psi\,h=\psi\,f^{-1}$ and $\exists_4\mathmakebox[12mm]{\shifttext{-2mm}{\bigtri {\psi\,f}{(\psi\,f)^2}{\psi\,h}}}$, or
            $\psi\,h\neq \psi\,f^{-1}$ and $\exists_2\shifttext{0pt}{\bigtri {\psi\,f}{(\psi\,f)^2}{\psi\,h}}$
        \item $(\psi\,f)^3=\psi\,h$.
        \qedhere
      \end{enumerate}
    \end{proof}
\end{lem}

The following two results are the fruit of
calculations occurring in the Appendix.
The first result shows that commutativity is preserved by $\psi$;
the second shows that, in the context of groups,
$\psi(fg)$ is equal to $(\psi\,f)(\psi\,g)$ or to $(\psi\,g)(\psi\,f)$.

\begin{lem}\label{result:psi_preserves_commutativity}
  Let $\G$ and $\H$ be groups,
  let $\Psi:\Sd\,\G\to\Sd\,\H$ be an isomorphism,
  and let $f$ and $g$ be morphisms of $\G$.
  Then $fg=gf$ if and only if $(\psi f)(\psi g)=(\psi g)(\psi f)$ 
    \begin{proof}
    First, suppose that $f$ or $g$ is equal to the identity arrow in $\G$.
    Then $fg=gf$ follows trivially, and we have
    $(\psi\,f)(\psi\,g)=(\psi\,g)(\psi\,f)$ because either $\psi\,f$
    or $\psi\,g$ must be equal to the identity arrow in $\H$.

    Now, suppose that $f$ and $g$ are not identity arrows.
    Then at least one of the following statements must be true:
      \begin{enumerate}
      \item $f=g$,
      \item $f=g^{-1}$,
      \item $f^2=g$,
      \item $f=g^2$,
      \item $f\neq g$ and $f\neq g^{-1}$ and $f^2\neq g$ and $f\neq g^2$.
      \end{enumerate}
    We will show that the Lemma holds in each of the above cases.
      \begin{itemize}
      \setlength\itemsep{2mm}
      \item In case (1), we have $fg=f^2=gf$ and
        $(\psi\,f)(\psi\,g)=(\psi\,f)(\psi\,f)=(\psi\,g)(\psi\,f)$.
      \item In case (2), we have $fg=ff^{-1}=\id=f^{-1}f=gf$ and
        \[(\psi\,f)(\psi\,g)=(\psi\,f)(\psi\,f)^{-1}
          =\id=(\psi\,f)^{-1}(\psi\,f)=(\psi\,g)(\psi\,f).\]
      \item In case (3), we have $fg=f^3=gf$ and
        \[(\psi\,f)(\psi\,g)=(\psi\,f)^3=(\psi\,g)(\psi\,f).\]
      \item Case (4) is analogous to case (3).
      \item {\onehalfspacing For case (5), recall from Proposition \ref{result:commutativity}
        that, under the given conditions on $f$ and $g$,
        we have $fg=gf$ if and only if there are an even number of 2-simplices
        of the form $\bigtri fgh$ for each $h$.
        Thus, the following are equivalent:
          \begin{enumerate}
          \setlength\itemsep{2mm}
          \item $fg=gf$
          \item There are an even number of $\bigtri fgh$ for each $h$
          \item There are an even number of $\shifttext{2mm}{\bigtri{\psi\,f}{\psi\,g}{\psi\,h}}$ for each $h$
          \item $(\psi\,f)(\psi\,g)=(\psi\,g)(\psi\,f)$.
          \qedhere
          \end{enumerate}\par}
      \end{itemize}
    \end{proof}
\end{lem}

\begin{prop}
  Let $\G$ and $\H$ be groups,
  and let $\Psi:\Sd\,\G\to\Sd\,\H$ be an isomorphism.
  For any pair $f,g$ of morphisms in $\G$
  there is equality
  $\psi(fg)=(\psi\,f)(\psi\,g)$ or
  $\psi(fg)=(\psi\,g)(\psi\,f)$.
    \begin{proof}
    Let $f$ and $g$ be some arrows in $\G$.
    We will use Proposition \ref{result:compoo}, which gives criteria for
    determining when the a given morphism $h$ is equal to one of the composites
    $fg$ and $gf$.
    When the hypotheses of Proposition \ref{result:compoo} are not satisfied,
    we resort to more elementary means.

    First, if $f$ (resp. $g$) is equal to the identity arrow of $\G$,
    then $\psi\,f$ (resp. $\psi\,g$) is the identity arrow of $\H$,
    in which case the proof follows easily.
    For example, if $f=\id$ then $\psi(fg)=\psi\,g=(\psi\,f)(\psi\,g)$.

    Now, suppose that $f$ and $g$ are not identity arrows.
    Then at least one of the following statements must be true:
      \begin{enumerate}
      \item $f=g$,
      \item $f=g^{-1}$,
      \item $f^2=g$,
      \item $f=g^2$,
      \item $f^2=g^{-1}$,
      \item $f^{-1}=g^2$,
      \item $f\neq g$ and $f\neq g^{-1}$ and $f^2\neq g$ and $f\neq g^2$
            and $f^2\neq g^{-1}$ and $f^{-1}\neq g^2$.
      \end{enumerate}
    We will show that the result holds in each of the above cases.
      \begin{itemize}
      \setlength\itemsep{2mm}
      \item In case (1), we have $\psi(fg)=\psi(f^2)=(\psi\,f)(\psi\,f)=(\psi\,f)(\psi\,g)$.
      \item In case (2), we have $\psi(fg)=\psi(\id)=\id=(\psi\,f)(\psi\,f)^{-1}=(\psi\,f)(\psi\,g)$.
      \item In case (3), we have $\psi(fg)=\psi(f^3)=(\psi\,f)(\psi\,f)^2=(\psi\,f)(\psi\,g)$.
      \item Case (4) is analogous to case (3).
      \item In case (5), we have
        \begin{align*}
        \psi(fg)&=\psi(f^{-1})=(\psi\,f)^{-1}=(\psi\,f)(\psi\,f)^{-2}\\
        &=(\psi\,f)((\psi\,f)^2)^{-1}=(\psi\,f)(\psi(g^{-1}))^{-1}=(\psi\,f)(\psi\,g).
        \end{align*}
      \item Case (6) is analogous to case (5).
      \end{itemize}
      Case (7) above is precisely the list of
      hypotheses for Proposition \ref{result:compoo}.
      To apply that Lemma, we will suppose that $h$ is a morphism
      in $\G$, and we will show that
      \[\text{$h$ equals $fg$ or $gf$}\iff\text{$\psi\,h$ equals $(\psi\,f)(\psi\,g)$ or $(\psi\,g)(\psi\,f)$},\]
      from which follows the stated result that $\psi(fg)$
      is equal to $(\psi\,f)(\psi\,g)$ or to $(\psi\,g)(\psi\,f)$.

      Suppose that conditions (7) above
      are satisfied by $f$ and $g$.
      This section's earlier results guarantee that $\psi(\id)=\id$, and
      that the following hold for all $h$ in $\G$:
        \begin{itemize}
          \item $\psi(h^{-1})=\psi(h)^{-1}$,
          \item $\psi(h^2)=\psi(h)^2$, and
          \item $\psi(h^3)=\psi(h)^3$.
        \end{itemize}
     It is a corollary of these facts, and of the bijectivity of $\psi$,
     that $f^i=g^j$ if and only if $(\psi\,f)^i=(\psi\,g)^j$
     for any $i,j$ in $\{-2,-1,0,1,2\}$.
     Therefore, having assumed that the conditions (7) are satisfied by $f$ and
     $g$, we deduce that the same conditions are satisfied by $\psi\,f$ and
     $\psi\,g$:
          \begin{gather*}
          \psi\,f\neq \psi\,g ,\hspace{5mm}
          \psi\,f\neq (\psi\,g)^{-1} ,\hspace{5mm}
          (\psi\,f)^2\neq \psi\,g ,\hspace{5mm}    \\
          \psi\,f\neq (\psi\,g)^2 ,\hspace{5mm}
         (\psi\,f)^2\neq (\psi\,g)^{-1} ,\hspace{5mm}
         (\psi\,f)^{-1}\neq (\psi\,g)^2.
         \end{gather*}

      The following points demonstrate that each case listed in the statement
      of Proposition \ref{result:compoo} applies equally as well
      to $f$ and $g$ as it does to $\psi\,f$ and $\psi\,g$:
        \begin{itemize}
          \item $f$ and $g$ commute if and only if
                $\psi\,f$ and $\psi\,g$ commute,
          \item $f^2=g^2$ if and only if
                $(\psi\,f)^2=(\psi\,g)^2$,
          \item $f^2=g^{-2}$ if and only if
                $(\psi\,f)^2=(\psi\,g)^{-2}$, and
          \item $f^2=\id$ if and only if $(\psi\,f)^2=\id$
                (and similarly for $g$).
        \end{itemize}
     We will work through the most complicated case explicitly.
     Suppose that $fg\neq gf$ and $f^2\neq g^{-2}$,
     as in case 4 subcase 3 of Proposition \ref{result:compoo}.
     Then the following are equivalent for
     any non-identity morphism $h$ in $\G$:
      \begin{enumerate}
      \doublespacing
        \item $h$ equals $fg$ or $gf$
        \item One of the following is satisfied:
                        \begin{enumerate}
                          \item $\exists_1\bigtri fgh$ and $\exists_1\shifttext{0pt}{\bigtri{f^{-1}}gh}$
                            and $\exists_1\bigtri f{g^{-1}}h$, or
                          \item $\exists_1\shifttext{0pt}{\bigtri{f^{-1}}{g^{-1}}h}$
                            and two of the following three existential statements hold:
                            \[\exists_2\bigtri fgh,
                                \hspace{25pt}
                              \exists_2\shifttext{0pt}{\bigtri {f^{-1}}gh},
                                \hspace{25pt}
                              \exists_2\bigtri f{g^{-1}}h.
                            \]
                        \end{enumerate}
        \item One of the following is satisfied:
                        \begin{enumerate}
                          \item $\exists_1\bigtri {\psi\,f}{\psi\,g}{\psi\,h}$ and $\exists_1\mathmakebox[15mm]{\shifttext{3mm}{\bigtri{(\psi\,f)^{-1}}{\psi\,g}{\psi\,h}}}$
                            and $\exists_1\mathmakebox[15mm]{\shifttext{-3mm}{\bigtri {\psi\,f}{(\psi\,g)^{-1}}{\psi\,h}}}$, or
                            \vspace{1mm}
                          \item $\exists_1\mathmakebox[18mm]{\shifttext{3mm}{\bigtri{(\psi\,f)^{-1}}{(\psi\,g)^{-1}}{\psi\,h}}}$
                            and two of the following three statements hold:
                            \[\exists_2\shifttext{0mm}{\bigtri {\psi\,f}{\psi\,g}{\psi\,h}},
                                \hspace{25pt}
                              \exists_2\mathmakebox[15mm]{\shifttext{3mm}{\bigtri{(\psi\,f)^{-1}}{\psi\,g}{\psi\,h}}},
                                \hspace{25pt}
                              \exists_2\shifttext{0mm}{\bigtri {\psi\,f}{(\psi\,g)^{-1}}{\psi\,h}}.
                            \]
                        \end{enumerate}
        \item $\psi\,h$ equals $(\psi\,f)(\psi\,g)$ or $(\psi\,g)(\psi\,f)$
      \end{enumerate}
                \onehalfspacing
                We have used the fact that 2-simplices of the form $\bigtri fgh$
                are in bijection with those of the form $\bigtri{\psi\,f}{\psi\,g}{\psi\,h}$,
                and similarly when $f$ (resp. $g$) has been replaced with $f^{-1}$ (resp. $g^{-1}$).
     The other cases listed in Proposition \ref{result:compoo}
     have proofs that are analogous to the above,
     thus in every case we may conclude that $h$ equals $fg$ or $gf$
     if and only if $\psi\,h$ equals $(\psi\,f)(\psi\,g)$
     or $(\psi\,g)(\psi\,f)$.
    \end{proof}
\end{prop}

\begin{cor}\label{psi(fg,gf)=(psi_f_psi_g,_psi_g_psi_f)}
  Let $\G$ and $\H$ be groups,
  and let $\Psi:\Sd\,\G\to\Sd\,\H$ be an isomorphism.
  For any pair $f,g$ of morphisms in $\G$
  there is equality
  \[\psi(\{fg,gf\})=\{(\psi\,f)(\psi\,g),(\psi\,g)(\psi\,f)\}.\]
    \begin{proof}
    From the previous Proposition, we know that the right hand side
    above is a subset of the left hand side:
    \[\psi(\{fg,gf\})=\{\psi(fg),\psi(gf)\}\subseteq\{(\psi\,f)(\psi\,g),(\psi\,g)(\psi\,f)\}.\]
    Let $\phi:\H\to\G$ denote the inverse map to $\psi$, induced
    by the isomorphism $\Phi:\Sd\,\H\to\Sd\,\G$ that is inverse
    to $\Psi$.
    Again by the previous Proposition, we have subset inclusion
    \begin{align*}
    \phi\bigl(\{(\psi\,f)(\psi\,g),(\psi\,g)(\psi\,f)\}\bigr)
    =&\bigl\{\phi\bigl((\psi\,f)(\psi\,g)\bigr),\phi\bigl((\psi\,g)(\psi\,f)\bigr)\bigr\}\\
    &\subseteq\{(\phi\,\psi\,f)(\phi\,\psi\,g),(\phi\,\psi\,g)(\phi\,\psi\,f)\}
    =\{fg,gf\}.
    \end{align*}
    It follows that $\psi$ is a bijection between $\{fg,gf\}$ and
    $\{(\psi\,f)(\psi\,g),(\psi\,g)(\psi\,f)\}.$
    \end{proof}
\end{cor}

We can now make use of the following result from group theory,
which is applicable because the category of small groups is equivalent to
the category of small single-object groupoids.

\begin{lem}[Bourbaki]\label{result:Bourbaki}
    \textnormal{
    \begin{displayquote}
    Let $\text{G}$, $\text{G}'$ be two groups,
    $f:\text{G}\to \text{G}'$ a mapping such that,
    for two arbitrary elements $x,y$ of $\text{G}$,
    $f(xy)=f(x)f(y)$ \textit{or} $f(xy)=f(y)f(x)$.
    It is proposed to prove that $f$ is a homomorphism of
    $\text{G}$ into $\text{G}'$ \textit{or} a homomorphism of
    $\text{G}$ into the opposite group $\text{G}'^0$
    (in other words, either $f(xy)=f(x)f(y)$ for \textit{every} ordered pair
    $(x,y)$ or $f(xy)=f(y)f(x)$ for \textit{every} ordered pair $(x,y)$).
    \end{displayquote}
    }
\end{lem}

The above is quoted directly from Bourbaki's book \cite{Bourbaki},
in ``Exercises for \S 4'', Problem 26, p.139.
For an elementary proof, see Lemma 4 in
Mansfield's paper \cite{Mansfield} on the group determinant.
We now state our theorem in the special case of single-object groupoids.

\begin{thm}\label{thm:groups}
  Let $\G$ and $\H$ be groups.
  and let $\Psi:\Sd\,\G\to\Sd\,\H$ be an isomorphism.
  Then there exists an isomorphism $P:\G\to\H$.
  \begin{proof}
    By the Bourbaki Lemma, $\psi$ is a possibly-contravariant
    group isomorphism.
    If $\psi$ is covariant then we set $P=\psi$,
    and if $\psi$ is contravariant then we set $P$
    equal to the map $\psi'$ (from Construction \ref{result:psi'})
    defined by $\psi'\,a=\psi\,a$ and
    $\psi'\,f=(\psi\,f)^{-1}$
    for objects $a$ and morphisms $f$ in $\G$.
    Thus, for non-identity morphisms $f$ in $\G$,
    \[ \raisebox{1pc}{$\displaystyle
          P\,f=\begin{cases}
            \psi\,f       &\text{if $\psi$ is covariant}\\
            \psi(f^{-1})  &\text{if $\psi$ is contravariant}
          \end{cases}$}.\qedhere  \]
  \end{proof}
\end{thm}

\subsection{Connected multi-object groupoids}\label{subsec:ConnectedMultiObj}

This subsection establishes a result analogous to the previous theorem
pertaining to connected groupoids that have multiple objects.
The combinatorics from the appendix will not be necessary;
this subsection's result is significantly easier to prove
than that of the previous subsection.

Let $\G$ and $\H$ be connected groupoids that have more
than one object each,
and let $\Psi:\Sd\,\G\to\Sd\,\H$ be an isomorphism.
Given an arrow $f$ of $\G$ that is not an endomorphism,
say that $\psi$ is \textit{covariant} at $f$ if
  \begin{equation}\label{eqn:p(domf)=dom(pf)}
    \psi(\dom\,f)=\dom(\psi\,f)~\text{ and }~
    \psi(\cod\,f)=\cod(\psi\,f).
  \end{equation}
Similarly, say that $\psi$ is \textit{contravariant} at a given
non-endomorphism $f$ if
  \begin{equation}\label{eqn:p(domf)=dom(pf)}
    \psi(\dom\,f)=\cod(\psi\,f)~\text{ and }~
    \psi(\cod\,f)=\dom(\psi\,f).
  \end{equation}

Recall from Construction \ref{result:psi'},
that there is a map $\Psi':\Sd\,\G\to\Sd\,\H$ defined as the composite
$\Psi\circ\alpha_\G$, and that $\psi'$ is the map $\G\to\H$
induced by $\Psi'$.
This map $\psi$ is characterized by the equalities $\psi'(a)=\psi(a)$
for objects $a$ and $\psi'(f)=\psi(f)^{-1}$ for morphisms $f$.

We may assume without loss of generality
that $\psi$ is covariant at some non-endomorphism
$f$ in $\G$; if this fails to be the case then $\psi$ is contravariant
everywhere, so we may consider $\psi'$ instead of $\psi$
to obtain a map that is covariant.
Assuming that $\psi$ is covariant at some non-endomorphism,
and that the groupoids $\G$ and $\H$ are connected, we will prove that $\psi$
is covariant everywhere,
that is, $\psi$ is a morphism of \textit{directed} graphs.

Note first that $\psi$ is covariant at $f$ if and only
if $\psi'$ is contravariant at $f$
Indeed, by the inversion Lemma \ref{result:Psi(-f)=-Psi(f)}
we have equality $\psi'\,f=(\psi\,f)^{-1}=\psi(f^{-1})$,
hence the domain of $\psi\,f$ equals the codomain of $\psi'\,f$,
and vice versa.

A consequence will be that $\psi$ is functorial: we will show that
there is equality $\psi(f\circ g)=(\psi\,f)\circ(\psi\,g)$
for every pair $f,g$ of morphisms in $\G$,
thereby establishing an analog to Theorem \ref{thm:groups} pertaining to
connected groupoids that have multiple objects.

\begin{prop}\label{result:psi_preserves_sequential/parallel}
{\singlespacing
  Let $f$ and $g$ be morphisms of $\G$ satisfying
    \begin{enumerate}
      \item neither $f$ nor $g$ is an endomorphism, and
      \item the sets $\{\dom\,f,\cod\,f\}$ and $\{\dom\,g,\cod\,g\}$
            are not disjoint, that is, $f$ and $g$ are not unrelated.
    \end{enumerate}
  If $f$ and $g$ are end-to-end,
  $f$ and $g$ are sequential if and only if $\psi\,f$
  and $\psi\,g$ are sequential.
  Similarly, if $f$ and $g$ are ends-to-ends,
  then $f$ and $g$ are parallel (resp. opposed)
  if and only if $\psi\,f$ and $\psi\,g$ are parallel (resp. opposed).
  \par}\onehalfspacing
  \begin{proof}
    Given the conditions (1) and (2),
    the morphisms $f$ and $g$ must be end-to-end or ends-to-ends.
    Recall from Proposition \ref{result:psi_preserves_ete/estes/eteo/eoteo/ur}
    that $f$ and $g$ are end-to-end if and only if $\psi\,f$
    and $\psi\,g$ are end-to-end, and similarly in the case where
    $f$ and $g$ are ends-to-ends.

    Suppose first that $f$ and $g$ are end-to-end sequential morphisms.
    By Proposition \ref{result:composability-criterion-for-end-to-end}, this is true
    if and only if there is a unique filler for the triangle $\filltri fg$.
    By Corollary \ref{result:fillers/3rd|->fillers/3rd},
    this is true if and only if
    there is a unique filler for $\filltri{\psi\,f}{\psi\,g}$.
    Applying Proposition \ref{result:composability-criterion-for-end-to-end} once
    more, this is true if and only if $\psi\,f$ and $\psi\,g$ are sequential.

    Now, suppose that $f$ and $g$ are non-sequential end-to-end morphisms,
    that is, if $f$ and $g$ are coinitial or coterminal.
    By Proposition \ref{result:composability-criterion-for-end-to-end},
    this is true if and only if
    there are multiple fillers for $\filltri fg$.
    By Corollary \ref{result:fillers/3rd|->fillers/3rd},
    $\filltri{\psi\,f}{\psi\,g}$.
    Applying Proposition \ref{result:composability-criterion-for-end-to-end} again,
    this is true if and only if $\psi\,f$ and $\psi\,g$
    are coterminal or coinitial.
    
    Suppose now that $f$ and $g$ are parallel ends-to-ends morphisms.
    From Proposition \ref{result:composability-criterion-for-ends-to-ends},
    we know that this occurs if and only if
    there are four fillers for $\filltri fg$.
    By Corollary \ref{result:fillers/3rd|->fillers/3rd},
    this is true if and only if there are four fillers for triangle
    $\filltri{\psi\,f}{\psi\,g}$ has four fillers.
    Applying Proposition \ref{result:composability-criterion-for-ends-to-ends}
    once more, this is true if and only if $\psi\,f$ and $\psi\,g$ are parallel.
    The proof is similar in case $f$ and $g$ are opposed.
  \end{proof}
\end{prop}

\begin{cor}\label{local_covariance->global_covariance}
  Let $f$ and $g$ be morphisms of $\G$
  satisfying conditions (1) and (2) from
  Proposition \ref{result:psi_preserves_sequential/parallel}.
  If $\psi$ is covariant at $f$, then it is covariant at $g$.
  It follows from connectedness
  that if $\psi$ is covariant at any non-endomorphism $f$,
  then $\psi$ is covariant at every non-endomorphism in $\G$.
  \begin{proof}
    Suppose first that $\psi$ is covariant at
    $f$ and that $f$ and $g$ are parallel.
    Because $\psi\,f$ and $\psi\,g$ are parallel, we have
      \[  \dom(\psi\,g)=\dom(\psi\,f)=\psi(\dom\,f)=\psi(\dom\,g)
      \]
    and similarly for codomains,
    showing that $\psi$ is covariant at $g$.
    The proof is similar in case $f$ and $g$
    are opposed or sequential.

    Suppose now that $\psi$ is covariant at $f$ and that
    $f$ and $g$ are coterminal or coinitial.
    Then $f$ and $g^{-1}$ are sequential,
    so $\psi$ is covariant at $g^{-1}$.
    It follows that $\psi$ is covariant at $g$ because
    $g^{-1}$ and $g$ are opposed.

    Now, we claim that covariance at any one non-endomorphism $f$
    is sufficient to guarantee covariance everywhere.
    Indeed, suppose that $\psi$ is covariant at $f$,
    writing $a=\dom\,f$ and $b=\cod\,f$.
    Letting $g:c\to d$ be any other non-endomorphism in $\G$,
    we will show that $\psi$ is covariant at $g$.
    We have proved above that this is true if $f$ and $g$
    are related; we still need to consider the case of unrelated morphisms,
    where the objects $a,b,c,d$ are all distinct.

    Because $\G$ is a connected groupoid, there exists
    some arrow $h:b\to c$.
    Then $\psi$ is covariant at $h$ because $f$ and $h$ are sequential,
    and it follows that $\psi$ is covariant at $g$ because $h$ and $g$
    are sequential.
  \end{proof}
\end{cor}

Thus we can assume without loss of generality
that if $\G$ and $\H$ are connected groupoids
and $\Psi:\Sd\,\G\to\Sd\,\H$ is an isomorphism,
then $\psi:\G\to\H$ is an isomorphism of directed graphs;
if this fails to be the case, we consider $\psi'$ instead of $\psi$.
It is now possible to demonstrate functorality.

\begin{prop}
  Assume that $\G$ and $\H$ are connected multi-object groupoids,
  $\Psi:\Sd\,\G\to\Sd\,\H$ is an isomorphism, and $\psi$ is covariant
  at every morphism.
  Then $\psi(f\circ g)=(\psi\,f)\circ(\psi\,g)$ for every pair
  $f,g$ of morphisms in $\G$.
  \begin{proof}
    We consider all possible forms a composite could take.

    Suppose first that $f$ and $g$ are sequential end-to-end morphisms,
    and write $h$ for the composite of $f$ and $g$.
    Recall Corollary \ref{result:encoding_sequential_fg=h},
    which says that $h$ is the composite of $f$ and $g$
    if and only if there is only one filler of $\filltri fg$, and
    $h$ is the third side of that filler.
    We use Corollary \ref{result:fillers/3rd|->fillers/3rd}
    to show that there is only one filler of
    $\shifttext{1mm}{\filltri{\psi\,f}{\psi\,g}}$,
    and $\psi\,h$ is the third side of that filler.
    It follows that $\psi\,h$ is the composite of $\psi\,f$ and $\psi\,g$.

    If $f$ and $g$ are opposed ends-to-ends morphisms then the equation
        \[\{\psi(fg),\psi(gf)\}=\{(\psi\,f)(\psi\,g),(\psi\,g)(\psi\,f)\}\]
    follows from Corollary \ref{result:encoding_opposed_fg=h},
    which says that there exactly two fillers of $\filltri fg$,
    and that these fillers have respective third sides equal to
    $fg$ and to $gf$.
    Indeed, by Corollary \ref{result:fillers/3rd|->fillers/3rd}
    there are exactly two fillers of $\shifttext{1mm}{\filltri{\psi\,f}{\psi\,g}}$,
    and these fillers have respective third sides equal to
    $(\psi\,f)(\psi\,g)$ and to $(\psi\,g)(\psi\,f)$,
    hence $\psi$ gives a bijection between the respective sets
      \[\text{$\{fg,gf\}$ and $\{(\psi\,f)(\psi\,g),(\psi\,g)(\psi\,f)\}$}\]
    of third sides.
    To see that $\psi(fg)=(\psi\,f)(\psi\,g)$, note that we cannot have
    $\psi(fg)=(\psi\,g)(\psi\,f)$ because $\psi(fg)$ is an endomorphism
    of $\cod(\psi\,f)$ whereas $(\psi\,g)(\psi\,f)$ is an endomorphism of
    $\dom(\psi\,f)$.

    {\onehalfspacing
    Suppose now that $f$ and $g$ are end-to-endo.
    Recall Proposition \ref{result:comp->o/<-o}, which states
    that $h$ is the composite of $f$ and $g$ if and only if
    $h$ is parallel to $f$ and $h$ is a third side of the triangle
    $\shifttext{1mm}{\filltri{f^{-1}}g}$.
    Thus, $\psi$ sends the composite of $f$ and $g$
    to a third side of
    $~\shifttext{6mm}{\filltri{(\psi\,f)^{-1}}{\psi\,g}}\hspace{3mm}$
    that is parallel to $\psi\,f$.
    It follows from Proposition \ref{result:comp->o/<-o}
    that this third side is the composite of $\psi\,f$ and $\psi\,g$,
    proving the desired result.
    Note that this same proof holds whether $\dom\,f=\dom\,g=\cod\,g$
    or $\cod\,f=\dom\,g=\cod\,g$.
    \par}

    Finally, suppose that $f$ and $g$ are both endomorphisms
    of some object $a$ in $\G$.
    Given that $\G$ is connected and has multiple objects,
    there exists some non-endomorphism $k$ in $\Mor(\G)$
    satisfying $\dom\,k=a$. We then have the following:
    \begin{minipage}[b]{.94\textwidth}
    \begin{align*}
      \psi(f\circ g)
      &= \psi(k^{-1}\circ k\circ f\circ g)                          \\
      &= \psi(k^{-1})\circ\psi(k\circ f\circ g)
        & \text{($k^{-1}$ and $k\circ f\circ g$ are ends-to-ends)}  \\
      &= \psi(k^{-1})\circ\psi(k\circ f)\circ\psi(g)
        & \text{($k\circ f$ and $g$ are end-to-endo)}               \\
      &= \psi(k^{-1})\circ\psi(k)\circ\psi(f)\circ\psi(g)
        & \text{($k$ and $f$ are end-to-endo)}                      \\
      &= \psi(f)\circ\psi(g)
        & \text{(because $\psi(k^{-1})$ equals $(\psi\,k)^{-1}$)}
    \end{align*}
    \end{minipage}\qedhere
  \end{proof}
\end{prop}

\begin{thm}\label{thm:multi-obj_groupoids}
Let $\G$ and $\H$ be connected multi-object groupoids,
and let $\Psi:\Sd\,\G\to\Sd\,\H$ be an isomorphism.
Then there exists an isomorphism $P:\G\to\H$.
  \begin{proof}
  Pick some non-endomorphism $f$ in $\Mor(\G)$.
  If $\psi$ is covariant at $f$, then it follows from Corollary
  \ref{local_covariance->global_covariance} that $\psi$ is
  covariant everywhere on $\G$.
  It then follows from the previous proposition that $\psi$
  is functorial, so we can set $P=\psi$.

  Suppose instead that $\psi$ is contravariant at $f$.
  It follows that $\psi'$ is covariant at $f$,
  where $\psi'$ is the map
    from Construction \ref{result:psi'}, defined by
    $\psi'\,a=\psi\,a$ and $\psi'\,f=(\psi\,f)^{-1}$.
  By the same argument as above, it follows that $\psi'$ is
  covariant everywhere on $\G$, hence $\psi'$ is functorial.
  In this case, we set $P=\psi'$.
  \end{proof}
\end{thm}

\subsection{Statement and proof of the groupoid isomorphism theorem}

It is now possible to prove this paper's main result.

\begin{thm}
  Let $\G$ and $\H$ be groupoids.
  Any isomorphism $\Psi:\Sd\,\G\to\Sd\,\H$ induces
  an isomorphism $P:\G\to\H$.
  \begin{proof}
    Suppose that $\Psi:\Sd\,\G\to\Sd\,\H$ is an isomorphism.
    Recall from Lemma \ref{result:connected/coproducts_preserved}
    that we have an isomorphism $\amalg_i(\Sd\,\C_i)\cong\Sd(\amalg_i\C_i)$
    for any set $\{\C_i\}$ of small categories,
    and that a category $\C$ is connected if and only if its subdivision
    $\Sd\,\C$ is connected.

    There exists some index sets $I$ and $J$, and some
    collections $\{\G_i\}_{i\in I}$
    and $\{\H_j\}_{j\in J}$ of connected groupoids such that
    $\G=\amalg_i\G_i$ and $\H=\amalg_j\H_j$.
    We obtain a composite isomorphism
      \[            \Gamma: \coprod_i\Bigl(\Sd\,\G_i\Bigr)
          \xrightarrow\cong \Sd\Bigl(\coprod_i\G_i\Bigr)
          \xrightarrow\Psi  \Sd\Bigl(\coprod_j\H_j\Bigr)
          \xrightarrow\cong \coprod_j\Bigl(\Sd\,\H_j\Bigr).
      \]
    The image of a connected category is connected.
    Thus, for each $i$ there is some $j$ such that the image
    $\Gamma(\Sd\,\G_i)$ is a subcategory of $\Sd\,\H_j$.
    The inverse isomorphism $\Gamma^{-1}$ must take
    $\Sd\,\H_j$ back into $\Sd\,\G_i$ because $\Sd\,\H_j$ is connected.
    The connected components of $\Sd\,\G$ are in bijection with the
    connected components of $\Sd\,\H$, so we can use $I$ to reindex the
    components $\{H_j\}$ of $\H$ so as to obtain an isomorphism
    $\Psi_i:\Sd\,\G_i\to\Sd\,\H_i$ for each $i$ in $I$.
    These maps $\Psi_i$ are the evident restrictions of $\Gamma$,
    satisfying $\Gamma=\amalg_i\Psi_i$.

    For each $\Psi_i$ we obtain an isomorphism $P_i:\G_i\xrightarrow\cong\H_i$,
    either from Theorem \ref{thm:groups}
    or Theorem \ref{thm:multi-obj_groupoids},
    depending on whether $\G_i$ has one object or many.
    Setting $P=\amalg_iP_i$ we have isomorphism
    $P:\amalg_i\G_i\to\amalg_i\H_i$ between $\G$ and $\H$.
  \end{proof}
\end{thm}



\subsection*{Acknowledgements}
First, I'd like to thank my mentor, Jonathan Rubin, who has been extremely
helpful with the organization of this paper, and who has
been responsible for much of the content.
Henry Chan and Claudio Gonzales offered
generous support with proofreading and notation.
Many thanks to Peter May for promoting the study of categorical subdivision
and for encouraging me to develop these ideas.



\appendix
\numberwithin{thm}{section}
\numberwithin{equation}{section}
\section{Proof of Propositions \ref{result:commutativity} and \ref{result:compoo}}
\label{appx:proof_of_compoo}

{\onehalfspacing
Here we will prove the final two Lemmas of Section
\ref{sec:encoding}.
We assume throughout that $f$ and $g$ are distinct non-identity
endomorphisms of some common object in a groupoid $\G$,
satisfying the identities $f\neq g$ and $f\neq g^{-1}$ and $f^2\neq g$ and $f\neq g^2$.

Proposition \ref{result:commutativity} states that, under the above assumptions
on $f$ and $g$, we have commutativity $fg=gf$ if and only if
for every third side $h$ of the triangle $\filltri fg$ there are an even
number of 2-simplices of the form $\shifttext{3mm}{\bigtri fgh}$.
Proposition \ref{result:compoo} gives criteria for a given third side $h$ of
the triangle $\filltri fg$ to satisfy one of the equations $h=fg$ and $h=gf$.
\par}

\subsection*{Proof of Proposition \ref{result:commutativity}.}

{\onehalfspacing

Write $\fills(f,g)$ for the set of non-degenerate 2-simplex fillers
of the triangle $\filltri fg$.
Our assumptions that $f\neq g$ and $f^2\neq g$ and $f\neq g^2$
guarantee that the set $\fills(f,g)$ contains six distinct 2-simplices
(see Lemma \ref{result:end-to-end:multiplicity=composites}).
We have assumed $f\neq g^{-1}$ to avoid the trivial case
where no non-identity morphisms $h$ satisfy $h=fg$ or $h=gf$.

Recall that the set $C(f,g)$ of formal composites consists
of the quadruples $(k,s,l,t)$ satisfying
  \begin{enumerate}
    \item $(k,l)$ is equal to $(f,g)$ or $(g,f)$, and
    \item $(s,t)$ is equal to $(1,1)$ or $(1,-1)$ or $(-1,1)$.
  \end{enumerate}
Our proof makes use of the bijection
    \[
    \zeta:C(f,g)\xrightarrow\cong\fills(f,g)
    \]
from Lemma \ref{result:end-to-end:multiplicity=composites},
defined as below:
      \begin{align*}
        (f,1,g,1)&\mapsto\barcons{f|g}          
      & (g,1,f,1)&\mapsto\barcons{g|f}        \\
        (f,-1,g,1)&\mapsto\barcons{f|f^{-1}g}   
      & (g,1,f,-1)&\mapsto\barcons{gf^{-1}|f} \\
        (f,1,g,-1)&\mapsto\barcons{fg^{-1}|g}   
      & (g,-1,f,1)&\mapsto\barcons{g|g^{-1}f} . 
      \end{align*}

As before, we write $\ev:C(f,g)\to\Mor(\G)$ for the evaluation
map that sends $(k,s,l,t)$ to $k^sl^t$.
By Lemma \ref{result:end-to-end:multiplicity=composites}, each
element $\gamma$ of $C(f,g)$ is sent by $\zeta$ to a 2-simplex of the form
{ \setlength\abovedisplayskip{10pt}
  \setlength\belowdisplayskip{10pt}
  \[\bigtri fg{\ev(\gamma)}.\]
  }

We define an equivalence relation $\sim$ on $C(f,g)$ by
  \[(k,s,l,t)\sim(k',s',l',t')\iff k^sl^t=k'^{s'}l'^{t'}.\]
We can use $\sim$ to define an undirected graph structure on $C(f,g)$
by placing an edge between two formal composites $\gamma,\psi\in C(f,g)$
if and only if $\gamma\sim\psi$ and $\gamma\neq\psi$.
Then $\gamma\sim\psi$ exactly when either
  \begin{enumerate}
    \item $\gamma=\psi$, or
    \item there is an edge between $\gamma$ and $\psi$.
  \end{enumerate}
Consequently, the components of this graph are precisely the fibers
of the evaluation map $\ev$.

This graph will be denoted $G(f,g)$, and denoted the \textit{graph of ev}.
Working with the connected components of $G(f,g)$,
we reduce the problem of counting 2-simplices
to a problem of determining fiber sizes.

Note that the graph structure of $G(f,g)$ is encoded by
the categorical structure of $\Sd\,\G$: two nodes $\gamma$ and $\psi$
have an edge between them if and only if the 2-simplices $\zeta(\gamma)$
and $\zeta(\psi)$, each of which fills the triangle $\filltri fg$,
have third sides that are equal.
Explicitly, we have a commutative diagram
  \[\xymatrix{
  C(f,g) \ar[r]^-\ev \ar[d]_\zeta^{\rotatebox{-90}{$\cong$}}
& \{\text{non-identity morphisms of $\G$}\} \ar[d]^{\theta^\G_1}\\
  \fills(f,g) \ar[r]_{3rd_{f,g}}
& simp_1(\Sd\,\G)
  }\]
where $\theta_1^\G$ is the bijection sending $d$ to $\barcons d$,
and where $3rd_{f,g}$ is the function that sends
each 2-simplex of the form $\bigtri fgd$ to the third side $\barcons d$.

To prove Proposition \ref{result:commutativity}, it will suffice to show that $fg=gf$
if and only if every connected component of $G(f,g)$ has even cardinality.
This is because the size of each connected component $\ev^{-1}\{h\}$
is equal to the number of non-degenerate 2-simplices of the form $\bigtri fgh$.

We will always draw the vertices of $G(f,g)$ in the following configuration
\par}
  \[\setlength{\tabcolsep}{1cm}
  \begin{tabular}{cc}
  $(f,1,g,1)$     &   $(g,1,f,1)$    \\\\
  $(f,-1,g,1)$    &   $(g,1,f,-1)$   \\\\
  $(f,1,g,-1)$    &   $(g,-1,f,1)$   
  \end{tabular}\]
and will usually suppress the labels on the vertices to points $*$.
Configurations in the graph of $\ev$ correspond to equational conditions.
Specifically, we have the table below, where the following are equivalent:
  \begin{enumerate}
    \item the upper equation holds
    \item the lower graph is a subgraph of $G(f,g)$
    \item a single edge of the lower graph is in $G(f,g)$.
  \end{enumerate}
\begin{dontbotheriftheequationisoverlong}
\label{table:subgraphs}
\begin{tabular}{|c|c|c|c|c|c|c|c|}\hline
Equation
& $fg=gf$
& $f^2=\id$
& $g^2=\id$
& $f^2=g^2$
& $fgf=g$
& $gfg=f$
& $fg^{-1}f=g$
\\\hline Subgraph
  & \sixgraph\sixabcdef\endsixgraph
  & \sixgraph\sixacbd\endsixgraph
  & \sixgraph\sixaebf\endsixgraph
  & \sixgraph\sixcedf\endsixgraph
  & \sixgraph\sixadbc\endsixgraph
  & \sixgraph\sixafbe\endsixgraph
  & \sixgraph\sixcfde\endsixgraph
\\\hline
\end{tabular}
\end{dontbotheriftheequationisoverlong}
The equivalence of conditions 1-3 above can be verified
by considering if-and-only-if statements such as
  \begin{align*}
      fg=gf &\iff
      \ev(f,1,g,1)=\ev(g,1,f,1) \\
      &\iff \ev(f,-1,g,1)=\ev(g,1,f,-1) \\
      &\iff \ev(f,1,g,-1)=\ev(g,-1,f,1).
  \end{align*}
If we suppose that $fg=gf$, then the graph
  \begin{equation}\label{eqn:sixgraph_fg=gf}
  \begin{minipage}{1cm}
    \sixgraph\sixabcdef\endsixgraph
    \end{minipage}
  \end{equation}
is a subgraph of $G(f,g)$.
It follows that each connected component of $G(f,g)$ must have
cardinality equal to 2, 4, or 6.
This is because equivalence classes in a coarser equivalence relation
are unions of equivalence classes in a finer one --
the equivalence relation illustrated by graph \ref{eqn:sixgraph_fg=gf}
is the finest possible among graphs $G(f,g)$ satisfying $fg=gf$.

We now prove the converse implication:
if each connected component of $G(f,g)$ has even cardinality,
then $fg=gf$.
Let $P(f,g)$ denote the partition of $C(f,g)$ defined by
  \[
  \Biggl\{
  \Bigl\{(f,1,g,1) ,(g,1,f,1) \Bigr\},
  \Bigl\{(f,-1,g,1),(g,1,f,-1)\Bigr\},
  \Bigl\{(f,1,g,-1),(g,-1,f,1)\Bigr\}
  \Biggr\}
  .\]
Note that we have commutativity $fg=gf$ if and only if
the evaluation map $\ev:C(f,g)\to\Mor(\G)$ factors through $P(f,g)$.
Indeed we have $fg=gf$ if and only if there is equality
  \[\ev(f,s,g,t)=f^sg^t=g^tf^s=\ev(g,t,f,s)\]
for every pair $s,t\in\{-1,1\}$,
if and only if each fiber $\ev^{-1}\{h\}$
is a union of some classes in $P(f,g)$.

Note that there are three possible partitions of 6 into even integers, namely
  \[
     6 \hspace{15pt}\text{and}\hspace{15pt}
     2+4 \hspace{15pt}\text{and}\hspace{15pt}
     2+2+2  .\]
We will prove first that if the partition of $C(f,g)$ into fibers of $\ev$
corresponds to either 6 or $2+4$, then $fg=gf$.
Supposing that $C(f,g)$ is partitioned as 6 or as $2+4$, there exists
some morphism $h$ in $\G$ such that $\ev^{-1}\{h\}$ has four or more elements.
It follows that $\ev^{-1}\{h\}$ must contain an entire class in $P(f,g)$.
This implies $f^sg^t=g^tf^s$ for some $s$ and $t$, and commutativity follows.

To complete our proof, suppose that the partition of $C(f,g)$ given by $\ev$
corresponds to $2+2+2$.
Then, the graph $G(f,g)$ consists of exactly three disjoint edges.
Our strategy will be to show that, under these conditions,
$f^2\neq\id$ and $g^2\neq\id$ and $f^2\neq g^2$.
We then use these facts to prove that $fg=gf$.

Suppose for contradiction that $f^2=\id$.
Then $G(f,g)$ must contain \[\sixgraph\sixacbd\endsixgraph\] as a subgraph,
hence $G(f,g)$ is exactly equal to the graph
  \[\sixgraph\sixacbd\draw[shorten >=-1mm, shorten <=-1mm](e) -- (f);\endsixgraph.\]
But the bottom edge implies $fg^{-1}=g^{-1}f$, and multiplying on the left
and right by $g$ implies $fg=gf$.
This, however, implies that all three horizontal edges must be in $G(f,g)$,
contradicting our assumption of a $2+2+2$ partition.
Therefore, $f^2\neq\id$.
The arguments for $g^2\neq\id$ and $f^2\neq g^2$ are entirely analogous.

Now, suppose for contradiction that $fg\neq gf$,
so that no edge below
  \[\sixgraph\sixabcdef\sixacbd\sixaebf\sixcedf\endsixgraph\]
is contained in the graph $G(f,g)$.
In this case, $G(f,g)$ must be a subgraph of
  \[\sixgraph\sixadbc\sixafbe\sixcfde\endsixgraph.\]
Since $G(f,g)$ must consist of three disjoint edges,
it must be of the form
  \[
    \raisebox{25pt}{\begin{minipage}{.1\textwidth}
    \centering
    \sixgraph
        \draw[shorten >=-2mm, shorten <=-2mm](a) -- (d);
        \draw[shorten >=-2mm, shorten <=-2mm](c) -- (f);
        \draw[shorten >=-1mm, shorten <=-1mm](b) -- (e);
        \endsixgraph
    \end{minipage}
    \hspace{15pt}\text{or}
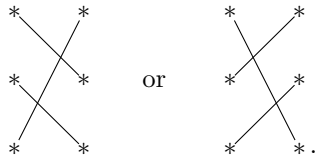
\hspace{15pt}
    \begin{minipage}{.1\textwidth}
    \centering
    \sixgraph
        \draw[shorten >=-2mm, shorten <=-2mm](b) -- (c);
        \draw[shorten >=-2mm, shorten <=-2mm](d) -- (e);
        \draw[shorten >=-1mm, shorten <=-1mm](a) -- (f);
        \endsixgraph
    \end{minipage}}
    .\]
In either case we have a contradiction, for
  \[\sixgraph
        \draw[shorten >=-2mm, shorten <=-2mm](a) -- (d);
    \endsixgraph\]
is a subgraph of $G(f,g)$ if and only if
  \[\sixgraph
        \draw[shorten >=-2mm, shorten <=-2mm](b) -- (c);
    \endsixgraph\]
is a subgraph of $G(f,g)$.
We conclude that if $G(f,g)$ is partitioned into connected components each
having size 2, then the assumption $fg\neq gf$ is contradictory.
\qed

\subsection*{Proof of Proposition \ref{result:compoo}.}
{\onehalfspacing
We will prove Proposition \ref{result:compoo} case-by-case,
making use of the graph $G(f,g)$ defined above.
Recall that the number of non-degenerate 2-simplices of the form $\bigtri fgh$
is equal to the size of the fiber $\ev^{-1}\{h\}$, where $ev:C(f,g)\to\Mor(\G)$
is the evaluation map defined by $(k,s,l,t)\mapsto k^sl^t$.
In addition to our previous assumptions that $f$ and $g$ are non-identity
endomorphisms of some common object in $\G$ satisfying
$f\neq g$ and $f\neq g^{-1}$ and $f^2\neq g$ and $f\neq g^2$,
we assume also that $f^2\neq g^{-1}$ and $f^{-1}\neq g^2$.
These last two assumptions guarantee that there are bijections
  \[ C(f^{-1},g)\cong\fills(f^{-1},g)
    ~\text{ and }~
     C(f,g^{-1})\cong\fills(f,g^{-1}) \]
as per Lemma \ref{result:end-to-end:multiplicity=composites}.
These two bijections will be useful in the proof of case 4.

As before, configurations in the graph of $\ev$ correspond
to equational conditions involving $f$ and $g$.
Recall that any distinct points belonging to the same
connected component of $G(f,g)$ must have an edge between them.
In other words, every full connected subgraph of $G(f,g)$ is a complete graph.

{\onehalfspacing

{\em\ul{Case 1, subcase 1:
Suppose that $f^2=\id=g^2$ and $fg=gf$.}}\hfill\\
Under these assumptions, the graph of $\ev$ consists of a
single connected component.
Indeed, $G(f,g)$ is the complete graph on six vertices, as below:
  \[\sixgraph\sixabcdef\sixadbc\sixafbe\sixacbd
    \sixaebf\sixcfde\sixcedf\endsixgraph.\]
Thus, every third side of $\filltri fg$ is equal to $fg$ and to $gf$.
In other words, all six 2-simplices in $\fills(f,g)$
are of the form $\bigtri fg{fg}$.
  \begin{align*}
    \text{$h=fg=gf$}
      &\iff
    \text{$\exists_6\bigtri fgh$}
\\[5pt]
    \text{$h\neq fg$ and $h\neq gf$}
      &\iff
    \text{$\nexists\bigtri fgh$}
  \end{align*}

{\em\ul{Case 1, subcase 2:
Suppose that $f^2=\id=g^2$ and $fg\neq gf.$}}\hfill\\
Consulting the table on page \pageref{table:subgraphs},
the condition $f^2=\id=g^2$ guarantees that
  \[\sixgraph\sixacbd\sixaebf\sixcedf\endsixgraph\]
is a subgraph of $G(f,g)$.
There can be no edges connecting the two halves of the
graph above; if there were such a connection, then
$G(f,g)$  would have only one connected component, in
which case it would be complete and we would have $fg=gf$.
Thus, $G(f,g)$ consists of two connected components,
and is exactly equal to the graph above.

The fibers $\ev^{-1}\{fg\}$ and $\ev^{-1}\{gf\}$
have three elements each, and every third side of $\filltri fg$
is equal to $fg$ or to $gf$.
In other words, every 2-simplex in $\fills(f,g)$
is of the form $\bigtri fg{fg}$ or $\bigtri fg{gf}$,
and we have the following:
  \begin{align*}
    \text{$h=fg$ or $h=gf$}
      &\iff
    \text{$\exists_3\bigtri fgh$}
\\[5pt]
    \text{$h\neq fg$ and $h\neq gf$}
      &\iff
    \text{$\nexists\bigtri fgh$}.
  \end{align*}

{\em\ul{Case 2, subcase 1:
Suppose that $f^2=g^2$ and $fg=gf$ and $f^2\neq\id$ and $g^2\neq\id$.}}\hfill\\
Because $f^2=g^2$ and $fg=gf$, and because each connected full subgraph
of $G(f,g)$ is complete, the below
  \[\sixgraph\sixabcdef\sixcfde\sixcedf\endsixgraph\]
must be a subgraph of $G(f,g)$.
There can be no edges connected the two components of the graph above,
for otherwise $G(f,g)$ would be the complete graph,
contradicting our assumptions that $f^2\neq\id$ and $g^2\neq\id$.

Thus, if $h=fg=gf$ then the fiber $\ev^{-1}\{h\}$ has two elements,
and if $h$ is not equal to $fg$ or $gf$ then the fiber
$\ev^{-1}\{h\}$ is empty or has four elements.
  \begin{align*}
    \text{$h=fg=gf$}
      &\iff
    \text{$\exists_2\bigtri fgh$}
\\[5pt]
    \text{$h\neq fg$ and $h\neq gf$}
      &\iff
    \text{$\nexists\bigtri fgh$ or $\exists_4\bigtri fgh$}
  \end{align*}
Explicitly, we have $h=fg=gf$ if and only if there are two 2-simplices
of the form $\bigtri fgh$, namely $\barcons{f|g}$ and $\barcons{g|f}$.

\clearpage
{\em\ul{Case 2, subcase 2:
Suppose that $f^2=g^2$ and $fg\neq gf$ and $f^2\neq\id$ and $g^2\neq\id$.}}\hfill\\
Because $fg\neq gf$ and $f^2\neq\id$ and $g^2\neq\id$,
the graph $G(f,g)$ cannot contain any of the edges below:
  \[\sixgraph\sixabcdef\sixacbd\sixaebf\endsixgraph.\]
We have $f^2=g^2$,
so the only question is what $(f,1,g,1)$ and $(g,1,f,1)$ are connected to.
The graph $G(f,g)$ must be equal to one of the two displayed below.
  \[\setlength{\tabcolsep}{.5cm}\begin{tabular}{cc}
      \sixgraph\sixcedf\sixadbc\sixafbe\endsixgraph
  &   \sixgraph\sixcedf\endsixgraph
  \\  \text{if $gfg=f$}
  &   \text{if $gfg\neq f$}
    \end{tabular}\]

If $gfg=f$ then the fibers of $fg$ and $gf$ both have size 3.
If $gfg\neq f$ then the fibers both have size 1.
If the fiber $\ev^{-1}\{h\}$ is non-empty and if $h$
is not equal to $fg$ or $gf$, then we must have
$gfg\neq f$ and fiber size 2.
  \begin{align*}
    \text{$h=fg$ or $h=gf$}
      &\iff
    \text{$\exists_1\bigtri fgh$ or $\exists_3\bigtri fgh$}
\\[5pt]
    \text{$h\neq fg$ and $h\neq gf$}
      &\iff
    \text{$\nexists\bigtri fgh$ or $\exists_2\bigtri fgh$}
  \end{align*}

{\em\ul{Case 3, subcase 1:
Suppose $f^2\neq g^2$ and $fg=gf$, and either $f^2=\id$ or $g^2=\id$.}}\hfill\\
The graph $G(f,g)$ contains one of the below two graphs as a subgraph.
  \[\setlength{\tabcolsep}{.5cm}\begin{tabular}{cc}
      \sixgraph\sixabcdef\sixacbd\sixadbc\endsixgraph
  &   \sixgraph\sixabcdef\sixaebf\sixafbe\endsixgraph
  \\  \text{if $f^2=\id$}
  &   \text{if $g^2=\id$}
  \\  \text{and $g^2\neq\id$}
  &   \text{and $f^2\neq\id$}
    \end{tabular}\]
Note that we cannot have both $f^2=\id$ and $g^2=\id$,
because $f^2\neq g^2$.
We claim that $G(f,g)$ is exactly equal to one
of the two graphs above.
Indeed, supposing for contradiction that $G(f,g)$ contains one of the above
as a \textit{proper} subgraph, we have only one connected component,
hence $G(f,g)$ is complete and our assumption $f^2\neq g^2$ is contradicted.

\clearpage
The above two graphs are isomorphic. In both cases,
if $h=fg=gf$ then the fiber of $h$ has four elements.
If $h\neq fg$ and $h\neq gf$,
then the fiber $\ev^{-1}\{h\}$ is empty or has two elements.
  \begin{align*}
    \text{$h=fg=gf$}
      &\iff
    \text{$\exists_4\bigtri fgh$}
\\[5pt]
    \text{$h\neq fg$ and $h\neq gf$}
      &\iff
    \text{$\nexists\bigtri fgh$ or $\exists_2\bigtri fgh$}
  \end{align*}

{\em\ul{Case 3, subcase 2:
Suppose $f^2\neq g^2$ and $fg\neq gf$, and either $f^2=\id$ or $g^2=\id$.}}\hfill\\
Then $G(f,g)$ is one of the four graphs displayed below.

  \[\setlength{\tabcolsep}{.5cm}\begin{tabular}{cccc}
      \sixgraph\sixacbd\sixafbe\sixcfde\endsixgraph
  &   \sixgraph\sixacbd\endsixgraph
  &   \sixgraph\sixaebf\endsixgraph
  &   \sixgraph\sixaebf\sixadbc\sixcfde\endsixgraph
  \\  \text{if $f^2=\id$}
  &   \text{if $f^2=\id$}
  &   \text{if $g^2=\id$}
  &   \text{if $g^2=\id$}
  \\  \text{and $fg^{-1}f=g$}
  &   \text{and $fg^{-1}f\neq g$}
  &   \text{and $fg^{-1}f\neq g$}
  &   \text{and $fg^{-1}f=g$}
    \end{tabular}\]
The two graphs in the middle are guaranteed by the equations
$f^2=\id$ and $g^2=\id$, respectively.
The graph on the far left is the only possible extension
of the graph next to it; any other
extension would contradict one of the assumptions
$f^2\neq g^2$ and $fg\neq gf$.
Similarly, the graph on the far right is the only
possible extension of the graph next to it.

Thus, if $h=fg$ or $h=gf$ then the fiber of $h$ has
three elements (if $fg^{-1}f=g$) or two elements (if $fg^{-1}f\neq g$).
On the other hand, if the fiber $\ev^{-1}\{h\}$ is non-empty
and if $h$ is not equal to $fg$ or $gf$, then we must have
$fg^{-1}f=g$ and fiber size 1.
  \begin{align*}
    \text{$h=fg$ or $h=gf$}
      &\iff
    \text{$\exists_2\bigtri fgh$ or $\exists_3\bigtri fgh$}
\\[5pt]
    \text{$h\neq fg$ and $h\neq gf$}
      &\iff
    \text{$\nexists\bigtri fgh$ or $\exists_1\bigtri fgh$}
  \end{align*}

{\em\ul{Case 4: Suppose $f^2\neq g^2$ and 
  $f^2\neq\id$ and $g^2\neq\id$.}}\hfill\\
Under these assumptions,
the number of non-degenerate 2-simplices of the form $\bigtri fgh$
does not completely determine whether the given
morphism $h$ is equal to one of $fg$ or $gf$.
Our strategy here is to consider not only fillers of the triangle
$\filltri fg$, but also to consider the fillers for
$\shifttext{2mm}{\filltri{f^{-1}}g}$ and $\filltri f{g^{-1}}$ and
$\shifttext{2mm}{\filltri{f^{-1}}{g^{-1}}}$.
Let $C'(f,g)$ denote the superset of $C(f,g)$
consisting of quadruples $(k,s,l,t)$ that satisfy
  \begin{itemize}
    \item $(k,l)$ is equal to $(f,g)$ or $(g,f)$, and
    \item $s$ and $t$ are both elements of the set $\{1,-1\}$.
  \end{itemize}
We generalize $G(f,g)$, constructing a graph $W(f,g)$ on the set
$C'(f,g)$, placing an edge between distinct elements $(k,s,l,t)$
and $(k',s',l',t')$ whenever the composites $k^s\circ l^t$ and
$k'^{s'}\circ l'^{t'}$ are equal in $\G$.
Thus, the connected components of $W(f,g)$
are the fibers of the evaluation map $\ev:C'(f,g)\to\Mor(\G)$
sending $(k,s,l,t)$ to $k^sl^t$.
We will always draw the vertices of $W(f,g)$ in the following configuration
  \[  \tikz{
    \path node (mid) {}
          node [above = 1.25cm of mid] {$(f,1,g,1)$}
          node [above = .5cm of mid] {$(g,1,f,1)$}
          node [left  = 3.75cm   of mid] {$(f,-1,g,1)$}
          node [left  = 1.25cm   of mid] {$(g,1,f,-1)$}
          node [right = 3.75cm   of mid] {$(f,1,g,-1)$}
          node [right = 1.25cm   of mid] {$(g,-1,f,1)$}
          node [below = 1.25cm of mid] {$(f,-1,g,-1)$}
          node [below = .5cm of mid] {$(g,-1,f,-1)$};
  }\]
and will usually suppress the labels on the vertices to points $*$.
Note that $G(f,g)$ is the full subgraph of $W(f,g)$
defined on the upper six vertices.
Similarly, $G(f^{-1},g)$, $G(f,g^{-1})$, and
$G(f^{-1},g^{-1})$ are (respectively) isomorphic to the
left, right, and lower T-shaped subgraphs of $W(f,g)$.
For example, the canonical inclusion
  \[G(f^{-1},g)\mono W(f,g)\]
is defined by sending $(f^{-1},s,g,t)$ to $(f,-s,g,t)$
and $(g,s,f^{-1},t)$ to $(g,s,f,-t)$.

The graph $W(f,g)$ is determined by its subgraphs $G(f^{\pm1},g^{\pm1})$,
which are in turn determined by $\Sd\,\G$.
By considering these four graphs simultaneously, one may
obtain a good deal of information concerning $W(f,g)$.

As with $G(f,g)$, any two points of $C'(f,g)$ belonging to
the same connected component of $W(f,g)$ must have an edge between them.
Configurations in $W(f,g)$ correspond to equational conditions.
We have the tables below, where the following are equivalent:
  \begin{enumerate}
    \item the upper equation holds
    \item the lower graph is a subgraph of $W(f,g)$
    \item a single edge of the lower graph is in $W(f,g)$.
  \end{enumerate}
\begin{dontbotheriftheequationisoverlong}\label{table:eightgraphs}
\begin{minipage}{2\textwidth}\centering
\begin{tabular}{|c|c|c|c|c|c|}\hline
Equation
& $fg=gf$
& $f^2=\id$
& $g^2=\id$
& $f^2=g^2$
& $f^2=g^{-2}$
\\\hline Subgraph
  & \eightgraph\eightfggf\endeightgraph
  & \eightgraph\eightffid\endeightgraph
  & \eightgraph\eightggid\endeightgraph
  & \eightgraph\eightffgg\endeightgraph
  & \eightgraph\eightffgginv\endeightgraph
\\\hline
\\[-2mm]\hline
  Equation
&
& $fgf=g$
& $gfg=f$
& $fg^{-1}f=g$
& $fgf=g^{-1}$
\\\hline
    Subgraph
  &
  & \eightgraph\eightfgfg\endeightgraph
  & \eightgraph\eightgfgf\endeightgraph
  & \eightgraph\eightfginvfg\endeightgraph
  & \eightgraph\eightfgfginv\endeightgraph
\\\hline
\end{tabular}
\end{minipage}
\end{dontbotheriftheequationisoverlong}

The above graphs have the following symmetries:
  \begin{itemize}
    \item There is an edge in $W(f,g)$ between $(k,s,l,t)$ and $(m,u,n,v)$
          if and only if there is an edge between $(l,-t,k,-s)$
          and $(n,-v,m,-u)$.
          This comes from inversion
            \[k^sl^t=m^un^v \iff l^{-t}k^{-s}=n^{-v}m^{-u}.\]
          The pairing between vertices $(k,s,l,t)$
          and $(l,-t,k,-s)$ is illustrated below:
            \[\eightgraph\eightfgfginv\eightfginvfg\endeightgraph.\]
    \item For each graph above, there is another graph
          obtained by interchanging the roles of $f$ and $g$.
          For example, the graphs corresponding to the equations
          $f^2=\id$ and $g^2=\id$ are one such pair.
  \end{itemize}

{\em\ul{Case 4, subcase 1:
Suppose $f^2\neq g^2$ and $fg=gf$ and $f^2\neq\id$ and $g^2\neq\id$.}}\hfill\\
Because $f^2\neq g^2$ and $f^2\neq\id$ and $g^2\neq\id$,
the graph $W(f,g)$ cannot contain any of the edges below:
  \[\eightgraph\eightffgg\eightffid\eightggid\endeightgraph.\]
The graph $W(f,g)$ must be equal to one of the two graphs below.
  \[\setlength{\tabcolsep}{.5cm}\begin{tabular}{cc}
      \eightgraph\eightfggf\endeightgraph
  &   \eightgraph\eightfggf\eightffgginv
        \draw[shorten >=-1.25mm, shorten <=-1.25mm](a) to [bend right=20] (h);
        \draw[shorten >=-1.25mm, shorten <=-1.25mm](b) to [bend left=20] (g);
      \endeightgraph
  \\  \text{if $f^2\neq g^{-2}$}
  &   \text{if $f^2=g^{-2}$}
    \end{tabular}\]
The graph on the left is guaranteed by the condition $fg=gf$.
The graph on the right is the only possible extension
that is consistent with the assumptions $f^2\neq g^2$
and $f^2\neq\id$ and $g^2\neq\id$.

Suppose that $h=fg=gf$.
Then in each of the two graphs above, the subgraph $G(f,g)$
contains two elements of the fiber of $h$,
namely $(f,1,g,1)$ and $(g,1,f,1)$.
Thus we have $\exists_2\bigtri fgh$.
Similarly, the left and right T-shaped subgraphs, which are
respectively isomorphic to $G(f^{-1},g)$ and $G(f,g^{-1})$,
each contain either two or four elements of the fiber of $h$
(depending on whether $f^2=g^{-2}$).
Therefore we have
  \[h=fg=gf\implies\exists_2\bigtri fgh
    ~\text{ and }~\exists_{\geq2}\bigtri{f^{-1}}gh
    ~\text{ and }~\exists_{\geq2}\bigtri f{g^{-1}}h
  .\]

On the other hand, suppose that $h$ is not equal to $fg$ or $gf$,
and that the fiber of $h$ is non-empty.
If $f^2=g^{-2}$ then $\ev^{-1}\{h\}$ must be one of the connected components
  \[\Bigl\{(f,-1,g,1),(g,1,f,-1)\Bigr\}~\text{ or }~\Bigl\{(f,1,g,-1),(g,-1,f,1)\Bigr\},\]
in which case we have (respectively) $\nexists\bigtri f{g^{-1}}h$ or
$\nexists\shifttext{3mm}{\bigtri {f^{-1}}gh}$.
If $f^2\neq g^{-2}$ then the fiber of $h$ must be one of the
connected components not equal to $\{(f,1,g,1),(g,1,f,1)\}$,
so we have $\nexists\bigtri f{g^{-1}}h$ or
$\nexists\shifttext{2mm}{\bigtri {f^{-1}}gh}$ or $\nexists\bigtri fgh$.
    \[\text{We conclude that }~
    h=fg=gf \iff
    \shifttext{3mm}{$\exists\bigtri fgh$}
    \text{ and }\exists\shifttext{4mm}{$\bigtri{f^{-1}}gh$}
    \text{ and }\exists\bigtri f{g^{-1}}h
    .\]

{\em\ul{Case 4, subcase 2: $f^2\neq g^2$ and
$fg\neq gf$ and $f^2\neq\id$ and $g^2\neq\id$ and $f^2=g^{-2}$.}}\hfill\\
Because $f^2\neq g^2$ and $fg\neq gf$ and $f^2\neq\id$ and $g^2\neq\id$,
the graph $W(f,g)$ cannot contain any of the edges below:
  \[\eightgraph\eightffgg\eightfggf\eightffid\eightggid\endeightgraph.\]
Therefore, because $f^2=g^{-2}$, the graph $W(f,g)$ must be equal to
one of the two graphs displayed below.
  \[\setlength{\tabcolsep}{.5cm}\begin{tabular}{cc}
      \eightgraph\eightffgginv\eightfginvfg\endeightgraph
  &   \eightgraph\eightffgginv\endeightgraph
  \\  \text{if $fg^{-1}f=g$}
  &   \text{if $fg^{-1}f\neq g$}
    \end{tabular}\]

Suppose first that $h=fg$ or $h=gf$.
In each case above, the fiber of $h$ must
be equal to one of the connected components
  \[\Bigl\{(f,1,g,1),(f,-1,g,-1)\Bigr\}~\text{ or }~\Bigl\{(g,1,f,1),(g,-1,f,-1)\Bigr\}.\]
These connected components are both contained in the intersection
of the left and right T-shaped subgraphs.
Those two subgraphs correspond (respectively) to $G(f^{-1},g)$ and
$G(f,g^{-1})$, so we have
$\exists_2\bigtri{f^{-1}}gh$ and $\exists_2\bigtri f{g^{-1}}h$.

Suppose now that $h$ is not equal to $fg$ or $gf$, and that
the fiber $\ev^{-1}\{h\}$ is non-empty.
If $fg^{-1}f=g$ then $\ev^{-1}\{h\}$ must be equal to one of the connected components
  \[\Bigl\{(f,-1,g,1),(g,-1,f,1)\Bigr\}~\text{ or }~\Bigl\{(g,1,f,-1),(f,1,g,-1)\Bigr\}.\]
If $fg^{-1}f\neq g$ then the fiber of $h$ is equal to one of the singleton sets
  \[\bigl\{(f,-1,g,1)\bigr\},~\bigl\{(g,-1,f,1)\bigr\},~\bigl\{(g,1,f,-1)\bigr\},~\bigl\{(f,1,g,-1)\bigr\}.\]
In either case, the intersection of $\ev^{-1}\{h\}$ with the left T-shaped
subgraph contains at most one point, and similarly for the intersection
with the right T subgraph.
Therefore we have the following:
  \[
    \text{$h=fg$ or $h=gf$}
      \iff
    \text{$\exists_2\bigtri{f^{-1}}gh$ and $\exists_2\bigtri f{g^{-1}}h$}.
  \]

{\em\ul{Case 4, subcase 3: $f^2\neq g^2$ and $fg\neq gf$ and
$f^2\neq\id$ and $g^2\neq\id$ and $f^2\neq g^{-2}$.}}\hfill\\
This is the most complicated case.
As before, our strategy is to calculate all possibilities for the graph
$W(f,g)$ that are consistent with the above assumptions on $f$ and $g$.

First note that the nine graphs described in the table 
on page \pageref{table:eightgraphs} have pairwise disjoint edges sets,
and that there are a total of 28 edges among those graphs.
Therefore, the complete graph on 8 vertices has edges given by the
disjoint union of the edge sets corresponding to the nine graphs displayed
in the table (because the complete graph on 8 vertices has 28 edges).

Next, recall that if $W(f,g)$ contains a single edge
from one of the nine graphs in the table, then $W(f,g)$
contains that entire graph as a subgraph.
Therefore, $W(f,g)$ is the union of some choice of subgraphs
from among the nine displayed in the table (where by ``union'' we mean
that the edges of $W(f,g)$ are obtained as a disjoint union of the
edge sets of some subgraphs displayed in the table).

\newcommand\T{\textnormal T}
\newcommand\F{\textnormal F}
A given graph from the table is contained in $W(f,g)$
if and only if the corresponding equation is satisfied.
Thus, the graph $W(f,g)$ determines a set of ``generating relations''
among the nine equations listed in the table.
Explicitly, $W(f,g)$ assigns a value of $\T/\F$ to each of
of the nine equations.
Given our assumptions that $f^2\neq g^2$ and $fg\neq gf$ and $f^2\neq\id$
and $g^2\neq\id$ and $f^2\neq g^{-2}$,
the graph $W(f,g)$ cannot contain any of the edges below:\\[-25pt]
  \[\begin{tikzpicture}
    \newdimen\clen
    \clen=1.2cm
    \newdimen\dlen
    \dlen=.75cm
          \node (a) at ( 0, \clen) {$*$};
          \node (b) at ( 0, \dlen) {$*$};
          \node (c) at (-\clen, 0) {$*$};
          \node (d) at (-\dlen, 0) {$*$};
          \node (e) at ( \clen, 0) {$*$};
          \node (f) at ( \dlen, 0) {$*$};
          \node (g) at ( 0,-\clen) {$*$};
          \node (h) at ( 0,-\dlen) {$*$};
      \eightffid\eightggid
      \draw [shorten >=-.75mm, shorten <=-.75mm] (a) -- (b);
      \draw [shorten >=-.75mm, shorten <=-.75mm] (c) -- (d);
      \draw [shorten >=-.75mm, shorten <=-.75mm] (e) -- (f);
      \draw [shorten >=-.75mm, shorten <=-.75mm] (g) -- (h);
      \draw[shorten >=-1.25mm, shorten <=-1.25mm](a) to [bend left=25] (g);
      \draw[shorten >=-1.00mm, shorten <=-1.00mm](b) -- (h);
      \draw[shorten >=-1.25mm, shorten <=-1.25mm](c) to [bend left=25] (e);
      \draw[shorten >=-1.00mm, shorten <=-1.00mm](d) -- (f);
    \end{tikzpicture}.\]\\[-15pt]
We need only to consider the remaining four equations
from the table on page \pageref{table:eightgraphs}:
  \begin{equation}\label{eqns:fourGenerators}
    fgf=g       
    ,\hspace{5mm}
    gfg=f       
    ,\hspace{5mm}
    fg^{-1}f=g  
    ,\hspace{5mm}
    fgf=g^{-1}  .
  \end{equation}
Below is the union of the graphs corresponding to the above four equations:
  \begin{equation}\label{supergraph}
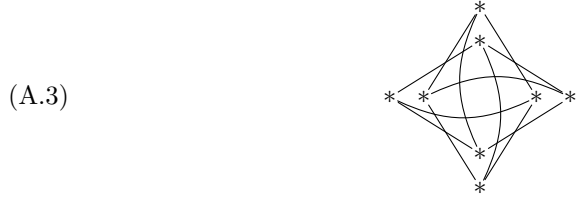
  \begin{minipage}[m]{.5\textwidth}\centering
  \begin{tikzpicture}
    \newdimen\clen
    \clen=1.2cm
    \newdimen\dlen
    \dlen=.75cm
        \node (a) at ( 0, \clen) {$*$};
        \node (b) at ( 0, \dlen) {$*$};
        \node (c) at (-\clen, 0) {$*$};
        \node (d) at (-\dlen, 0) {$*$};
        \node (e) at ( \clen, 0) {$*$};
        \node (f) at ( \dlen, 0) {$*$};
        \node (g) at ( 0,-\clen) {$*$};
        \node (h) at ( 0,-\dlen) {$*$};
        \eightfgfg\eightgfgf
        \draw[shorten >=-1.25mm, shorten <=-1.25mm](c) to [bend right=25] (f);
        \draw[shorten >=-1.25mm, shorten <=-1.25mm](d) to [bend left=25] (e);
        \draw[shorten >=-1.25mm, shorten <=-1.25mm] (a) to [bend right=25] (h);
        \draw[shorten >=-1.25mm, shorten <=-1.25mm](b) to [bend left=25] (g);
    \end{tikzpicture}.\end{minipage}\end{equation}
The graph $W(f,g)$ must be a subgraph of the union above, for any edge in
$W(f,g)$ belongs to a graph corresponding to one of the generating
relations (\ref{eqns:fourGenerators}).

Each assignment of $\T/\F$ to the four
generating relations (\ref{eqns:fourGenerators})
corresponds to a subgraph of the above (\ref{supergraph}).
Given that each connected, full subgraph of $W(f,g)$ must be a complete graph,
only some assignments of $\T/\F$
correspond to actual possibilities for the graph $W(f,g)$.
Thus, only some subgraphs of the union above could possibly
be equal to $W(f,g)$.

We consider the space $S=\{\T,\F\}^4$ of quadruples of $\T$'s
and $\F$'s, where the four coordinates correspond respectively
to each of the equations (\ref{eqns:fourGenerators}).
Given an element $\gamma$ of $S$, the corresponding subgraph of
(\ref{supergraph}) is given by the union of
graphs corresponding to those equations (\ref{eqns:fourGenerators})
for which $\gamma$'s value is $\T$.
For example, the assignment $(\F,\F,\F,\F)$ corresponds
to the graph with eight vertices and no edges,
the assignment $(\T,\T,\T,\T)$ corresponds to the
union (\ref{supergraph}), and the assignment
$(\T,\F,\F,\F)$ gives the graph associated
to the equation $fgf=g$ (see the table on page \pageref{table:eightgraphs}).

We will say that a graph is \textit{valid} if
each of its connected components is a complete graph.
The graph $W(f,g)$ must be a subgraph of (\ref{supergraph}) that is valid,
and it and must correspond to one of the elements of $S$.
We claim that the following six elements of $S$ are the only ones
whose corresponding graphs are valid:
  \begin{equation}\label{ScorrespondingToValid}
  \begin{tabular}[b]{c|c|c}
    (\F,\F,\F,\F)
  & (\T,\F,\F,\F)   
  & (\F,\T,\F,\F)   \\\hline
    (\F,\F,\T,\F)   
  & (\F,\F,\F,\T)   
  & (\F,\F,\T,\T)   
  \end{tabular}.\end{equation}
The elements of $S$ above correspond
(respectively) to the following six graphs.
    \begin{equation}\label{grid:sixeightgraphs}
  \setlength{\tabcolsep}{.5cm}\begin{tabular}[m]{c|c|c}
    \eightgraph\endeightgraph
  & \eightgraph\eightfgfg\endeightgraph
  & \eightgraph\eightgfgf\endeightgraph
  \\
  & \text{$fgf=g$}
  & \text{$gfg=f$}
  \\\hline
    \eightgraph\eightfginvfg\endeightgraph
  & \eightgraph\eightfgfginv\endeightgraph
  & \eightgraph\eightfginvfg\eightfgfginv\endeightgraph
  \\
    \text{$fg^{-1}f=g$}
  & \text{$fgf=g^{-1}$}
  & \text{and $fgf=g^{-1}$}
  \\
  &&\text{$fg^{-1}f=g$}
  \end{tabular}\end{equation}
It is clear that any connected component in any of the above graphs
is complete: each connected component has size at most two.

It remains to show that the elements of $S$ that are \textit{not}
displayed in Table (\ref{ScorrespondingToValid}) all correspond to
graphs that are non-valid.
We have seen that the graph (\ref{supergraph}), which corresponds to the
assignment $(\T,\T,\T,\T)$, is not valid.
The remaining elements of $S$ are displayed below, together with
their corresponding subgraphs.
  \[\begin{tabular}{c|c|c}
    \eightgraph\eightfgfg\eightfginvfg\endeightgraph
  & \eightgraph\eightgfgf\eightfginvfgtwo\endeightgraph
  & \eightgraph\eightfgfg\eightgfgf\eightfginvfg\endeightgraph
    \\
    (\T,\F,\T,\F) 
  & (\F,\T,\T,\F) 
  & (\T,\T,\T,\F) 
    \\\hline
    \eightgraph\eightfgfg\eightfgfginvtwo\endeightgraph
  & \eightgraph\eightgfgf\eightfgfginv\endeightgraph
  & \eightgraph\eightfgfg\eightgfgf\eightfgfginv\endeightgraph
    \\
    (\T,\F,\F,\T) 
  & (\F,\T,\F,\T) 
  & (\T,\T,\F,\T) 
    \\\hline
    \eightgraph\eightfgfg\eightfginvfg\eightfgfginvtwo\endeightgraph
  & \eightgraph\eightgfgf\eightfginvfgtwo\eightfgfginv\endeightgraph
  & \eightgraph\eightfgfg\eightgfgf\endeightgraph
    \\
    (\T,\F,\T,\T) 
  & (\F,\T,\T,\T) 
  & (\T,\T,\F,\F) 
  \end{tabular}\]
By inspection, each of the above graphs has a connected component
that is not complete.
Thus, none of the graphs displayed above are valid,
so $W(f,g)$ must be equal to one
of the six graphs (\ref{grid:sixeightgraphs}).

We can now claim that $h$ equals $fg$ or $gf$
if and only if one of the following is satisfied:
      \begin{enumerate}
        \doublespacing
        \item $\exists_1\bigtri fgh$ and
              $\exists_1\shifttext{1mm}{\bigtri{f^{-1}}gh}$ and
              $\exists_1\bigtri f{g^{-1}}h$, or
        \item $\exists_1\shifttext{1mm}{\bigtri{f^{-1}}{g^{-1}}h}$ and
              two of the following three hold:
                \[\text{$\exists_2\bigtri fgh$ or
                        $\exists_2\shifttext{1mm}{\bigtri {f^{-1}}gh}$ or
                        $\exists_2\bigtri f{g^{-1}}h$.}\]
      \end{enumerate}
The proof of this claim is given below in two steps: ``only if'' and ``if.''

\subsection*{Only if} Suppose that $h$ is equal to $fg$ or to $gf$.
If $W(f,g)$ is equal to one of the graphs on the left hand side
of the table (\ref{grid:sixeightgraphs}), corresponding
to assignment $(\F,\F,\F,\F)$ or $(\F,\F,\T,\F)$,
then condition (1) is satisfied:
each of the subgraphs $G(f,g)$, $G(f^{-1},g)$ and $G(f,g^{-1})$ of $W(f,g)$
contains exactly one element of the fiber $\ev^{-1}\{h\}$.

If $W(f,g)$ is equal to one of the graphs in the center or on
the right hand side of the table, corresponding to
$(\T,\F,\F,\F)$ or $(\F,\T,\F,\F)$ or $(\F,\F,\F,\T)$ or $(\F,\F,\T,\T)$,
then condition (2) is satisfied:
the subgraph $G(f^{-1},g^{-1})$ of $W(f,g)$ contains one element of the fiber of $h$,
and exactly two of the subgraphs $G(f,g)$, $G(f^{-1},g)$ and $G(f,g^{-1})$
contain two elements of the fiber of $h$.

\subsection*{If}
If $W(f,g)$ is equal to one of the two graphs on the left of the table
(\ref{grid:sixeightgraphs}), and if condition (1)
is satisfied, then $h=fg$ or $h=gf$.
If $W(f,g)$ is equal to one of the other four graphs,
then (1) is false for every $h$.

If $W(f,g)$ is equal to one of the four graphs on the right or in the middle of
the table, and if condition (2) is satisfied, then we must have $h=fg$ or $h=gf$.
If $W(f,g)$ is equal to one of the two graphs on the left of the table,
then condition (2) is false for every $h$.
\qed

\par}

\end{document}